\documentclass[11pt, oneside]{article}   	
\usepackage{geometry}                		
\usepackage{graphicx}				
\usepackage{color}
\usepackage{amssymb}
\usepackage{amsmath}
\usepackage{amsthm}
\usepackage{bm}
\usepackage{aliases}
\usepackage[utf8]{inputenc}
\usepackage[english]{babel}
\usepackage{lineno}
\usepackage{listings}
\usepackage{textcomp}
\usepackage{cleveref}
\usepackage{stmaryrd}
\usepackage{siunitx}
\usepackage{commath}

\title{Globally constraint-preserving FR/DG scheme for Maxwell's equations at all orders}
\author{
Arijit Hazra\footnote{TIFR Center for Applicable Mathematics, Bangalore, India. Email: \texttt{arijit@tifrbng.res.in}}, \ 
Praveen Chandrashekar\footnote{TIFR Center for Applicable Mathematics, Bangalore, India. Email: \texttt{praveen@tifrbng.res.in}}, \ and \ 
Dinshaw S. Balsara\footnote{Dept. of Physics, Univ. of Notre Dame, USA. Email: \texttt{dbalsara@nd.edu}}}
\date{}							

\begin{document}
\maketitle

\lstset{
   language=bash,
   keywordstyle=\bfseries\ttfamily\color[rgb]{0,0,1},
   identifierstyle=\ttfamily,
   commentstyle=\color[rgb]{0.133,0.545,0.133},
   stringstyle=\ttfamily\color[rgb]{0.627,0.126,0.941},
   showstringspaces=false,
   basicstyle=\small,
   numberstyle=\footnotesize,
   numbers=none,
   stepnumber=1,
   numbersep=10pt,
   tabsize=2,
   breaklines=true,
   prebreak = \raisebox{0ex}[0ex][0ex]{\ensuremath{\hookleftarrow}},
   breakatwhitespace=false,
   aboveskip={0.1\baselineskip},
    columns=fixed,
    upquote=true,
    extendedchars=true,
    backgroundcolor=\color[rgb]{0.9,0.9,0.9}
}

\begin{abstract}
Computational electrodynamics (CED), the numerical solution of Maxwell's equations, plays an incredibly important role in several problems in science and engineering. High accuracy solutions are desired, and the discontinuous Galekin (DG) method is one of the better ways of delivering high accuracy in numerical CED. Maxwell's equations have a pair of involution constraints for which mimetic schemes that globally satisfy the constraints at a discrete level are highly desirable. Balsara and K\"appeli (2018) presented a von Neumann stability analysis of globally constraint-preserving DG schemes for CED  up to fourth order. That paper was focused on developing the theory and documenting the superior dissipation and dispersion of DGTD schemes in media with constant permittivity and permeability. In this paper we present working DGTD schemes for CED that go up to fifth order of accuracy and analyze their performance when permittivity and permeability vary strongly in space.

	Our DGTD schemes achieve constraint preservation by collocating the electric displacement and magnetic induction as well as their higher order modes in the faces of the mesh.  Our first finding is that at fourth and higher orders of accuracy, one has to evolve some zone-centered modes in addition to the face-centered modes. It is well-known that the limiting step in DG schemes causes a reduction of the optimal accuracy of the scheme; though the schemes still retain their formal order of accuracy with WENO-type limiters. In this paper we document simulations where permittivity and permeability vary by almost an order of magnitude without requiring any limiting of the DG scheme. This very favorable second finding ensures that DGTD schemes retain optimal accuracy even in the presence of large spatial variations in permittivity and permeability. We also study the conservation of electromagnetic energy in these problems. Our third finding shows that the electromagnetic energy is conserved very well even when permittivity and permeability vary strongly in space; as long as the conductivity is zero. 
\end{abstract}
\section{Introduction}
The numerical solution of Maxwell's equations plays an incredibly important role in the computational solution of many problems in science and engineering. The Finite Difference Time Domain (FDTD) method, originally proposed by Yee~\cite{Yee1966} and greatly developed since the seminal papers by Taflove and Brodwin~\cite{Taflove1975}, \cite{Taflove1975a}, has been the mainstay for computational electrodynamics (CED) chiefly because it globally preserves the involution constraints that are inherent in Maxwell's equations. Several modern texts and reviews document the development of FDTD (Taflove and Hagness~\cite{Taflove2016}, Taflove, Oskooi and Johnson 2013, Gedney~\cite{Gedney2011}). As a result, given the enormously well-cited works by Yee and Taflove, it is very desirable to include global preservation of involution constraints into more modern schemes for CED. 

	FDTD had been formulated well before the explosion in modern higher order Godunov schemes, which started with the pioneering work of van Leer~\cite{VanLeer1977},~\cite{vanLeer1979}. Since that work, there has been a strong drive to include the physics of wave propagation into the numerical solution of hyperbolic systems; CED being a case in point. Two strong strains of early effort to develop higher order Godunov schemes for CED include finite-volume time-domain (FVTD) methods (Munz et al.~\cite{Munz2000}, Ismagilov~\cite{Ismagilov2015}, Barbas and Velarde~\cite{Barbas2015}) and discontinuous Galerkin time-domain (DGTD) methods (Cockburn \& Shu~\cite{Cockburn1989a},~\cite{Cockburn1998};  Cockburn et al.~\cite{Cockburn1990}, Hesthaven and Warburton~\cite{Hesthaven2002a}). There has also been a strong effort in the engineering CED community to design DGTD schemes for CED (Chen \& Liu~\cite{Chen2009}, Ren et al.~\cite{Ren2017}, Angulo et al. 2015) and some of those methods indeed use locally constraint-preserving bases. However, none of those DGTD methods incorporated the very beneficial globally constraint-preserving aspect of FDTD. For that next phase of evolution, one had to wait for developments that emerged in the field of numerical MHD and are now rapidly finding their way into CED.
	
	In MHD, one evolves Faraday's law in addition to the equations of computational fluid dynamics (Brecht et al.~\cite{Brecht1981}, Evans and Hawley~\cite{Evans1988}, DeVore~\cite{DeVore1991}, Dai and Woodward~\cite{Dai1998}, Ryu et al.~\cite{Ryu1998}, Balsara and Spicer~\cite{Balsara1999}). While studying adaptive mesh refinement and numerical schemes for MHD, advances were made in constraint-preserving reconstruction of magnetic fields (Balsara~\cite{Balsara2001},~\cite{Balsara2004},~\cite{Balsara2009}, Balsara and Dumbser~\cite{Balsara2015}, Xu et al.~\cite{Xu2016}). This made it possible to start with the face-centered magnetic induction fields in the Yee-type mesh and specify them at all locations within a computational zone. The edge-centered electric fields that are inherent to a constraint-preserving update of Faraday's law would then have to be multidimensionally upwinded. This multidimensional upwinding was achieved by using a newly-designed multidimensional Riemann solver (Balsara~\cite{Balsara2010}, \cite{Balsara2012a}, \cite{Balsara2014a}, \cite{Balsara2015b}, Balsara, Dumbser and Abgrall~\cite{Balsara2014}, Balsara and Dumbser~\cite{Balsara2015a}, Balsara and Nkonga~\cite{Balsara2017b}). These twin innovations, consisting of constraint-preserving reconstruction of vector fields, and multidimensional Riemann solvers, permitted a logically complete description of numerical MHD. Along the way, a third innovation in ADER (Arbitrary DERivatives in space and time) timestepping schemes was added which, while not essential, greatly simplified the accurate temporal evolution of MHD variables (Dumbser et al.~\cite{Dumbser2008},~\cite{Dumbser2013}, Balsara et al.~\cite{Balsara2009a},~\cite{Balsara2013}). DG schemes for the induction equation that were globally constraint-preserving were also devised by Balsara and K\"appeli~\cite{Balsara2017a}. The stage was now set for migrating these innovations back again to CED.
	
	FVTD schemes for CED that were based on the above-mentioned three innovations were developed in Balsara et al. (\cite{Balsara2016}, \cite{Balsara2017}, \cite{Balsara2018}). The constraint-preservation was accomplished by making a constraint-preserving reconstruction of the magnetic induction and the electric displacement. Unlike FDTD that operates on a pair of staggered control volumes, the present methods operate on the same control volume, see Fig. 1 from Balsara et al.~\cite{Balsara2017}. To ensure the mimetic preservation of the constraints, the primal variables were taken to be the facially collocated normal components of the electric displacement and the magnetic induction; where both vector fields were collocated on the same faces. The facially collocated magnetic induction evolves in response to the edge collocated electric fields, yielding a discrete representation of Faraday's law. The facially collocated electric displacement evolves in response to the edge collocated magnetic fields, yielding a discrete representation of Ampere's law. The resulting methods were indeed globally constraint-preserving in that the magnetic induction always remains divergence-free at all locations on the mesh and the electric displacement always satisfies the constraint imposed by Gauss' Law at any location on the mesh. All these above-mentioned advances that were made to mimetic FVTD schemes have been recently embedded into the DGTD schemes that we report on below.
	
	Globally constraint-preserving DGTD schemes were also explored in Balsara and K\"apelli (2018, BK henceforth). BK found that higher order DGTD schemes were almost totally free of dispersion error, having a dispersion error that was almost ~75 times smaller than FDTD. Since the methods were based on Riemann solvers, some dissipation is inevitable. Even so, BK found that the higher order DGTD schemes were almost free of dissipation even when electromagnetic waves spanned only a few zones. The first paper by BK was strongly focused on von Neumann stability analysis of DGTD schemes for CED. The present paper extends this study in several ways, which we list in the rest of this paragraph. First, BK focused on DGTD schemes up to fourth order, whereas the present paper presents flux reconstruction (FR) based DGTD schemes\footnote{We will use DGTD to refer to the current scheme though strictly it is a combination of FR and DG schemes.} that go up to fifth order, and the reconstruction schemes at fourth and fifth order presented in this work are new. This drive to higher order enables us to generalize an observation from BK who found that from fourth order and upwards one has to include some volumetrically-evolved modes in addition to the facially-evolved modes. Second, BK did not touch on the topic of limiting DGTD schemes, specifically because they realized that the non-linear hybridization of DGTD schemes for CED should be done with the utmost carefulness, if even it is needed. The limiting of any DG scheme always diminishes the optimal accuracy of a DG scheme. This reduction occurs to a greater or lesser extent based on whether the limiter is applied more or less aggressively, respectively. BK did find that DGTD schemes did not need any limiting up to fourth order but they only tried situations where the permittivity and permeability were constant. In this paper we report the very favorable finding that DGTD schemes don't seem to require any limiting even when the permittivity and permeability vary by almost an order of magnitude. This result is quite useful; however, it is predicated on the assumption that the conductivity is zero. Dealing with non-zero conductivity will be the topic of a subsequent study. Third, BK realized that when conductivity is zero the Maxwell's equations conserve a quadratic energy. While this energy conservation was not built into the scheme, BK found the very desirable result that quadratic energy was effectively conserved by fourth order DGTD schemes when the waves spanned only a few zones. This paper extends this study to fifth order DGTD schemes. Furthermore, we study the energy conservation of high order DGTD when the permittivity and permeability vary with space. A fourth offering in this paper is a proof that DGTD schemes for CED that are constructed according to the principles in BK and this paper are indeed energy stable.
	
	The rest of the paper is organized as follows. Section~(\ref{sec:maxwell}) introduces the Maxwell's equations and the simplified 2-D model that we consider in this paper. Section~(\ref{sec:approx}) explains the polynomial spaces used to approximate the solution variables and section~(\ref{sec:rec3}) explains the divergence-free reconstruction scheme at fourth order of accuracy. The numerical descritization of the Maxwell's equations using flux reconstruction and DG method is shown in section~(\ref{sec:scheme}) together with constraint preserving property. In section~(\ref{sec:stab}), we perform the stability analysis of the semi-discrete scheme at first order and show the dissipative character coming from the 1-D and 2-D Riemann solvers. Section~(\ref{sec:res}) presents many test cases to demonstrate the performance of the scheme. Finally, the Appendix explains the reconstruction scheme at other orders and also briefly discusses the Riemann solvers.
\section{Maxwell's equations}
\label{sec:maxwell}
The Maxwell's equations are a system of linear partial differential equations that model the wave propagation behaviour of electric and magnetic fields in free space and material media. They can be written in vector form as
\[
\df{\bm{B}}{t} + \nabla \times \bm{E} = 0, \qquad \df{\bm{D}}{t} - \nabla       \times \bm{H} = -\bm{J}
\]
where
\begin{center}
\begin{tabular}{ll}
$\bm{B}$ = magnetic flux density & $\bm{D}$ = electric flux density \\
$\bm{E}$ = electric field & $\bm{H}$ = magnetic field \\
& $\bm{J}$ = electric current density
\end{tabular}
\end{center}
The fields are related to one another by constitutive laws
\[
\bm{B} = \mu \bm{H}, \qquad \bm{D} = \varepsilon \bm{E}, \qquad \bm{J} = \sigma \bm{E} \qquad \mu, \varepsilon \in \re^{3 \times 3} \textrm{ symmetric}
\]
where the coefficients
\[
\begin{aligned}
\varepsilon &= \textrm{permittivity tensor} \\
\mu &= \textrm{magnetic permeability tensor} \\
\sigma &= \textrm{conductivity}
\end{aligned}
\]
are material properties and are in general tensorial functions of spatial coordinates. In free space $\varepsilon = \varepsilon_0 = 8.85 \times 10^{-12}$ F/m and $\mu = \mu_0 = 4\pi \times 10^{-7}$ . The divergence of the electric flux gives the electric charge density $\rho$, which itself obeys a conservation law
\[
\nabla \cdot \bm{D} = \rho, \qquad \df{\rho}{t} + \nabla \cdot \bm{J} = 0
\]
Moreover, since magnetic monopoles have never been observed in nature, the divergence of the magetic flux must be zero $\nabla\cdot\B = 0$, which is an additional constraint that must be satisfied by the solution. Note that if the initial condition is divergence-free, then under the time evolution induced by the Maxwell's equations, the divergence remains zero at future times also.

In the present work, we will consider a 2-D model of the Maxwell's equations (TE polarization) for which the equations can be written in Cartesian coordinates as
\begin{eqnarray}
\df{D_x}{t} - \df{H_z}{y} &=& 0 \\
\df{D_y}{t} + \df{H_z}{x} &=& 0 \\
\df{B_z}{t} + \df{E_y}{x} - \df{E_x}{y} &=& 0
\end{eqnarray}
where
\[
(E_x, E_y) = \frac{1}{\varepsilon} (D_x,D_y), \qquad H_z = \frac{1}{\mu} B_z
\]
and $\mu$, $\varepsilon$ are scalars which may depend on spatial coordinates. We will consider a constraint on the divergence of the electric flux density $\D$ instead of the magnetic field, which has only one component in the above model. The above system of three PDE can be written in conservation form
\begin{equation}
\df{\con}{t} + \df{\fx}{x} + \df{\fy}{y} = 0
\label{eq:claw}
\end{equation}
where
\[
\con = \begin{bmatrix}
D_x \\
D_y \\
B_z \end{bmatrix}, \qquad
\fx = \begin{bmatrix}
0 \\
H_z \\
E_y \end{bmatrix} = \begin{bmatrix}
0 \\
\frac{1}{\mu} B_z \\
\frac{1}{\varepsilon} D_y \end{bmatrix}, \qquad
\fy = \begin{bmatrix}
-H_z \\
0 \\
-E_x \end{bmatrix} = \begin{bmatrix}
-\frac{1}{\mu} B_z \\
0 \\
-\frac{1}{\varepsilon} D_x \end{bmatrix}
\]
This is a system of hyperbolic conservation laws for which a Riemann problem can be solved exactly to determine the fluxes required in the numerical schemes like finite volume and DG method. This is explained in the Appendix. While for simplicity, we consider scalar material properties, the constraint preserving nature of the scheme holds for general material properties. In fact, everything we describe in this paper holds for the general case, and the only additional change required is to use the Riemann solvers for the general case, which are explained e.g., in \cite{Balsara2017}.

In the absence of currents, the Maxwell's equations conserve the total energy provided there is no net gain of energy at the boundaries or if we have periodic boundaries. For the 2-D model that we consider in this work, we can first show that the following additional conservation law holds
\[
\df{}{t} \left[ \frac{1}{2\varepsilon}(D_x^2 + D_y^2) + \frac{1}{2\mu} B_z^2 \right] + \df{}{x}(H_z E_y) - \df{}{y}(H_z E_x) = 0
\]
which implies that the quantity
\[
\tote(t) = \int_\Omega \left[ \frac{1}{2\varepsilon}(D_x^2 + D_y^2) + \frac{1}{2\mu} B_z^2 \right] \ud x \ud y
\]
which is the total energy, is conserved under periodic boundary conditions or if the net flux at the boundaries is zero.

\section{Approximation spaces}
\label{sec:approx}
We would like to approximate the vector field $\D$ such that its divergence is zero inside the cell if there is no charge density. If there is some electric charge, then the divergence of $\D$ must match this charge density, but we do not deal with this case in the present work. The approach we take to ensure divergence-free property is to use the divergence-free reconstruction ideas of Balsara~\cite{Balsara2001},~\cite{Balsara2004},~\cite{Balsara2009} which makes use of known values of the normal component of $\D$ on the faces of the cell and then reconstruct the vector field inside the cell by enforcing appropriate constraints on its divergence. 

\begin{figure}
\begin{center}
\begin{tabular}{ccc}
\includegraphics[width=0.32\textwidth]{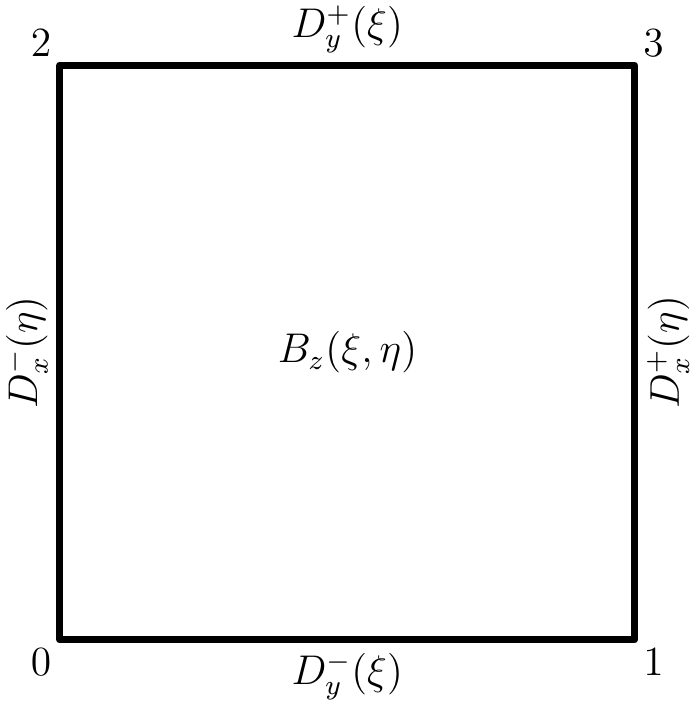} &
\includegraphics[width=0.32\textwidth]{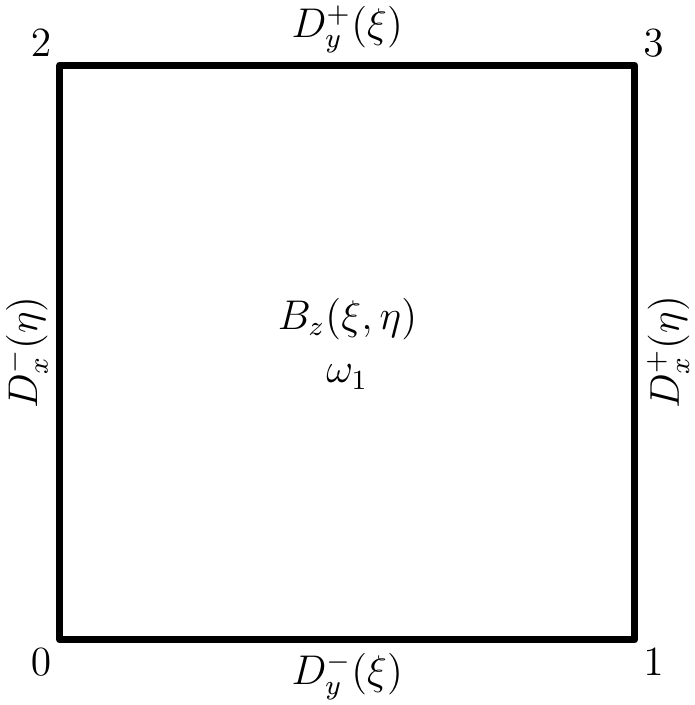} &
\includegraphics[width=0.32\textwidth]{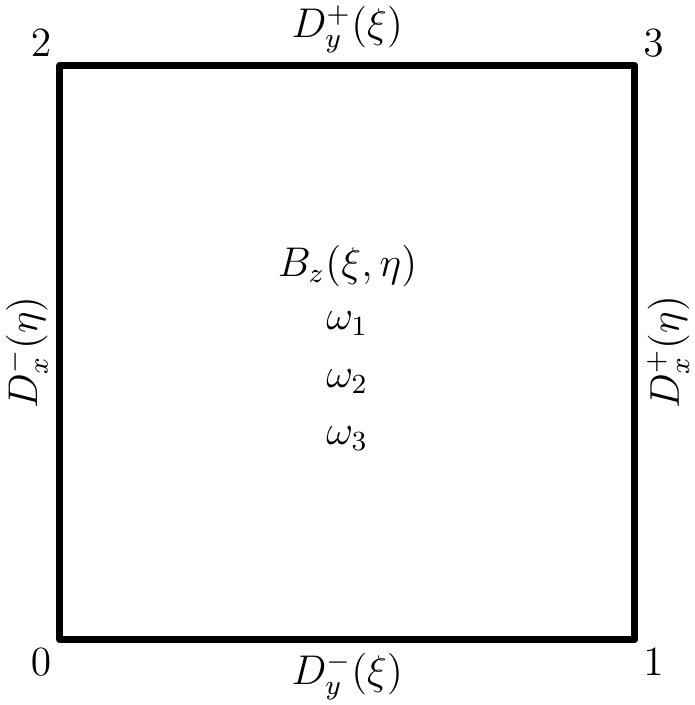} \\
(a) & (b) & (c)
\end{tabular}
\caption{Storage of solution variables: (a) $k=0,1,2$ (b) $k=3$ (c) $k=4$. For $k \ge 3$, we need extra information in addition to face solution.}
\label{fig:cell}
\end{center}
\end{figure}
We will approximate the normal components of $\D$ on the faces by one dimensional polynomials of degree $k \ge 0$. We map each cell to the reference cell $[-\half,+\half] \times [-\half,+\half]$ with coordinates $(\xi,\eta)$. Let $\poly_k(\xi)$ denote one dimensional polynomials of degree at most $k$ in the variable $\xi$, and $\poly_k(\xi,\eta)$ denote two dimensional polynomials of degree at most $k$. On the two vertical faces of a cell, the normal component is given by
\[
D_x^\pm(\eta) = \sum_{j=0}^k a_j^\pm \phi_j(\eta) \in \poly_k(\eta)
\]
while on the two horizontal faces, the corresponding normal component is given by
\[
D_y^\pm(\xi) = \sum_{j=0}^k b_j^\pm \phi_j(\xi) \in \poly_k(\xi)
\]
The location of these polynomials is illustrated in Figure~(\ref{fig:cell}). The basis functions $\phi_j$  are mutually orthogonal polynomials given by
\[
\phi_0(\xi) = 1, \quad
\phi_1(\xi) = \xi, \quad
\phi_2(\xi) = \xi^2 - \otw, \quad
\phi_3(\xi) = \xi^3 - \ttw \xi, \quad
\phi_4(\xi) = \xi^4 - \tft \xi^2 + \ttt
\]
\[
\phi_5(\xi) = \xi^5 - \frac{5}{18} \xi^3 + \frac{5}{336} \xi, \quad \textnormal{etc.}
\]
The magnetic field $B_z$ will be approximated inside each  cell by two dimensional polynomials of degree $k$ given by
\[
B_z(\xi,\eta) = \sum_{i=0}^{N(k)-1} \alpha_{i} \Phi_i(\xi,\eta) \in \poly_k(\xi,\eta), \qquad N(k) = \half(k+1)(k+2)
\]
where the two dimensional basis functions are given by
\begin{equation}
\begin{aligned}
\Phi_i \in \{ &1, \ \phi_1(\xi), \ \phi_1(\eta), &  \textrm{at second order} \\
& \phi_2(\xi), \ \phi_1(\xi)\phi_1(\eta), \ \phi_2(\eta), & \textrm{at third order} \\
& \phi_3(\xi), \ \phi_2(\xi)\phi_1(\eta), \ \phi_1(\xi)\phi_2(\eta), \ \phi_3(\eta), &  \textrm{at fourth order} \\
& \phi_4(\xi), \ \phi_3(\xi)\phi_1(\eta), \ \phi_2(\xi)\phi_2(\eta), \ \phi_1(\xi)\phi_3(\eta), \ \phi_4(\eta) \} & \textrm{at fifth order}
\end{aligned}
\label{eq:basis2d}
\end{equation}
Figure~(\ref{fig:cell}) shows the location of the above solution polynomials. By Gauss Theorem
\[
0 = \int_C \nabla\cdot \D \ud x \ud y = \int_{\partial C} \D \cdot \un
\]
which implies that
\begin{equation}
(a_0^+ - a_0^-) \dy + (b_0^+ - b_0^-) \dx = 0
\label{eq:compat}
\end{equation}
The above constraint will be satisfied by the initial condition, and the update scheme we devise will ensure that it is satisfied at future times also. Note that the above constraint depends only on the face averages of the solution variables stored on the faces.

The previous paragraphs have shown us how the modes that are the primal variables in our DG scheme are collocated (for the most part) at the faces of the mesh. While this is needed in order to formulate a globally constraint-preserving scheme, we should realize that we are actually interested in the solution of a PDE. Because of the Cauchy problem, the time-evolution of Maxwell's equations, just like the time-evolution of any hyperboloic PDE, relies on having all the spatial gradients. The reconstruction strategy that we describe below ensures that we can start with the modes at the skeleton (facial) mesh and obtain from it the variation of the electric displacement and magnetic induction at all locations on the mesh in a manner that is consistent with the involution constraint. Using the information of $D_x^\pm$, $D_y^\pm$ on the faces which are one dimensional polynomials of degree $k$, we have to reconstruct $\D = \D(\xi,\eta) \in \dfpoly_k(\xi,\eta)$ inside each cell such that the following conditions are satisfied.
\begin{enumerate}
\item The normal components of $\D(\xi,\eta)$ match the known values on the faces
\[
D_x(\pm\shalf,\eta) = D_x^\pm(\eta), \quad \forall \eta \in [-\shalf,+\shalf], \qquad D_y(\xi,\pm\shalf) = D_y^\pm(\xi), \quad \forall \xi \in [-\shalf,+\shalf]
\]

\item The divergence of $\D$ is zero everywhere inside the cell
\[
\nabla\cdot\D(\xi,\eta) = 0, \qquad \forall \xi, \eta \in [-\shalf,+\shalf]
\]
\end{enumerate}
The polynomial space $\dfpoly_k$ will be chosen so that the reconstruction problem is uniquely solvable. Note that we would like to have $(k+1)$'th order accurate approximations inside the cell which implies that $\poly_k(\xi,\eta) \subset \dfpoly_k(\xi,\eta)$ must be satisfied. However the space $\dfpoly_k$ must be bigger than $\poly_k$ in order to be able to satisfy the matching conditions on the faces and the divergence-free condition inside the cells. 
The precise form of the polynomial $\D(\xi,\eta)$ and the solution of the above reconstruction problem at various orders will be explained in the next section and in Appendix. For $k=0,1,2$ (upto third order accuracy), the reconstruction problem can be solved using the information of normal components on the faces but for $k \ge 3$, we require additional information from inside the cells to solve the reconstruction problem. For $k=3$ we specify an additional cell moment $\omega_1$ while for $k=4$, we specify three additional cell moments, $\omega_1, \omega_2, \omega_3$, see Figure~(\ref{fig:cell}).
\section{Divergence-free reconstruction of $\D$ inside a cell}
\label{sec:rec3}
In this section, we explain how to reconstruct the field $\D$ inside the cell given the values of the normal component on the faces of the cell, and in such a way that $\nabla\cdot\D = 0$. We explain the procedure for the case $k=3$ which leads to a fourth order approximation. The lower orders can be obtained from the fourth order solution and the reader can consult the Appendix. The fifth order case $(k=4)$ is also detailed in the Appendix. For solving the reconstruction problem, it is useful to note down the following results related to the 1-D orthogonal polynomials.
\[
\phi_1(\pm\shalf) = \pm \half, \quad \phi_2(\pm\shalf) = \frac{1}{6}, \quad \phi_3(\pm\shalf) = \pm \frac{1}{20}, \qquad \phi_4(\pm\shalf) = \frac{1}{70}
\]
\[
\phi_1'(\xi)=1, \quad \phi_2'(\xi) = 2 \phi_1(\xi), \quad \phi_3'(\xi) = 3 \phi_2(\xi) + \frac{1}{10}, \quad \phi_4'(\xi) = 4 \phi_3(\xi) + \frac{6}{35} \phi_1(\xi)
\]
\[
\phi_5'(\xi) = 5\phi_4(\xi)+\frac{5}{21}\phi_2(\xi)+\frac{1}{126}
\]
We will assume the following polynomial form for the vector field $\D$ inside the cell
\begin{align*}
D_x(\xi,\eta) = &\ a_{00} + a_{10}\phi_1(\xi) + a_{01}\phi_1(\eta) + a_{20}\phi_2(\xi) + a_{11} \phi_1(\xi)\phi_1(\eta) + a_{02}\phi_2(\eta) + \\
&\ a_{30} \phi_3(\xi) + a_{21}\phi_2(\xi)\phi_1(\eta) + a_{12}\phi_1(\xi)\phi_2(\eta) + a_{03}\phi_3(\eta) + a_{40}\phi_4(\xi) +\\
&\ a_{31}\phi_3(\xi)\phi_1(\eta) + a_{22}\phi_2(\xi)\phi_2(\eta) + a_{13}\phi_1(\xi)\phi_3(\eta) \\
D_y(\xi,\eta) = &\ b_{00} + b_{10}\phi_1(\xi) + b_{01}\phi_1(\eta) + b_{20}\phi_2(\xi) + b_{11}\phi_1(\xi)\phi_1(\eta) + b_{02}\phi_2(\eta) + \\
&\ b_{30}\phi_3(\xi) + b_{21}\phi_2(\xi)\phi_1(\eta) + b_{12}\phi_1(\xi) \phi_2(\eta) + b_{03} \phi_3(\eta) +  b_{31}\phi_3(\xi)\phi_1(\eta) + \\
&\ b_{22}\phi_2(\xi)\phi_2(\eta) + b_{13}\phi_1(\xi)\phi_3(\eta) + b_{04}\phi_4(\eta)
\end{align*}
Note that $D_x$ has the form of a polynomial $\poly_4(\xi,\eta)$ except that the terms corresponding to $\eta^4$ are not included. Similarly, $D_y$ belongs to $\poly_4(\xi,\eta)$ except for the term $\xi^4$ which is not included. Such polynomial spaces to approximate vector fields in a divergence conforming manner were introduced in~\cite{Brezzi1987} and are called BDFM polynomials. Both the components completely include $\poly_3(\xi,\eta)$ polynomials. Matching the cell solution to the face solution, we get the following 16 equations
\begin{align*}
 a_{00} \pm \shalf a_{10} + \tfrac{1}{6}a_{20} \pm \tfrac{1}{20}a_{30} + \tfrac{1}{70}a_{40} \; & = a_0^\pm \\
 a_{01} \pm \shalf a_{11} + \tfrac{1}{6}a_{21} \pm \tfrac{1}{20}a_{31}                       \; & = a_1^\pm \\
 a_{02} \pm \shalf a_{12} + \tfrac{1}{6}a_{22}        	                                  \; & = a_2^\pm \\
 a_{03} \pm \shalf a_{13}  		               	                                  \; & = a_3^\pm \\
 b_{00} \pm \shalf b_{01} + \tfrac{1}{6}b_{02} \pm \tfrac{1}{20}b_{03} + \tfrac{1}{70}b_{04} \; & = b_0^\pm \\
 b_{10} \pm \shalf b_{11} + \tfrac{1}{6}b_{12} \pm \tfrac{1}{20}b_{13}                       \; & = b_1^\pm \\
 b_{20} \pm \shalf b_{21} + \tfrac{1}{6}b_{22}        	                                  \; & = b_2^\pm \\
 b_{30} \pm \shalf b_{31}         	                                                  \; & = b_3^\pm
\end{align*}   
The divergence of the vector field $\D$ is a polynomial of degree 3 and making it zero inside the cell yields the following set of 10 equations
\begin{align*}
(a_{10}+\tfrac{1}{10}a_{30})\dy + (b_{01}+\tfrac{1}{10}b_{03})\dx &= 0 \\
(2 a_{20}+\tfrac{6}{35}a_{40})\dy + (b_{11}+b_{13}/10)\dx &= 0 \\
(a_{11}+a_{31}/10)\dy + (2 b_{02}+\tfrac{6}{35} b_{04})\dx &= 0 \\
 3 a_{30} \dy + b_{21} \dx &= 0 \\
2a_{21}\dy+2b_{12}\dx &= 0 \\
 a_{12}\dy+3b_{03}\dx &= 0 \\
 4 a_{40}\dy +   b_{31} \dx &= 0 \\
 3 a_{31}\dy + 2 b_{22} \dx &= 0 \\
 2 a_{22}\dy + 3 b_{13} \dx &=  0 \\
   a_{13}\dy + 4 b_{04} \dx &= 0
\end{align*}
The first equation in the above set is redundant since it is contained in the other equations due to the constraint~(\ref{eq:compat}). Ignoring this equation, we can solve for some of the coefficients $a_{ij}$, $b_{ij}$ in terms of the face solution as follows:

\boxed{
\begin{minipage}[h]{0.49\textwidth}
\begin{align*}
a_{00}=& \shalf (a_0^- + a_0^+)  + \tfrac{1}{12} (b_1^+ - b_1^-) \tfrac{\dx}{\dy} \\
a_{10}=& a_0^+ - a_0^- + \tfrac{1}{30}(b_2^+ - b_2^-) \tfrac{\dx}{\dy} \\ 
a_{20}=& -\shalf (b_1^+ - b_1^-) \tfrac{\dx}{\dy} + \tfrac{3}{140}(b_3^+ - b_3^-) \tfrac{\dx}{\dy}\\
a_{30}=& -\tfrac{1}{3}(b_2^+ - b_2^-) \tfrac{\dx}{\dy} \\
a_{03}=& \shalf (a_3^- + a_3^+) \\
a_{12}=& a_2^+ - a_2^- \\
a_{13}=& a_3^+ - a_3^- \\
a_{40}=& -\tfrac{1}{4}(b_3^+ - b_3^-) \tfrac{\dx}{\dy}
\end{align*}
\end{minipage}
\hspace{1em}
\begin{minipage}[h]{0.49\textwidth}
\begin{align*}
b_{00}=& \shalf (b_0^- + b_0^+) + \tfrac{1}{12} (a_1^+ - a_1^-) \tfrac{\dy}{\dx} \\
b_{01}=& b_0^+ - b_0^- + \tfrac{1}{30}(a_2^+ - a_2^-) \tfrac{\dy}{\dx} \\
b_{02}=& -\shalf (a_1^+ - a_1^-) \tfrac{\dy}{\dx} + \tfrac{3}{140}(a_3^+ - a_3^-) \tfrac{\dy}{\dx} \\
b_{30}=& \shalf (b_3^- + b_3^+) \\
b_{03}=& -\tfrac{1}{3} (a_2^+ - a_2^-) \tfrac{\dy}{\dx}  \\
b_{21}=& b_2^+ - b_2^- \\
b_{31}=& b_3^+ - b_3^- \\
b_{04}=& -\tfrac{1}{4} (a_3^+ - a_3^-) \tfrac{\dy}{\dx} \\
\end{align*}
\end{minipage}
}

\noindent
The remaining coefficients satisfy the following 9 equations

\begin{minipage}[h]{0.32\textwidth}
\[
\begin{aligned}
a_{01} + \tfrac{1}{6} a_{21} &= \shalf (a_1^+ + a_1^-), \\
a_{02} + \tfrac{1}{6} a_{22} &= \shalf (a_2^+ + a_2^-), \\
a_{11} + \tfrac{1}{10}a_{31} &= a_1^+ - a_1^-,
\end{aligned}
\]
\end{minipage}
\begin{minipage}[h]{0.32\textwidth}
\[
\begin{aligned}
b_{10} + \tfrac{1}{6} b_{12} &= \shalf (b_1^+ + b_1^-), \\
b_{20} + \tfrac{1}{6} b_{22} &= \shalf (b_2^+ + b_2^-), \\
b_{11} + \tfrac{1}{10}b_{13} &= b_1^+ - b_1^-, 
\end{aligned}
\]
\end{minipage}
\begin{minipage}[h]{0.32\textwidth}
\[
\begin{aligned}
2a_{21}\dy+2b_{12}\dx &=& 0 \\
 3 a_{31}\dy + 2 b_{22} \dx &=& 0 \\
 2 a_{22}\dy + 3 b_{13} \dx &=& 0
\end{aligned}
\]
\end{minipage}

\noindent
and we have more unknowns than equations. We can set the following coefficients which are not needed for fourth order accuracy to zero
\[
a_{31} = b_{22} = a_{22} = b_{13} = 0
\]
and we further obtain the solution for the following coefficients
\begin{equation}
\boxed{
a_{11} = a_1^+ - a_1^-,  \qquad a_{02}= \shalf (a_2^+ + a_2^-),  \qquad
b_{11} = b_1^+ - b_1^-, \qquad b_{20}= \shalf (b_2^+ + b_2^-)}
\end{equation}
The remaining unknowns satisfy the following set of equations
\begin{eqnarray}
\label{eq:bdfm31}
a_{01} + \frac{1}{6}a_{21} &=& \half (a_1^- + a_1^+) =: r_1 \\
\label{eq:bdfm32}
b_{10} + \frac{1}{6} b_{12} &=& \half (b_1^- + b_1^+) =: r_2 \\
\label{eq:bdfm33}
b_{12} \dx + a_{21} \dy &=& 0
\end{eqnarray}
We have four unknowns but only three equations. We cannot make any further assumptions regarding these coefficients since they are all at or below third degree, and we must retain all of them in order to get fourth order accuracy. The only way to complete the reconstruction is to provide an additional equation. Let us assume that we know the value of $\omega_1$ such that
\begin{equation}
b_{10} - a_{01}  = \omega_1
\label{eq:omg}
\end{equation}
Note that $\omega$ provides information about the mean value of the curl of the vector field in the cell. Then we can solve the equations to obtain
\begin{equation}
\boxed{
\begin{aligned}
a_{01} = \frac{1}{1+\frac{\dy}{\dx}}\left[r_1\frac{\dy}{\dx}+r_2 - \omega_1 \right], & \qquad
b_{10} = \omega_1 + a_{01} \\
a_{21} = 6(r_1 - a_{01}), & \qquad
b_{12} = 6(r_2 - b_{10})
\end{aligned}
}
\end{equation}
This completes the reconstruction of $\D$ inside the cell. Note that we had to introduce a cell moment to complete the reconstruction and the face solution alone is not sufficient to do this. We will know the value of $\omega_1$ from the initial condition and we have to device a scheme to evolve it forward in time which is explained in section~(\ref{sec:four}).
\section{Numerical scheme}
\label{sec:scheme}
We are now in a position to explain the constraint preserving scheme. Recall that we have several solution polynomials, some of which are stored on the faces and some are stored inside the cells. The basic solution variables have been illustrated in Figure~(\ref{fig:cell}). The solution polynomials $D_x(\xi,\eta)$, $D_y(\xi,\eta)$ are not independent and are obtained by the divergence-free reconstruction process described in section~(\ref{sec:rec3}) and in the Appendices. 
\begin{enumerate}
\item The normal component of $\D$ stored on the faces will be evolved by a flux reconstruction scheme applied on each face.
\item At fourth and fifth orders, we have additional quantities $\omega_i$ which are located inside the cells and we will devise a DG scheme for these quantities.
\item The magnetic flux density has only one component which is stored inside the cells and this will be evolved by a standard DG scheme.
\end{enumerate}
The FR and DG schemes require some numerical fluxes that are obtained from 1-D and 2-D Riemann problems. We give a short summary of these numerical fluxes in the Appendices~(\ref{sec:riem1d})-(\ref{sec:riem2d}).
\subsection{Flux reconstruction scheme for $\D$ on the faces}
Let us first describe the evolution scheme for the solution stored on the faces which is the normal component of $\D$. While we could use a DG scheme for this purpose, and this has been done by other researchers, in this work we will employ the flux reconstruction scheme which is also a high order numerical method for approximating the solutions of consevation laws. Like the spectral difference method~\cite{Kopriva1996},~\cite{Liu2006}, the flux reconstruction scheme~\cite{Huynh2007} is based on the differential formulation of conservation laws. By contrast, the DG schemes are based on an integral formulation. The basic idea is to first locally approximate the flux by a continuous polynomial using the piecewise discontinuous solution polynomial and some numerical fluxes coming from a Riemann solver. The solution is then updated to next time level using a collocation approach which avoids quadratures, which makes the method very efficient especially for 3-D problems. The construction of the continuous flux polynomial involves certain correction functions for which there are many possible choices available in the literature. Huynh~\cite{Huynh2007} proposed Radau polynomials as correction functions and later a more general family of correction functions were developed in~\cite{Vincent2011a} based on energy stability arguments in a Sobolev norm. This general correction function contains a parameter $c$ that is allowed take values in a certain interval and hence generates an infinite family of possible correction functions all of which lead to stable schemes. It has been discovered that FR schemes are equivalent to other high order schemes like spectral difference~\cite{Huynh2007} and certain types of nodal DG schemes~\cite{Huynh2007},~\cite{Vincent2011a},~\cite{DeGrazia2014},~\cite{Mengaldo2016} by choosing the correction functions appropriately, i.e., by choosing the parameter $c$. The nodal DG type schemes are recovered by using $c=0$ and this choice also leads to the most accurate numerical schemes~\cite{Castonguay2012},~\cite{Vincent2011}. The FR scheme has also been developed for advection-diffusion problems including Navier-Stokes equations and we refer the reader to the review article~\cite{Witherden2016a} for more references.
\begin{figure}
\begin{center}
\includegraphics[width=0.5\textwidth]{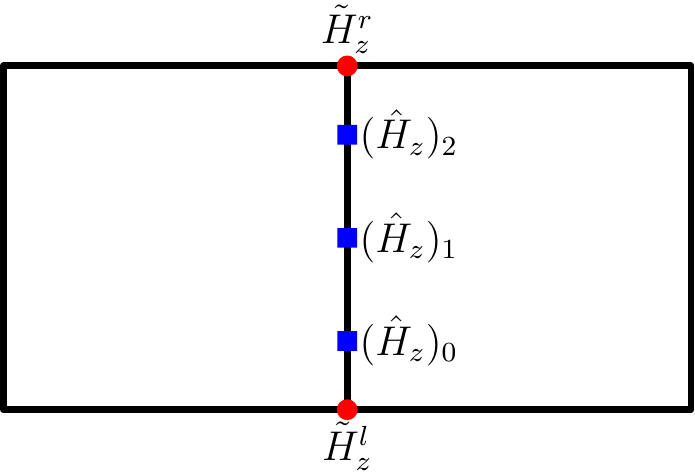}
\caption{Location of solution points of FR scheme on a vertical face for $k=2$. The blue squares are Gauss-Legendre quadrature points.}
\label{fig:face}
\end{center}
\end{figure}

Let us consider a vertical face on which $D_x$ is approximated by a polynomial of degree $k$ and we want to construct a scheme to update this information to the next time level. Note that $D_x$ evolves due to the $y$ derivatives of the magnetic field $B_z$ and hence the equation for $D_x$ can be discretized by a 1-D scheme on the vertical faces. Let us choose $k+1$ Gauss-Legendre quadrature points $\{\eta_i, 0 \le i \le k\}$ on the face, see Figure~(\ref{fig:face}), and let $\ell_j$, $j=0,1,\ldots,k$ be the corresponding Lagrange polynomials given by
\[
\ell_j(\eta) = \prod_{\substack{i=0 \\ i \ne j}}^k \left( \frac{\eta - \eta_i}{\eta_j - \eta_i} \right)
\]
Let us first interpolate $H_z$ using the Lagrange polynomials as follows
\[
H_z^\disc(\eta) = \sum_{j=0}^k (\hHz)_j \ell_j(\eta)
\]
The quantities $(\hHz)_j$ are values at the GL nodes as illustrated in Figure~(\ref{fig:face}); while we have a unique value of $D_x$ at each GL node of a vertical face, the other quantities $D_y$, $B_z$ are possibly discontinuous. Note that $D_y$ at the GL points are obtained from the reconstructed polynomial on the two cells sharing the face and $B_z$ is obtained by another polynomial inside the two neighbouring cells.  Since we have a 1-D Riemann problem at each GL point, we can compute $(\hHz)_j$ from a 1-D Riemann solver. This innovation adds a new aspect to FR schemes for involution constrained problems. This aspect is not present while solving usual conservation laws and arises because we are applying the FR scheme on the faces rather than over cells. Similarly, we will construct an interpolant $H_z^\disc(\xi)$ on each of the horizontal faces in the mesh. At any vertex, if we evaluate the flux $H_z^\disc$, we will get four different values from the four faces meeting at the vertex, and they need not agree with one another. In the next step, we will correct each of the interpolants to make them continuous at the vertices. 

At each vertex, we will have four different states that come together and define a 2-D Riemann problem and solution of this problem is briefly described in section~(\ref{sec:riem2d}).  Assume that we have computed the value of $H_z$ at all the vertices of the mesh using the 2-D Riemann solver. So at the bottom (l) and top (r) vertices, see Figure~(\ref{fig:face}), we know the unique values of $H_z$ which are given by the multidimensional Riemann solvers as $\tHz^l$, $\tHz^r$, respectively. We correct the above interpolant $H_z^\disc$ as follows
\[
H_z^\cont(\eta) = H_z^\disc(\eta) + [\tHz^l - H_z^\disc(-\shalf)] g_l(\eta) + [\tHz^r - H_z^\disc(+\shalf)] g_r(\eta)
\]
where the correction functions $g_l$, $g_r$ are polynomials of degree $k+1$ and have the property
\[
g_l(-\shalf) = g_r(+\shalf) = 1, \qquad g_l(+\shalf) = g_l(-\shalf) = 0
\]
This implies that
\[
H_z^\cont(-\shalf) = \tHz^l, \qquad H_z^\cont(+\shalf) = \tHz^r
\]
and hence $H_z^\cont(\eta)$ is a polynomial of degree $k+1$ and is continuous at the vertices. Following Huynh~\cite{Huynh2007}, the correction functions will be taken to be Radau polynomials of degree $k+1$ which also corresponds to taking $c=0$ in the general class of functions derived in~\cite{Vincent2011a}. The correction functions are given by
\[
g_l(\eta) = \frac{(-1)^k}{2}[ L_k(2\eta) - L_{k+1}(2\eta)], \qquad g_r(\eta) = \frac{1}{2}[ L_k(2\eta) + L_{k+1}(2\eta)]
\]
where $L_k : [-1,+1] \to \re$ is the Legendre polynomial of degree $k$. Note that we have to do a scaling of $\eta$ since in our convention we take $\eta \in [-\half,+\half]$. The normal component $D_x$ can now be updated by using a collocation approach
\begin{equation}
\df{D_x}{t}(\eta_i) = R_i := \frac{1}{\dy} \df{H_z^\cont}{\eta}(\eta_i), \qquad 0 \le i \le k
\label{eq:fr}
\end{equation}
Since $H_z^\cont$ is of degree $k+1$, the right hand side is of degree $k$ which agrees with the degree of the solution polynomial of $D_x$; hence the collocation scheme completely specifies the update for the facial solution.

The above discussion of the FR scheme shows that it is natural to use nodal basis functions which makes the collocation scheme easy to implement. However, in our approach, we actually represent $D_x(\eta)$ using orthogonal polynomials where solution coefficients are modal values and not nodal values, and hence the above nodal collocation update cannot be directly used. The use of orthogonal basis functions is convenient to write the solution of the divergence-free reconstruction problem in a simple form. The scheme for modal coefficients can be easily obtained by making a simple transformation. Let $V \in \re^{(k+1)\times(k+1)}$ be the Vandermonde matrix of the orthogonal polynomials corresponding to the GL nodes, i.e.,
\[
V_{ij} = \phi_j(\eta_i), \qquad 0 \le i,j \le k
\]
Then the update of the modal coefficients $a = [a_0, a_1, \ldots, a_k]^\top$ of $D_x(\eta)$ can be performed using the following equation
\[
\dd{a}{t} = V^{-1} R, \qquad R = [R_0, R_1, \ldots, R_k]^\top
\]
The Vandermonde matrix is also used to evaluate the polynomial $D_x(\eta)$ at the GL points. This matrix is common to each face and so that $V$ and it's inverse can be computed once in a pre-processing stage.
\subsection{Fourth order scheme for $\D$}
\label{sec:four}
At fourth order of accuracy ($k=3$), the face solution $\D$ does not completely determine the divergence-free reconstruction inside the cell. We have to specify $\omega_1$ as an additional information so that the reconstruction problem can be solved, which is defined as
\[
\omega_1 = b_{10} - a_{01} = 12\inttd [ D_y(\xi,\eta) \xi - D_x(\xi,\eta) \eta ]\ud\xi\ud\eta
\]
We will derive an evolution equation for $\omega_1$ using the induction equation. Since
\[
\frac{1}{12} \dd{\omega_1}{t} = \inttd \left(\df{D_y}{t} \xi - \df{D_x}{t} \eta \right) \ud\xi\ud\eta 
= -\inttd \left(\frac{1}{\dx} \df{H_z}{\xi} \xi + \frac{1}{\dy} \df{H_z}{\eta}\eta \right) \ud\xi\ud\eta
\]
Performing an integration by parts and using a numerical flux on the faces which is based on a 1-D Riemann solver, we obtain a semi-discrete DG scheme
\begin{eqnarray*}
\frac{1}{12}\dd{\omega_1}{t} 
&=&  - \frac{1}{\dx} \left[ \half \intod \hHz^{x-} \ud \eta + \half \intod \hHz^{x+} \ud \eta - \inttd H_z \ud \xi \ud \eta \right] \\
&& - \frac{1}{\dy} \left[ \half \intod \hHz^{y-} \ud \xi + \half \intod \hHz^{y+} \ud \xi - \inttd H_z \ud \xi \ud \eta \right]
\end{eqnarray*}
where the superscripts $x-$, $x+$ denotes the left and right faces of the cell, and $y-$, $y+$ denotes the bottom and top faces of the cell, and the fluxes $\hHz$ are obtained from the 1-D Riemann solver. The face integrals will be computed using $(k+1)$-point GL quadrature and the cell integral will be computed using tensor product of the same quadrature rule. Note that the quadrature points on the faces correspond to the solution points used in the FR scheme described in previous section; hence the flux $\hHz$ used in the $\omega_1$ equation is also used in the FR scheme described in the previous section.
\subsection{Fifth order scheme for $\D$}
\label{sec:five}
At fifth order of accuracy ($k=4$), the face solution $\D$ does not completely determine the divergence-free reconstruction inside the cell. We have to specify three cell moments $\omega_1, \omega_2,\omega_3$ as a additional information so that the reconstruction problem can be solved as shown in Appendix~(\ref{sec:rec4}). The update of $\omega_1$ has already been explained in previous section. The update equations for $\omega_2$ and $\omega_3$ can be derived in similar way. By definition
\[
\omega_2 = b_{20} - a_{11} = \inttd [180 D_y(\xi,\eta) \phi_2(\xi) - 144 D_x(\xi,\eta) \phi_1(\eta)\phi_1(\xi) ] \ud\xi\ud\eta
\]
and the semi-discrete scheme for the time evolution of $\omega_2$  is given by
\begin{align*}
\dd{\omega_2}{t} 
=&  - \frac{180}{\dx} \left[-\frac{1}{6} \int_{-\half}^\half \hHz^{x-} \ud \eta + \frac{1}{6} \int_{-\half}^\half \hHz^{x+} \ud \eta - \int_{-\half}^\half\int_{-\half}^\half 2 H_z \xi \ud \xi \ud \eta \right] \\
& - \frac{144}{\dy} \left[ \half \int_{-\half}^\half \hHz^{y-} \xi \ud \xi + \half \int_{-\half}^\half \hHz^{y+} \xi \ud \xi - \int_{-\half}^\half\int_{-\half}^\half H_z \xi \ud \xi \ud \eta \right]
\end{align*}
Similarly, we have
\[
\omega_3 = b_{11} - a_{02} = \inttd [144 D_y(\xi,\eta) \phi_1(\eta) \phi_1(\xi) - 180 D_x(\xi,\eta) \phi_2(\eta)] \ud\xi\ud\eta
\]
whose semi-discrete time evolution scheme is given by
\begin{align*}
\dd{\omega_3}{t} 
=&  - \frac{144}{\dx} \left[ \frac{1}{2} \int_{-\half}^\half \hHz^{x-} \eta \ud \eta + \frac{1}{2} \int_{-\half}^\half \hHz^{x+} \eta \ud \eta - \int_{-\half}^\half\int_{-\half}^\half H_z \eta \ud \xi \ud \eta \right] \\
& - \frac{180}{\dy} \left[-\frac{1}{6} \int_{-\half}^\half \hHz^{y-} \ud \xi + \frac{1}{6} \int_{-\half}^\half \hHz^{y+} \ud \xi - \int_{-\half}^\half\int_{-\half}^\half 2 H_z \eta \ud \xi \ud \eta \right]
\end{align*}
The face integrals will be computed using $(k+1)$-point GL quadrature and the cell integral will be computed using tensor product of the same quadrature rule. Note that the quadrature points on the faces correspond to the solution points used in the FR scheme; hence the flux $\hHz$ used in the $\omega_2, \omega_3$ equation is also used in the FR scheme applied on the faces.
\subsection{Discontinuous Galerkin method for $B_z$ inside cells}
In the 2-D model of Maxwell's equations that is considered in this paper, there is only one component of $\B$ so that we do not have to consider any constraint on this quantity. The magnetic flux $B_z$ is approximated by a two dimensional  polynomial $\poly_k$  inside each cell and we apply a standard DG scheme for this quantity, given by
\[
\begin{aligned}
\inttd \df{B_z}{t} \Phi_i(\xi,\eta) \ud\xi \ud\eta + & \inttd \left[ - \frac{1}{\dx} E_y \df{\Phi_i}{\xi} + \frac{1}{\dy} E_x \df{\Phi_i}{\eta} \right] \ud\xi\ud\eta \\
+ & \frac{1}{\dx} \intod \hEy^{x+} \Phi_i(+\shalf,\eta)\ud\eta - \frac{1}{\dx} \intod \hEy^{x-} \Phi_i(-\shalf,\eta)\ud\eta \\
- & \frac{1}{\dy} \intod \hEx^{y+} \Phi_i(\xi,+\shalf)\ud\xi + \frac{1}{\dy} \intod \hEx^{y-} \Phi_i(\xi,-\shalf)\ud\xi = 0
\end{aligned}
\]
where the test functions $\Phi_i$, $i=0,1,\ldots,N(k)-1$ are the basis functions of $\poly_k(\xi,\eta)$ as given in~(\ref{eq:basis2d}), $\hEy^{x-}$, $\hEy^{x+}$ are the values on the left and right faces obtained from the 1-D Riemann solver, and, $\hEx^{y-}$, $\hEx^{y+}$ are the values on bottom and top faces obtained from the 1-D Riemann solver. The integral inside the cell is evaluated using a tensor product of $(k+1)$-point GL quadrature while the face integrals are evaluated using $(k+1)$-point GL quadrature. The numerical fluxes used on the faces are common to the FR scheme and schemes for $\omega_i$ at fourth and fifth order accuracy.
\paragraph{Remark}
Note that in 3-D, we would approximate $\B$ in the same way as we approximate $\D$, i.e., the normal components of $\B$ are approximated on the faces, and the value inside the cell is obtained by a divergence-free reconstruction process. The evolution of $\B$ would then also follow similar approach as used for $\D$.
\subsection{Compatibility condition}
We have completely specified the semi-discrete scheme for all the variables which leads to a system of ODE. The update in time will be performed by standard time integration schemes. To solve the reconstruction problem, we must ensure that the compatibility condition~(\ref{eq:compat}) will be satisfied by the solution at future times also, assuming that it is satisfied by the initial condition. Such a scheme will then be refered to as being constraint preserving. Consider any cell $C$; using the flux reconstruction scheme~(\ref{eq:fr}) and $(k+1)$-point GL quadrature to integrate the normal component of $\D$ on the cell faces, we get
\[
\dd{}{t} \int_{\partial C} (\D \cdot \un)\ud s =  \sum_{i=0}^k \left[ \df{H_z^{\cont,x+}}{\eta}(\eta_i) - \df{H_z^{\cont,x-}}{\eta}(\eta_i) \right] \varpi_i + 
  \sum_{i=0}^k \left[ - \df{H_z^{\cont,y+}}{\xi}(\xi_i) + \df{H_z^{\cont,y-}}{\xi}(\xi_i) \right] \varpi_i
\]
where the $\varpi_i$ are the GL quadrature weights. The terms of the form $\df{H_z^\cont}{\xi}$, $\df{H_z^\cont}{\eta}$ are polynomials of degree $k$ and the quadrature is exact for such a polynomial, so that we can replace the sums with integrals
\begin{eqnarray*}
\dd{}{t} \int_{\partial C} (\D \cdot \un)\ud s &=& \intod \df{H_z^{\cont,x+}}{\eta}(\eta)\ud\eta  - \intod \df{H_z^{\cont,x-}}{\eta}(\eta) \ud\eta \\
&& - \intod   \df{H_z^{\cont,y+}}{\xi}(\xi)\ud\xi + \intod \df{H_z^{\cont,y-}}{\xi}(\xi) \ud \xi \\
&=& [(\tHz)_3 - (\tHz)_1] - [(\tHz)_2 - (\tHz)_0] - [(\tHz)_3 - (\tHz)_2] + [(\tHz)_1 - (\tHz)_0] \\
&=& 0
\end{eqnarray*}
where the subscripts on $\tHz$ denote the vertices of the cell, see Figure~(\ref{fig:cell}). This implies that
\begin{equation}
\dd{}{t}[(a_0^+ - a_0^-)\dy + (b_0^+ - b_0^-)\dx] = 0
\label{eq:c1}
\end{equation}
Hence under any time integration scheme, the compatibility condition~(\ref{eq:compat}) will be satisfied by our scheme assuming it holds for the initial condition. We see that the critical property required to achieve constraint preservation was to discretize the PDE on the faces and to use a unique value of $H_z$ at the vertices of the cells which comes from a 2-D Riemann solver.
\section{Energy stability analysis}
\label{sec:stab}
We will consider the energy stability of the first order scheme for constant $\varepsilon$ and $\mu$ with periodic boundary conditions. Balsara and K\"appeli~\cite{Balsara2018a} have performed Fourier stability analysis of fully discrete schemes and have derived CFL numbers for different time integration schemes. Here we perform direct energy stability analysis of the semi-discrete scheme. To aid in the proof, we introduce the usual $(i,j)$ indexing notation for the cells and half indices will be used to denote the faces and vertices.  At first order, our solution variables consist of face averages of normal components of $\D$ and cell average of $B_z$. The scheme is given by
\[
\dd{}{t} (D_x)_{\iph,j} = \frac{ \Hziphjph - \Hziphjmh }{\Delta y}, \qquad \dd{}{t}(D_y)_{i,\jph} = - \frac{\Hziphjph - \Hzimhjph}{\Delta x}
\]
\[
\dd{}{t}(B_z)_{i,j} = - \frac{\Eyiphj - \Eyimhj}{\Delta x} + \frac{\Exijph - \Exijmh}{\Delta y}
\]
Define the {\em total energy}
\[
\tote_h^*(t) = \sum_i \sum_j \frac{1}{2\varepsilon} (D_x)_{\iph,j}^2 \Delta x \Delta y + \sum_i \sum_j \frac{1}{2\varepsilon} (D_y)_{i,\jph}^2 \Delta x \Delta y + \sum_i \sum_j \frac{1}{2\mu} (B_z)_{i,j}^2 \Delta x \Delta y
\]
Then using the above scheme, the rate of change of energy is given by
\begin{eqnarray*}
\dd{\tote_h^*}{t} &=& \sum_i \sum_j \Exiphj [ \Hziphjph - \Hziphjmh ] \Delta x \\
&& - \sum_i \sum_j \Eyijph [\Hziphjph - \Hzimhjph] \Delta y \\
&& + \sum_i \sum_j \Hzij \left[ - \frac{\Eyiphj - \Eyimhj}{\Delta x} + \frac{\Exijph - \Exijmh}{\Delta y} \right] \Delta x \Delta y \\
&=:& \mathcal{P}_1 + \mathcal{P}_2 + \mathcal{P}_3
\end{eqnarray*}
To simplify the analysis below, define the average and difference operators
\[
\jumpx{\cdot}_{i,j} = (\cdot)_{\iph,j} - (\cdot)_{\imh,j}, \qquad \jumpy{\cdot}_{i,j} = (\cdot)_{i,\jph} - (\cdot)_{i,\jmh}
\]
\[
\avgx{\cdot}_{i,j} = \half[ (\cdot)_{\iph,j} + (\cdot)_{\imh,j}], \qquad \avgy{\cdot}_{i,j} = \half[(\cdot)_{i,\jph} + (\cdot)_{i,\jmh}]
\]
Using summation by parts we can write
\[
\mathcal{P}_1 = -\sum_i \sum_j \Hziphjph \jumpy{E_x}_{\iph,\jph} \Delta x, \quad
\mathcal{P}_2 = \sum_i \sum_j \Hziphjph \jumpx{E_y}_{\iph,\jph} \Delta y
\]
\[
\mathcal{P}_3 = \sum_i \sum_j \Eyiphj \jumpx{H_z}_{\iph,j} \Delta y - \sum_i \sum_j \Exijph \jumpy{H_z}_{i,\jph} \Delta x
\]
Let us use the fluxes obtained from a Riemann solver, see e.g.~\cite{Balsara2017} and also the Appendix,
\begin{eqnarray*}
\Exijph &=& \frac{1}{4}\left[ \Eximhj + \Exiphj + \Eximhjpo + \Exiphjpo \right] + \frac{\mu c}{2}\jumpy{H_z}_{i,\jph} \\
&=& \half[ \avgy{E_x}_{\imh,\jph} + \avgy{E_x}_{\iph,\jph}] + \frac{\mu c}{2}\jumpy{H_z}_{i,\jph}
\end{eqnarray*}
\begin{eqnarray*}
\Eyiphj &=& \frac{1}{4}\left[ \Eyijmh + \Eyijph + \Eyipojmh + \Eyipojph \right] - \frac{\mu c}{2}\jumpx{H_z}_{\iph,j} \\
&=& \half[ \avgx{E_y}_{\iph,\jmh} + \avgx{E_y}_{\iph,\jph} ] - \frac{\mu c}{2}\jumpx{H_z}_{\iph,j}
\end{eqnarray*}
\[
\begin{aligned}
\Hziphjph = & \ \frac{1}{4} \left[ \Hzij + \Hzipoj + \Hzijpo + \Hzipojpo \right] \\
& + \frac{\varepsilon c}{2} \jumpy{E_x}_{\iph,\jph} - \frac{\varepsilon c}{2} \jumpx{E_y}_{\iph,\jph}
\end{aligned}
\]
Note that the fluxes consist of a central part and some additional terms that depend on jumps in the solution variables. Then
\begin{eqnarray*}
\mathcal{P}_1 &=& - \sum_i \sum_j \half [ \avgy{H_z}_{i,\jph} + \avgy{H_z}_{i+1,\jph} ] \jumpy{E_x}_{\iph,\jph} \Delta x \\
&& - \frac{\varepsilon c}{2} \sum_i \sum_j [ \jumpy{E_x}_{\iph,\jph} -  \jumpx{E_y}_{\iph,\jph} ]  \jumpy{E_x}_{\iph,\jph} \Delta x \\
&=:& \mathcal{P}_4 + \mathcal{P}_5
\end{eqnarray*}

\begin{eqnarray*}
\mathcal{P}_2 &=&  \sum_i \sum_j \half [ \avgx{H_z}_{\iph,j} + \avgx{H_z}_{\iph,j+1} ] \jumpx{E_y}_{\iph,\jph} \Delta y \\
&& + \frac{\varepsilon c}{2} \sum_i \sum_j [ \jumpy{E_x}_{\iph,\jph} -  \jumpx{E_y}_{\iph,\jph} ]  \jumpx{E_y}_{\iph,\jph} \Delta y \\
&=:& \mathcal{P}_6 + \mathcal{P}_7
\end{eqnarray*}

\begin{eqnarray*}
\mathcal{P}_3 &=& \sum_i \sum_j \half[ \avgx{E_y}_{\iph,\jmh} + \avgx{E_y}_{\iph,\jph} ] \jumpx{H_z}_{\iph,j} \Delta y \\
&& - \sum_i \sum_j \half[ \avgy{E_x}_{\imh,\jph} + \avgy{E_x}_{\iph,\jph}] \jumpy{H_z}_{i,\jph} \Delta x \\
&& - \frac{\mu c}{2} \sum_i \sum_j \left\{ [ \jumpx{H_z}_{\iph,j} ]^2 \Delta y + [ \jumpy{H_z}_{i,\jph} ]^2 \Delta x \right\} \\
&=:& \mathcal{P}_8 + \mathcal{P}_9 + \mathcal{D}_1
\end{eqnarray*}
Note that $\mathcal{D}_1 \le 0$. Now consider
\begin{eqnarray*}
\mathcal{P}_4 + \mathcal{P}_9 &=& - \sum_i \sum_j \half [ \avgy{H_z}_{i,\jph} + \avgy{H_z}_{i+1,\jph} ] \jumpy{E_x}_{\iph,\jph} \Delta x \\
&& - \sum_i \sum_j \half[ \avgy{E_x}_{\imh,\jph} + \avgy{E_x}_{\iph,\jph}] \jumpy{H_z}_{i,\jph} \Delta x \\
&=& -\half \sum_i \sum_j [ \avgy{H_z}_{i,\jph} \jumpy{E_x}_{\iph,\jph} + \avgy{E_x}_{\iph,\jph} \jumpy{H_z}_{i,\jph} ] \Delta x \\
&& - \half \sum_i\sum_i  \avgy{H_z}_{i+1,\jph} \jumpy{E_x}_{\iph,\jph}  \Delta x \\
&& - \half \sum_i\sum_i \avgy{E_x}_{\imh,\jph} \jumpy{H_z}_{i,\jph} \Delta x
\end{eqnarray*}
In the third term on the right, shift the $i$ index back by one to obtain
\begin{eqnarray*}
\mathcal{P}_4 + \mathcal{P}_9 
&=& -\half \sum_i \sum_j [ \avgy{H_z}_{i,\jph} \jumpy{E_x}_{\iph,\jph} + \avgy{E_x}_{\iph,\jph} \jumpy{H_z}_{i,\jph} ] \Delta x \\
&& - \half \sum_i\sum_i  [ \avgy{H_z}_{i,\jph} \jumpy{E_x}_{\imh,\jph} + \avgy{E_x}_{\imh,\jph} \jumpy{H_z}_{i,\jph} ] \Delta x \\
&=& 0
\end{eqnarray*}
This follows because the term in each sum is a perfect difference. We show this for the first term.
\begin{eqnarray*}
&& \avgy{H_z}_{i,\jph} \jumpy{E_x}_{\iph,\jph} + \avgy{E_x}_{\iph,\jph} \jumpy{H_z}_{i,\jph} \\
&=& \frac{\Hzij + \Hzijpo}{2} [ \Exiphjpo - \Exiphj] + \frac{\Exiphj + \Exiphjpo}{2} [ \Hzijpo - \Hzij] \\
&=& \Hzijpo \Exiphjpo - \Hzij \Exiphj
\end{eqnarray*}
Similarly, we can show that $\mathcal{P}_6 + \mathcal{P}_8 = 0$. If $\Delta x = \Delta y = h$ then
\[
\mathcal{D}_2 := \mathcal{P}_5 + \mathcal{P}_7 = - \frac{\varepsilon c}{2} \sum_i \sum_j [ \jumpy{E_x}_{\iph,\jph} -  \jumpx{E_y}_{\iph,\jph} ]^2 h \le 0
\]
Hence we obtain
\[
\dd{\tote_h^*}{t} = \mathcal{D}_1 + \mathcal{D}_2 \le 0
\]
If $\Delta x \ne \Delta y$, then we cannot prove that $\mathcal{D}_2 \le 0$. In this case we have to modify the flux $\Hziphjph$ slightly as follows
\[
\begin{aligned}
\Hziphjph = & \ \frac{1}{4} \left[ \Hzij + \Hzipoj + \Hzijpo + \Hzipojpo \right] \\
& + \frac{\varepsilon c h}{2\Delta y} \jumpy{E_x}_{\iph,\jph} - \frac{\varepsilon c h}{2\Delta x} \jumpx{E_y}_{\iph,\jph}
\end{aligned}
\]
where $h$ could be defined as $h = \max(\Delta x, \Delta y)$. Then
\[
\mathcal{D}_2 = - \frac{\varepsilon c h}{2} \sum_i \sum_j \left[ \frac{\jumpy{E_x}_{\iph,\jph}}{\Delta y} -  \frac{\jumpx{E_y}_{\iph,\jph}}{\Delta x} \right]^2 \Delta x \Delta y \le 0
\]
and we can again prove energy stability. Note that the dissipation term $\mathcal{D}_2$  is created due to the curl of the electric field which is physically meaningful for the Maxwell model. We also observe that the dissipation is due to the additional terms in the numerical flux involving the jumps in solution variables, and the central part of the flux would lead to energy conservation.
\paragraph{Remark} The energy $\tote_h^*$ we have analyzed above is the energy of the solution on the faces and is not the true energy, which is defined as
\[
\tote_h = \int_\Omega \left[ \frac{1}{2\varepsilon}(D_x^2 + D_y^2) + \frac{1}{2\mu} B_z^2 \right] \ud x \ud y
\]
where $D_x, D_y$ inside the integral are obtained by the divergence-free reconstruction scheme. Note that, using the inequality $ab \le (a^2 + b^2)/2$, we get
\begin{eqnarray*}
\int_\Omega D_x^2 \ud x \ud y &=& \sum_i \sum_j \left[ \frac{1}{3} (D_x)_{\imh,j}^2 + \frac{1}{3} (D_x)_{\iph,j}^2 + \frac{1}{3} (D_x)_{\imh,j} (D_x)_{\iph,j} \right] \Delta x \Delta y \\
&\le& \sum_i \sum_j \left[ \frac{1}{2} (D_x)_{\imh,j}^2 + \frac{1}{2} (D_x)_{\iph,j}^2 \right] \Delta x \Delta y \\
&=& \sum_i \sum_j  (D_x)_{\iph,j}^2 \Delta x \Delta y
\end{eqnarray*}
and, using the inequality $ab \ge -(a^2 + b^2)/2$, we get
\[
\int_\Omega D_x^2 \ud x \ud y \ge \frac{1}{3} \sum_i \sum_j  (D_x)_{\iph,j}^2 \Delta x \Delta y
\]
with similar results for the $D_y$ component. Hence it follows that
\[
\tote_h \le \tote_h^* \le 3 \tote_h
\]
and so $\tote_h^*$ is an equivalent energy norm.
\section{Numerical results}
\label{sec:res}
Our numerical scheme belongs to the class of so called RKDG methods where a DG/FR scheme is used for spatial discretization and the resulting system of ODE are solved using a Runge-Kutta scheme. For degrees $k=0,1,2$, we use the first, second and third order strong stability preserving RK schemes~\cite{Shu1988},~\cite{Shu1989}, respectively, while for $k=3$ and $k=4$ we use the 5-stage, 4-th order strong stability preserving RK scheme~\cite{Spiteri2002},~\cite{Spiteri2003} or the classical fourth order RK scheme. The time step is computed from the CFL number which is defined as
\[
\cfl = \max\left\{ \max\frac{c\Delta t}{\dx}, \max\frac{c\Delta t}{\dy} \right\}
\]
where $c = \frac{1}{\sqrt{\mu\varepsilon}}$ is the speed of light, and the inner maximum is taken over the whole mesh. The $\cfl$ numbers have been derived in Balsara \& K\"appeli~\cite{Balsara2018a} using Fourier stability analysis.

We will measure the error in the solution using $L^1$ and $L^2$ norms. These norms are defined as follows for vector and scalar functions
\[
\norm{\D}_{L^1} = \frac{1}{|\Omega|}\int_\Omega \norm{\D} \ud x \ud y, \qquad \norm{\D}_{L^2} = \left( \frac{1}{|\Omega|} \int_\Omega \norm{\D}^2 \ud x \ud y \right)^\half, \qquad \norm{\D} = \sqrt{D_x^2 + D_y^2}
\]
\[
\norm{B_z}_{L^1} = \int_\Omega |B_z| \ud x \ud y, \qquad \norm{B_z}_{L^2} = \left( \int_\Omega B_z^2 \ud x \ud y \right)^\half
\]
and the integrals are computed using a tensor product of $(k+2)$-point Gauss-Legendre quadrature. Note that we measure the error norm of the solution polynomials relative to the reference or exact solutions and not just the error in the cell average value. In particular, the error in $\D$ is measured based on the reconstructed field inside the cells.
\subsection{Plane wave propagation}
This test case describes the propagation of a plane electromagnetic wave in vacuum. The purpose of this test case is to check the accuracy of our numerical method since we know the exact solution.  The simulation is performed in a square domain of $[-0.5, 0.5]\times[-0.5, 0.5]~\si{\square \m}$ divided in \numlist{16;32;64;128} square cells in each direction with 
periodic boundary conditions. The simulation is conducted for a time duration of \SI{3.5}{\nano \s}. The initial condition of $\B$ and $\D$ field is specified from magnetic vector potential $\bm{A}(x,y,t)$ and electric  vector potential $\bm{C}(x,y,t)$ and using the relationships $\B = \nabla \times \bm{A}$ and $\D = c \epsilon_0\nabla \times \bm{C}$.  The magnetic and electric vector potentials are given by
\[
\bm{A}(x,y,t) = \frac{1}{2\pi}\sin[2\pi(x+y-\sqrt{2}ct)] \hat{e}_y, \qquad
 \bm{C}(x,y,t) = -\frac{1}{2\pi\sqrt{2}}\sin[2\pi(x+y-\sqrt{2}ct)] \hat{e}_z
\]
The convergence of the error for different degree and grid sizes are shown in tables~(\ref{tab:pwavek1finalNONE})-(\ref{tab:pwavek4finalNONE}). We observe that with degree $k$ solution on the faces, all the quantities converge at the rate of $O(h^{k+1})$ under mesh refinement demonstrating that optimal accuracy is achieved by our method.
\begin{table}
 \begin{center} 
 \begin{tabular}{|c|c|c|c|c|c|c|c|c|}
 \hline
 $N_x\times N_y$ & $\|\D^h-\D\|_{L^1}$ & Ord & $\|\D^h-\D\|_{L^2}$ &Ord & $\|B_z^h-B_z\|_{L^1}$ & Ord & $\|B_z^h-B_z\|_{L^2}$ & Ord\\ 
 \hline 
 $ 16 \times 16 $ & 9.5684e-05 & --- &1.0646e-04 & --- & 3.6205e-02 & --- & 4.1073e-02 & --- \\ 
$ 32 \times 32 $ & 1.7320e-05 & 2.47 & 1.9014e-05 & 2.49 & 6.6895e-03 & 2.44 & 7.4826e-03 & 2.46\\ 
$ 64 \times 64 $ & 3.7425e-06 & 2.21 & 4.0522e-06 & 2.23 & 1.4458e-03 & 2.21 & 1.6167e-03 & 2.21\\ 
$ 128 \times 128 $ & 8.9361e-07 & 2.07 & 9.6327e-07 & 2.07 & 3.4400e-04 & 2.07 & 3.8629e-04 & 2.07 \\ 
 \hline 
 \end{tabular} 
 \end{center} 
 \caption{Plane wave test, degree=$1$: convergence of final error} 
 \label{tab:pwavek1finalNONE} 
 \end{table}

\begin{table}
 \begin{center} 
 \begin{tabular}{|c|c|c|c|c|c|c|c|c|}
 \hline
 $N_x\times N_y$ & $\|\D^h-\D\|_{L^1}$ & Ord & $\|\D^h-\D\|_{L^2}$ &Ord & $\|B_z^h-B_z\|_{L^1}$ & Ord & $\|B_z^h-B_z\|_{L^2}$ & Ord\\ 
 \hline 
 $ 16 \times 16 $ & 3.0156e-05 & --- &3.3378e-05 & --- & 1.1640e-02 & --- & 1.2951e-02 & --- \\ 
$ 32 \times 32 $ & 3.6268e-06 & 3.06 & 4.0113e-06 & 3.06 & 1.4044e-03 & 3.05 & 1.5605e-03 & 3.05\\ 
$ 64 \times 64 $ & 4.4793e-07 & 3.02 & 4.9534e-07 & 3.02 & 1.7378e-04 & 3.01 & 1.9310e-04 & 3.01\\ 
$ 128 \times 128 $ & 5.5780e-08 & 3.01 & 6.1684e-08 & 3.01 & 2.1671e-05 & 3.00 & 2.4074e-05 & 3.00 \\ 
 \hline 
 \end{tabular} 
 \end{center} 
 \caption{Plane wave test, degree=$2$: convergence of error} 
 \label{tab:pwavek2finalNONE} 
 \end{table}

\begin{table}
 \begin{center} 
 \begin{tabular}{|c|c|c|c|c|c|c|c|c|}
 \hline
 $N_x\times N_y$ & $\|\D^h-\D\|_{L^1}$ & Ord & $\|\D^h-\D\|_{L^2}$ &Ord & $\|B_z^h-B_z\|_{L^1}$ & Ord & $\|B_z^h-B_z\|_{L^2}$ & Ord\\ 
 \hline 
 $ 16 \times 16 $ & 3.0557e-07 & --- &3.5095e-07 & --- & 1.9458e-04 & --- & 2.5101e-04 & --- \\ 
$ 32 \times 32 $ & 1.1040e-08 & 4.79 & 1.3428e-08 & 4.71 & 1.1275e-05 & 4.11 & 1.4590e-05 & 4.10\\ 
$ 64 \times 64 $ & 5.0469e-10 & 4.45 & 6.1548e-10 & 4.45 & 6.7924e-07 & 4.05 & 8.9228e-07 & 4.03\\ 
$ 128 \times 128 $ & 2.6834e-11 & 4.23 & 3.3945e-11 & 4.18 & 4.1802e-08 & 4.02 & 5.5449e-08 & 4.01 \\ 
 \hline 
 \end{tabular} 
 \end{center} 
 \caption{Plane wave test, degree=$3$: convergence of error} 
 \label{tab:pwavek3finalNONE} 
 \end{table}

 \begin{table}
 \begin{center} 
 \begin{tabular}{|c|c|c|c|c|c|c|c|c|}
 \hline
 $N_x\times N_y$ & $\|\D^h-\D\|_{L^1}$ & Ord & $\|\D^h-\D\|_{L^2}$ &Ord & $\|B_z^h-B_z\|_{L^1}$ & Ord & $\|B_z^h-B_z\|_{L^2}$ & Ord\\ 
 \hline
$8\times 8$ & 2.4982e-07 & --- & 2.8435e-07&---&2.5514e-04&---&3.2612e-04&--- \\
$16\times 16$ & 6.3357e-09 & 5.30&7.2655e-09&5.29&6.4340e-06&5.31&8.5532e-06&5.25\\
$32\times 32$ & 1.7113e-10 & 5.21&2.0807e-10&5.13&1.9021e-07&5.08&2.5521e-07&5.07\\
$64\times 64$ & 5.1213e-12 & 5.06&6.4133e-12&5.02&5.9026e-09&5.01&7.9601e-09&5.00\\
$128\times 128$ & 1.5949e-13 & 5.00&2.0054e-13&5.00&1.8412e-10&5.00&2.4890e-10&5.00\\
 \hline 
 \end{tabular} 
 \end{center} 
 \caption{Plane wave test, degree=$4$: convergence of final error}               
 \label{tab:pwavek4finalNONE} 
 \end{table}





\subsection{Compact Gaussian electromagnetic pulse incident on a refractive disk}
This test case deals with scattering interaction of a compact electromagnetic pulse impinging upon a dielectric disc.  Simulations are performed in a domain $[-7.0, 7.0]\times[-7.0, 7.0]~\si{\square \m}$ upto the time 23.3 ns. A dielectric disc of radius \SI{0.75}{\m} is located at the center of the computational domain.  The initial condition is given by $\B = \nabla \times \bm{A}$ and $\D=c \epsilon_0 \nabla \times \bm{C}$ where the following magnetic and electric vector potential are used
\begin{align*}
\bm{A}(x,y,t)&=\frac{\lambda}{2\pi}\sin[2\pi(x+y)] e^{-\frac{(x-a)^2+(y-b)^2}{\chi^2}}\hat{e}_y,\\
 \bm{C}(x,y,t)&=-\frac{\lambda}{2\pi\sqrt{2}}\sin[2\pi(x+y)] e^{-\frac{(x-a)^2+(y-b)^2}{\chi^2}} \hat{e}_z
\end{align*}
where wavelength of the electromagnetic beam $\lambda = \SI{1.5}{\m}$, $\chi= \SI{1.5}{\m}$ and $(a,b)={(-2.5,2.5)}~\si{\m}$. The relative permittivity is taken as
\[
\epsilon_r(x,y) = 5.0 - 4.0 \tanh\left(\frac{\sqrt{x^2+y^2}-0.75}{0.08} \right)
\]
\begin{figure}
\begin{center}
\begin{tabular}{ccc}
\includegraphics[width=0.33\textwidth]{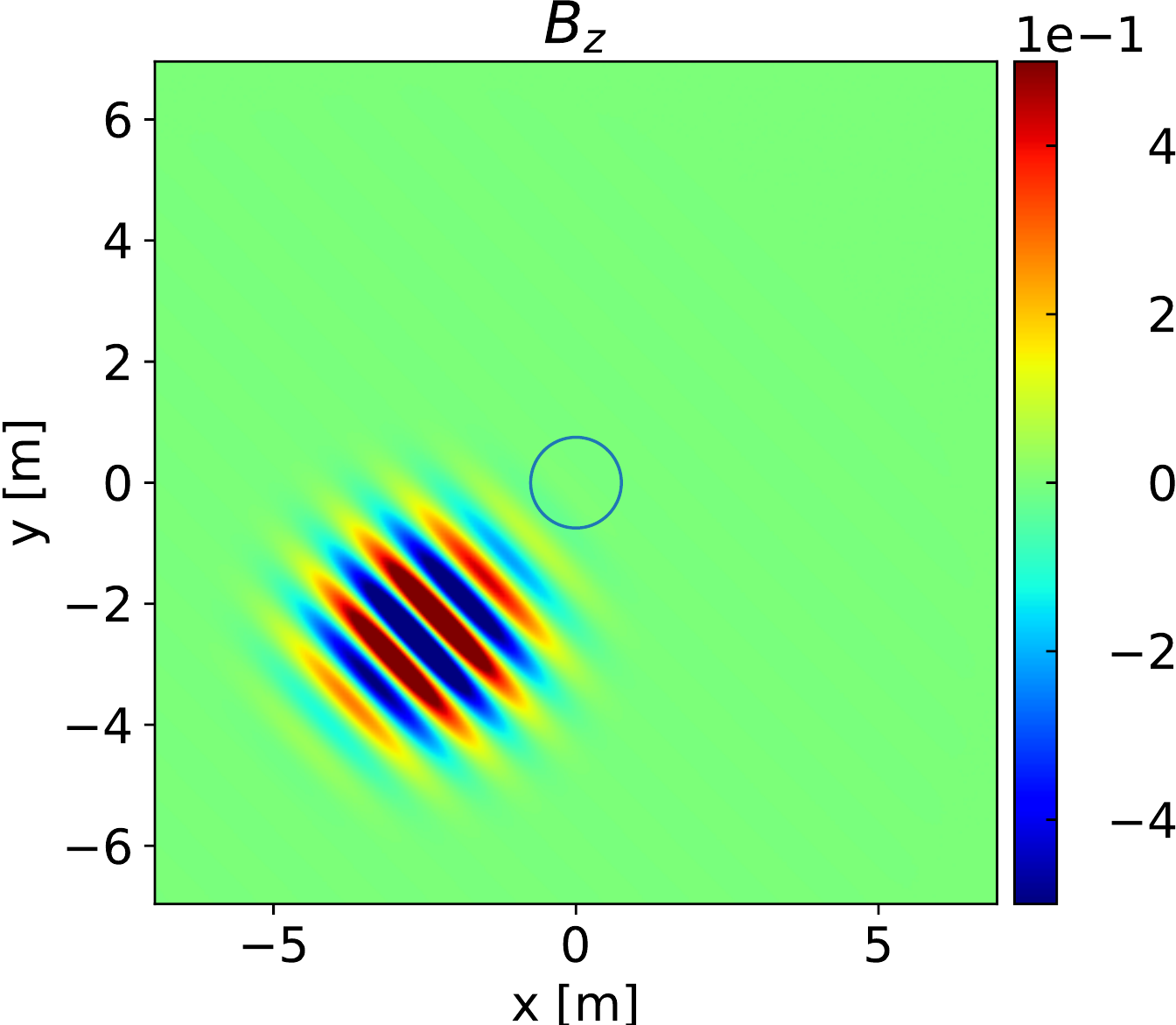} &
\includegraphics[width=0.33\textwidth]{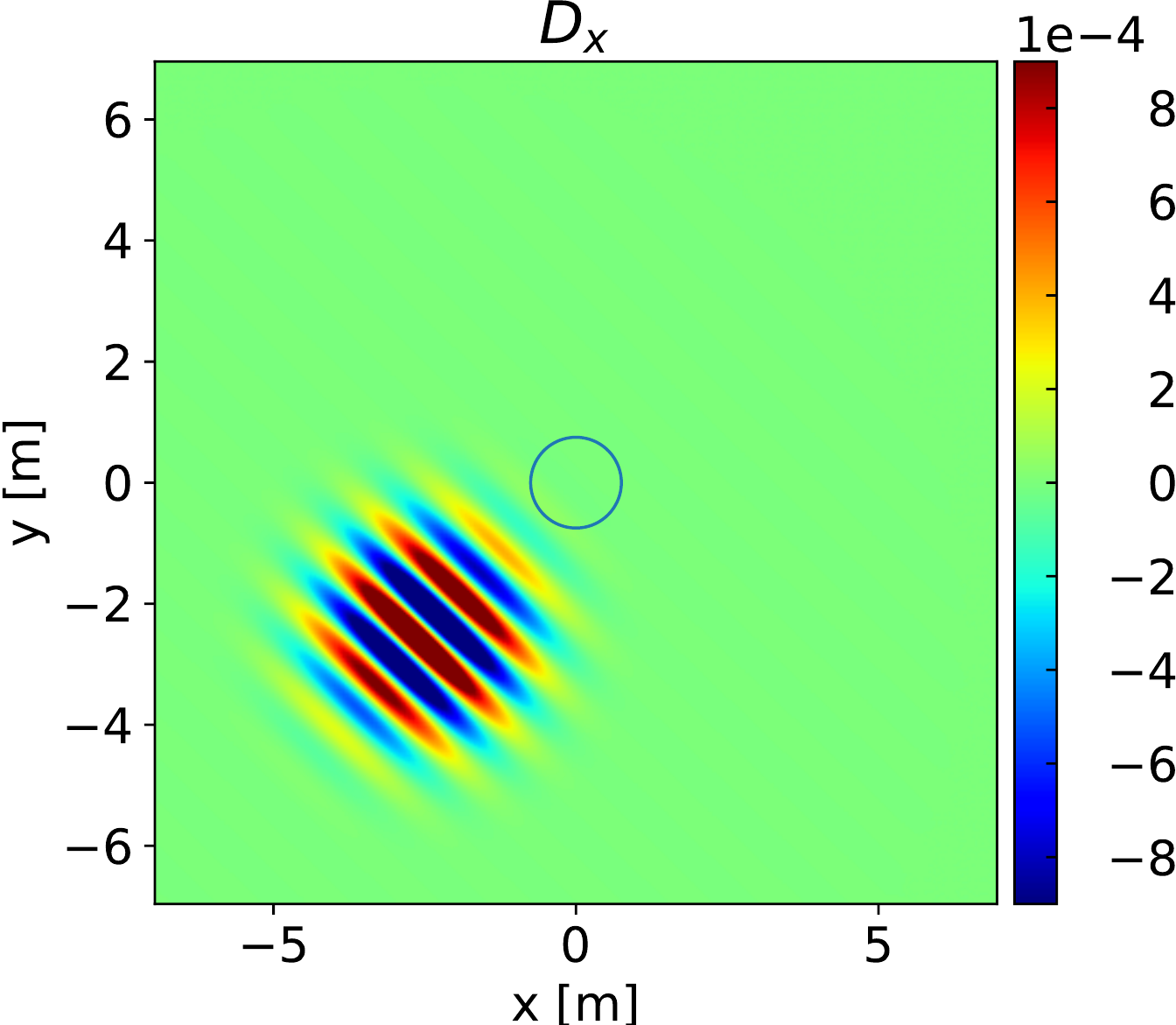} &
\includegraphics[width=0.33\textwidth]{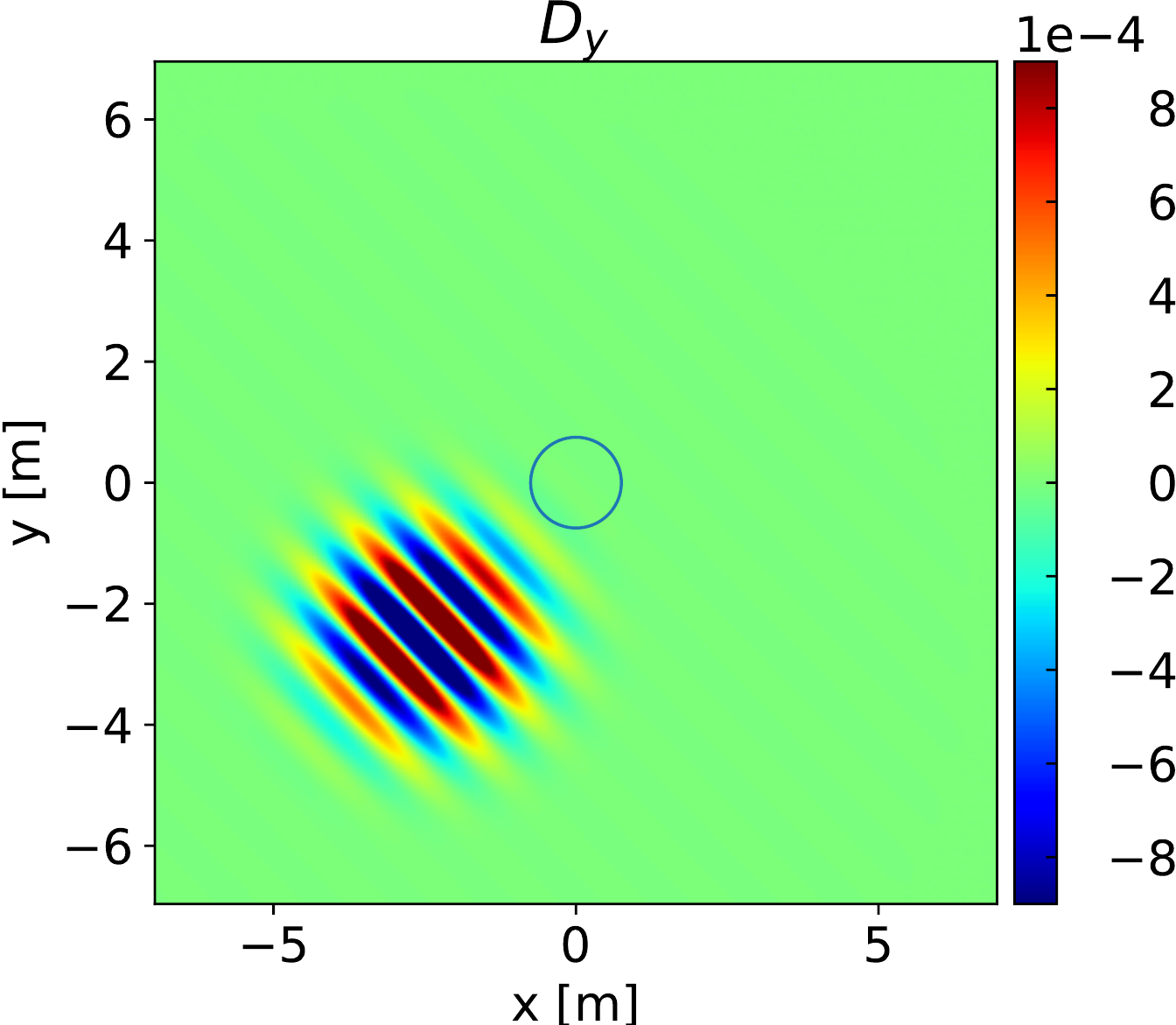} \\
\includegraphics[width=0.33\textwidth]{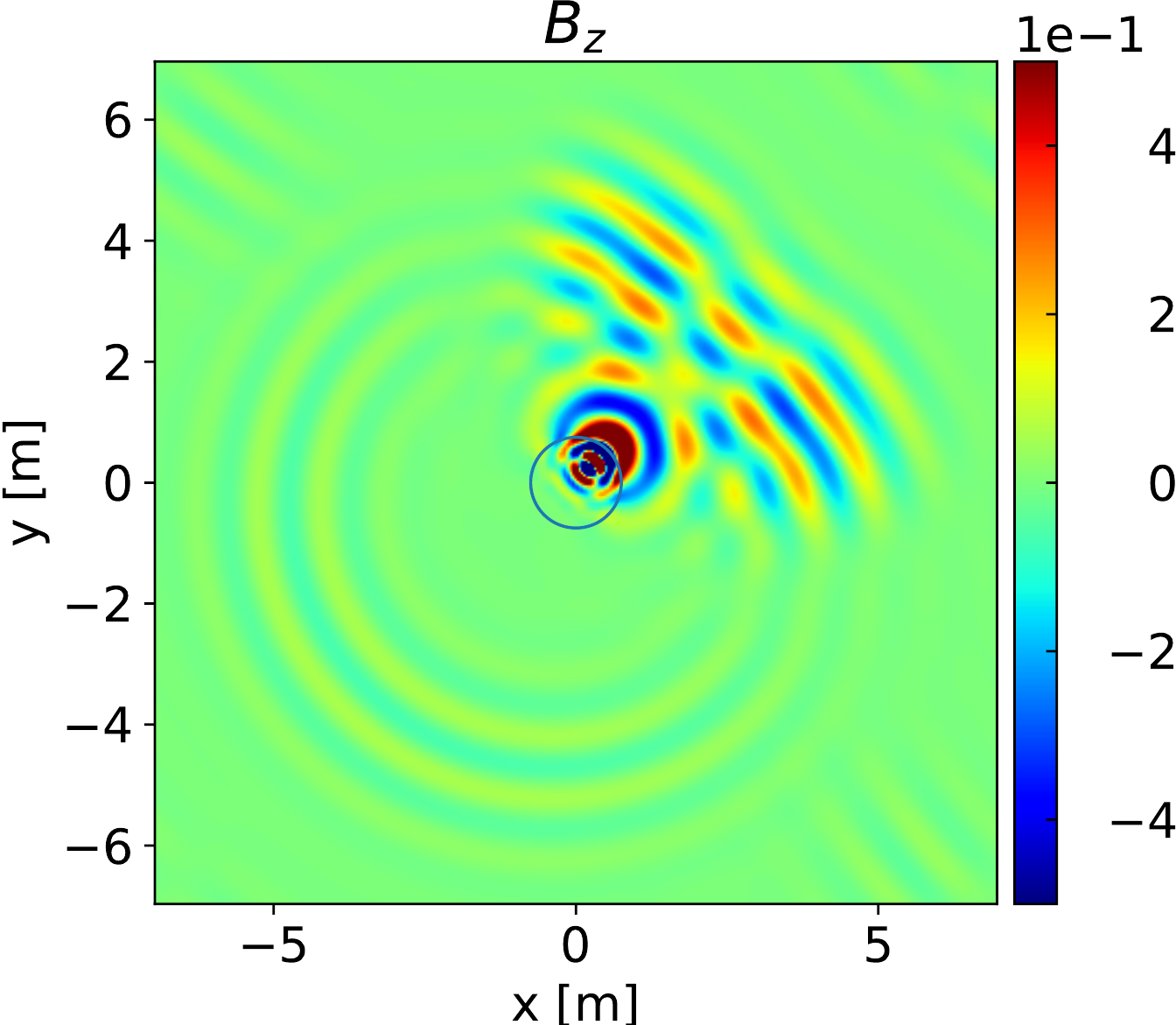} &
\includegraphics[width=0.33\textwidth]{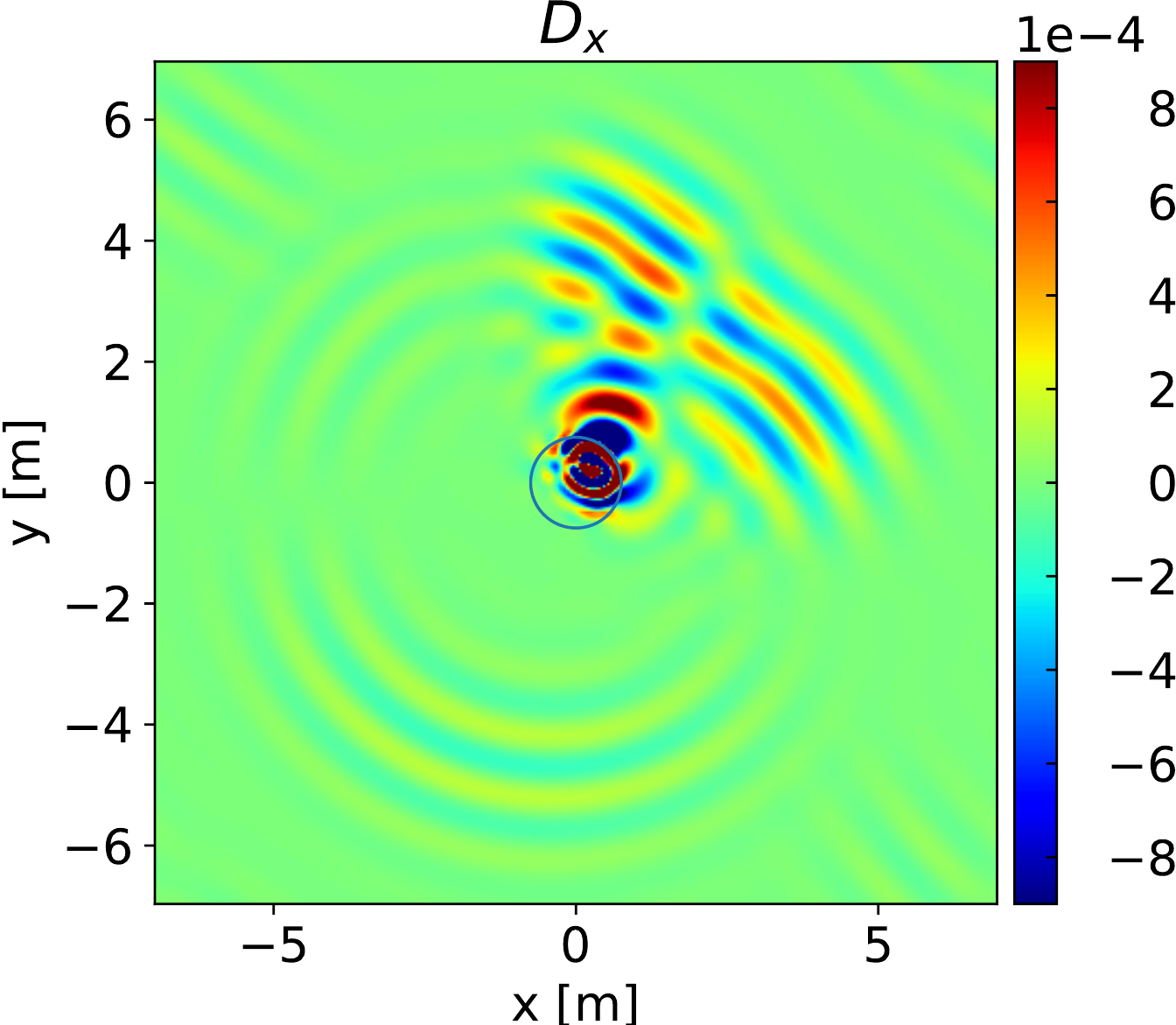} &
\includegraphics[width=0.33\textwidth]{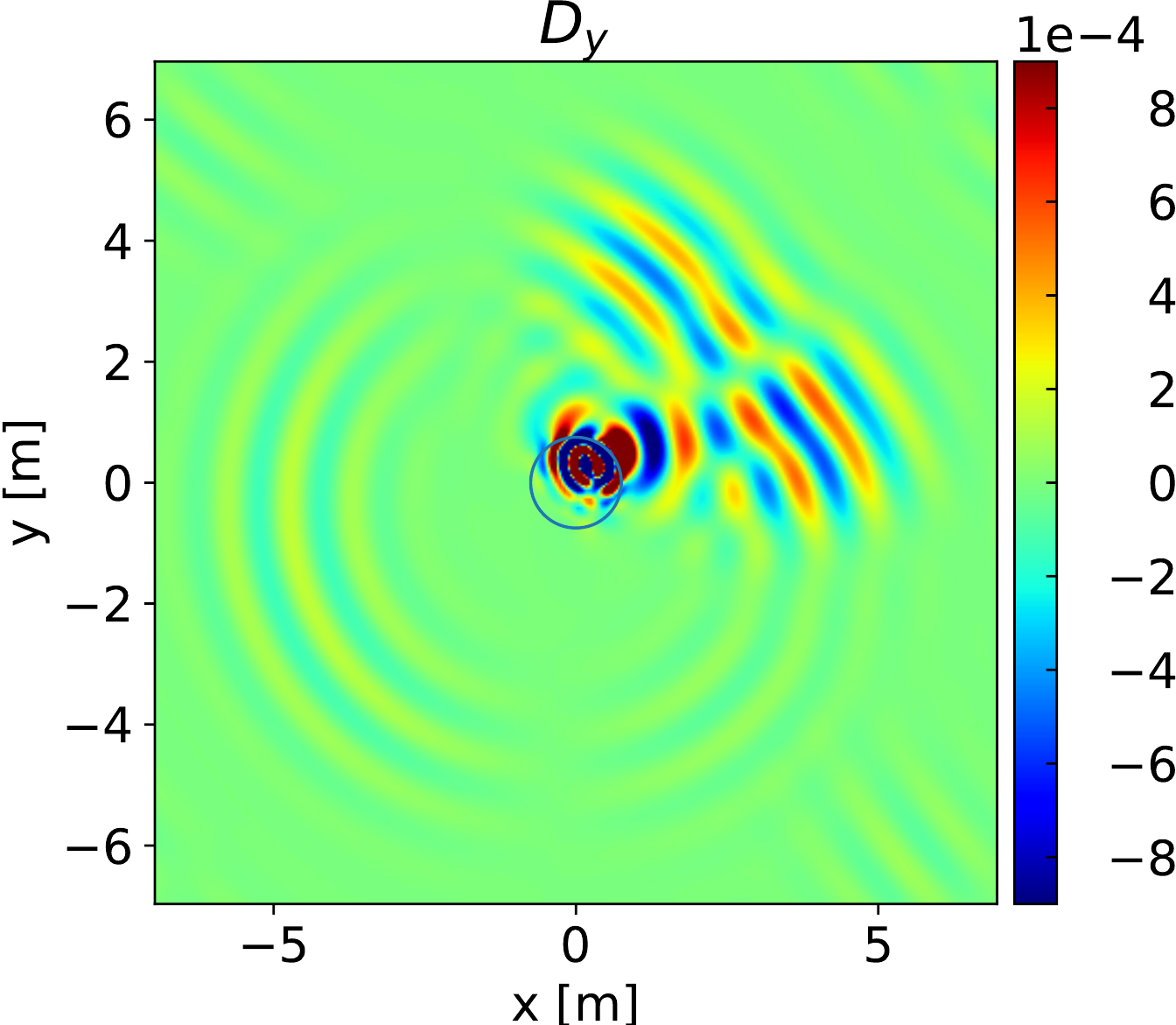} \\
\includegraphics[width=0.33\textwidth]{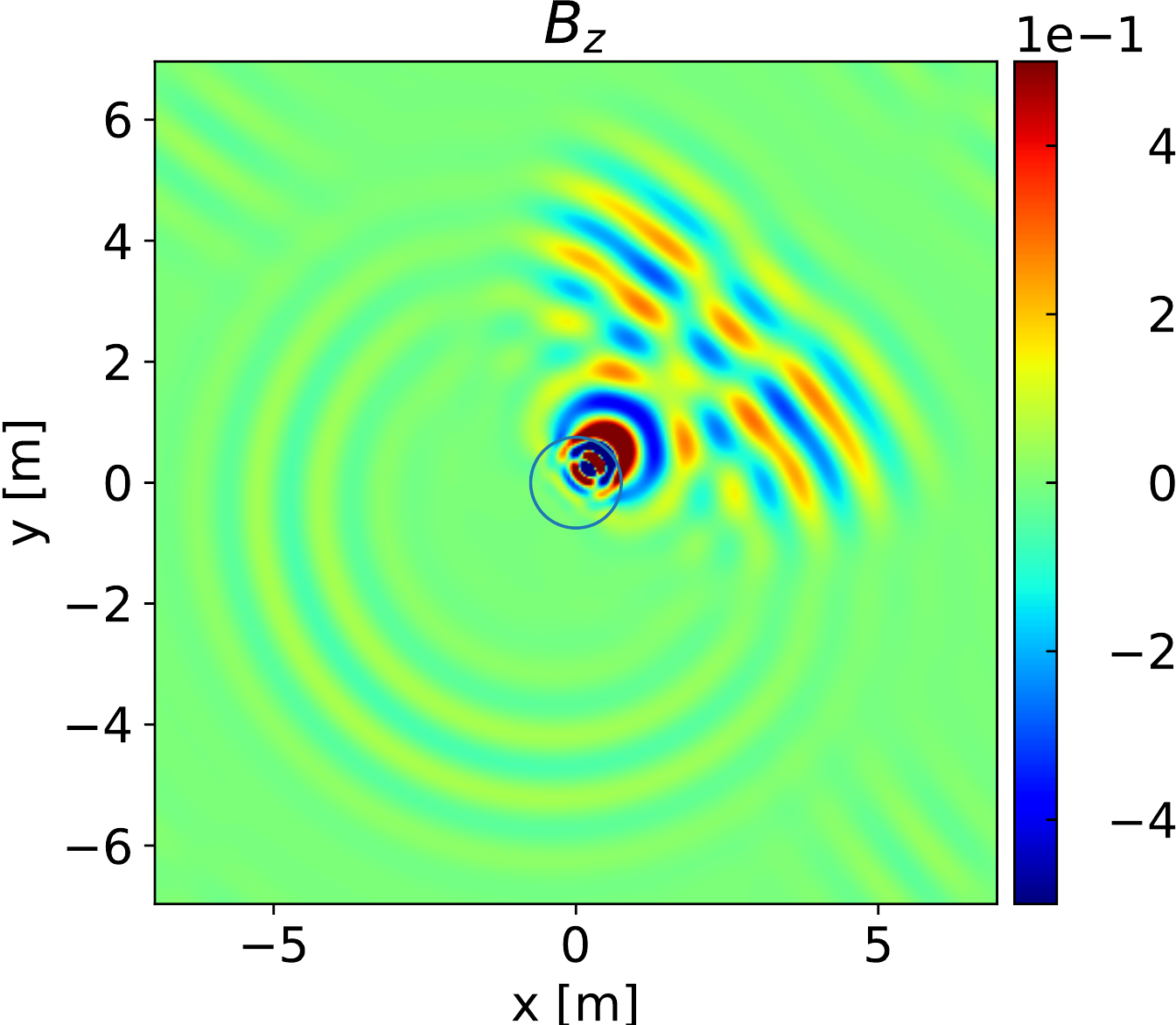} &
\includegraphics[width=0.33\textwidth]{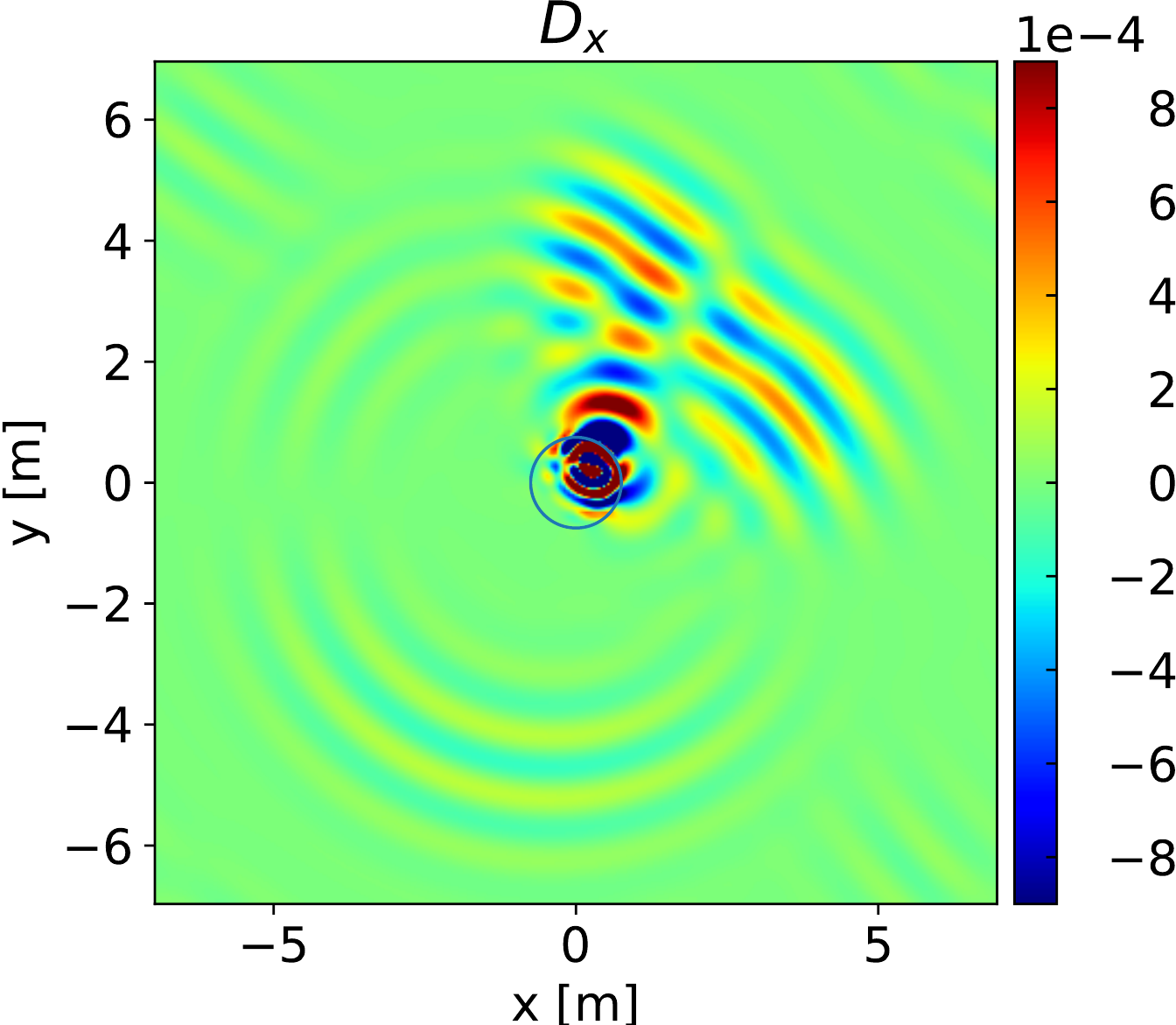} &
\includegraphics[width=0.33\textwidth]{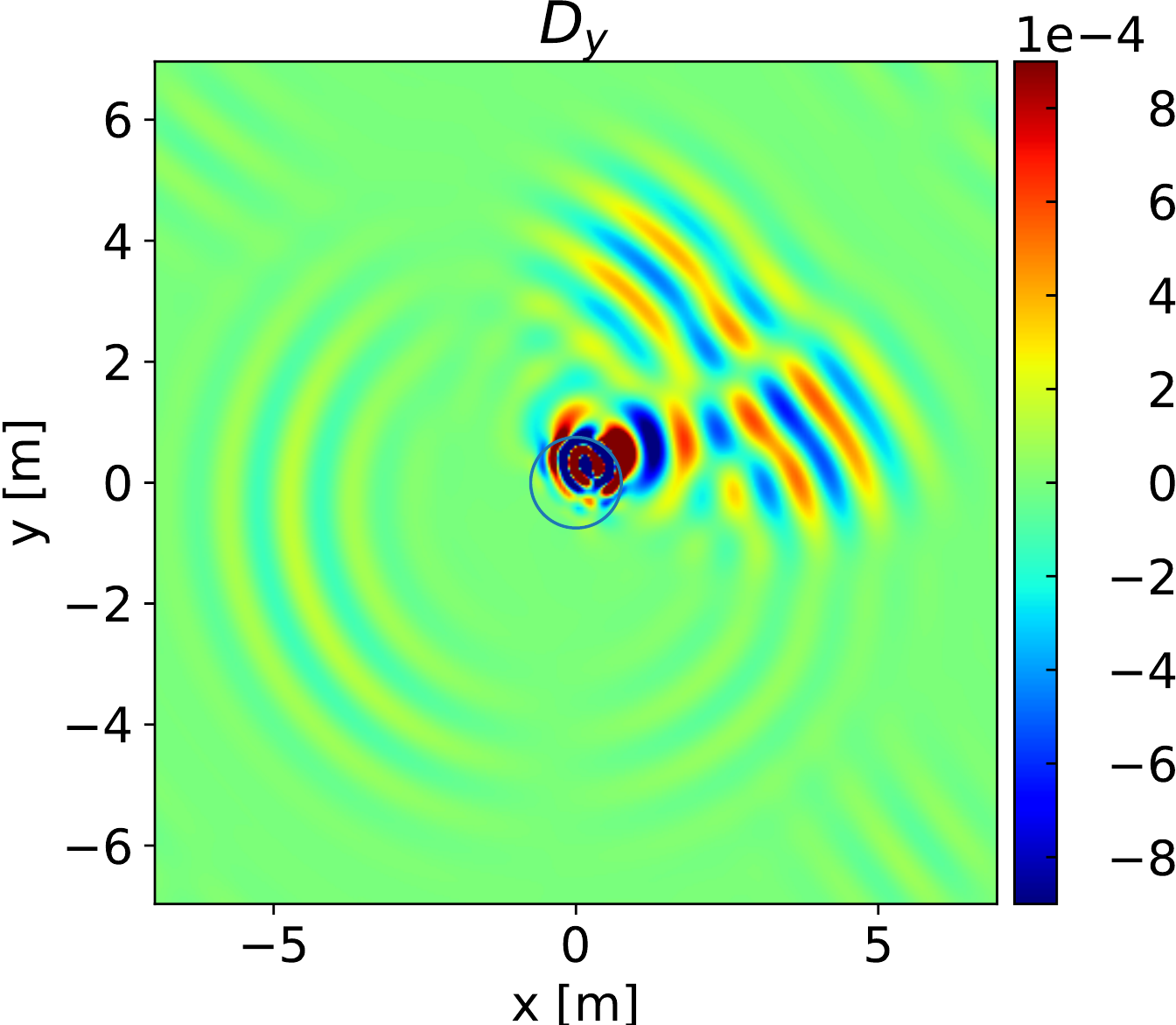}
\end{tabular}
\caption{Compact Gaussian electromagnetic pulse incident on a refractive disk using $200 \times 200$ cells. Top row: initial condition, middle row: $k=3$, bottom row: $k=4$}
\label{fig:gpulsecont}
\end{center}
\end{figure}


\begin{figure}
\begin{center}
\includegraphics[width=0.49\textwidth]{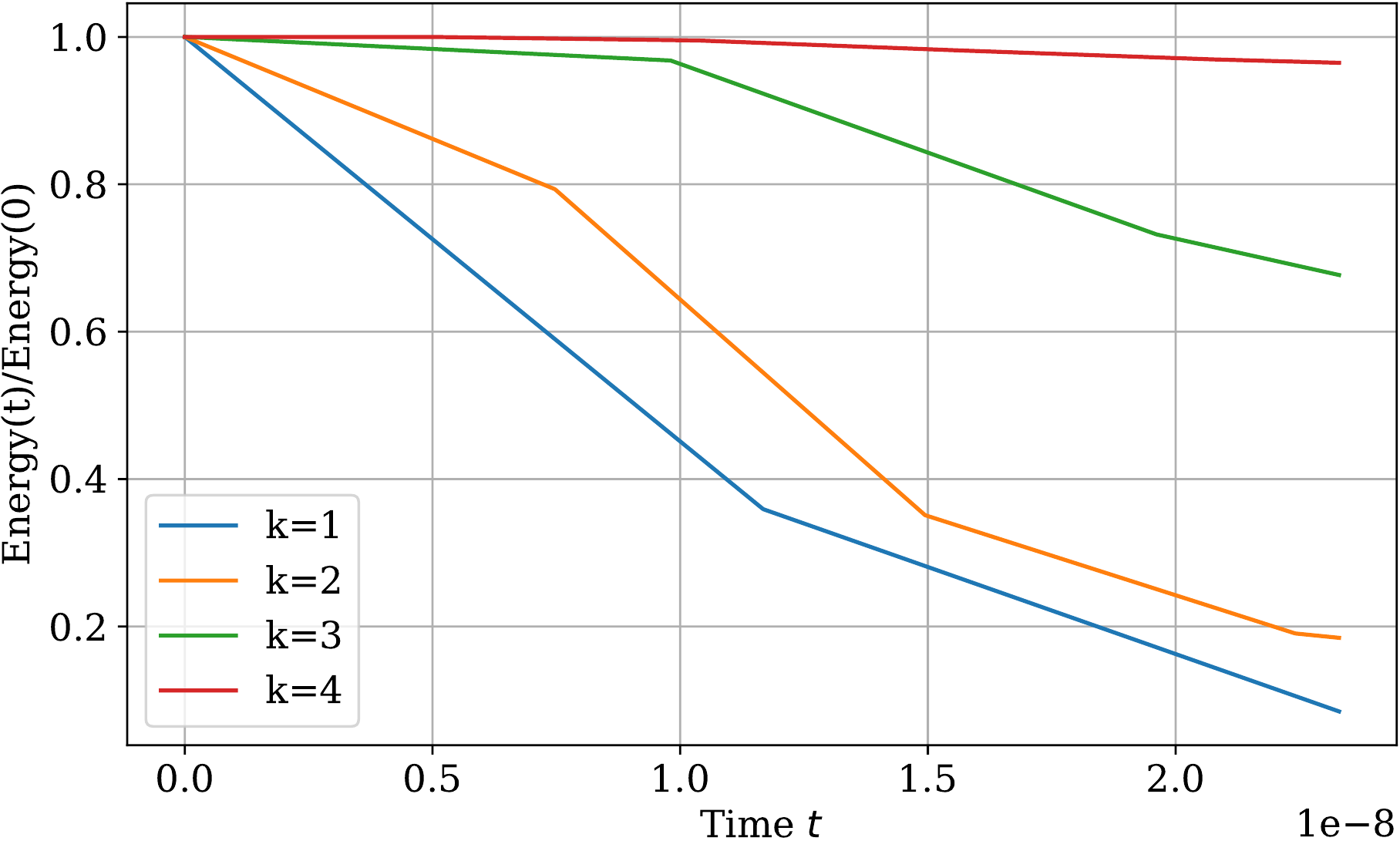} 
\includegraphics[width=0.49\textwidth]{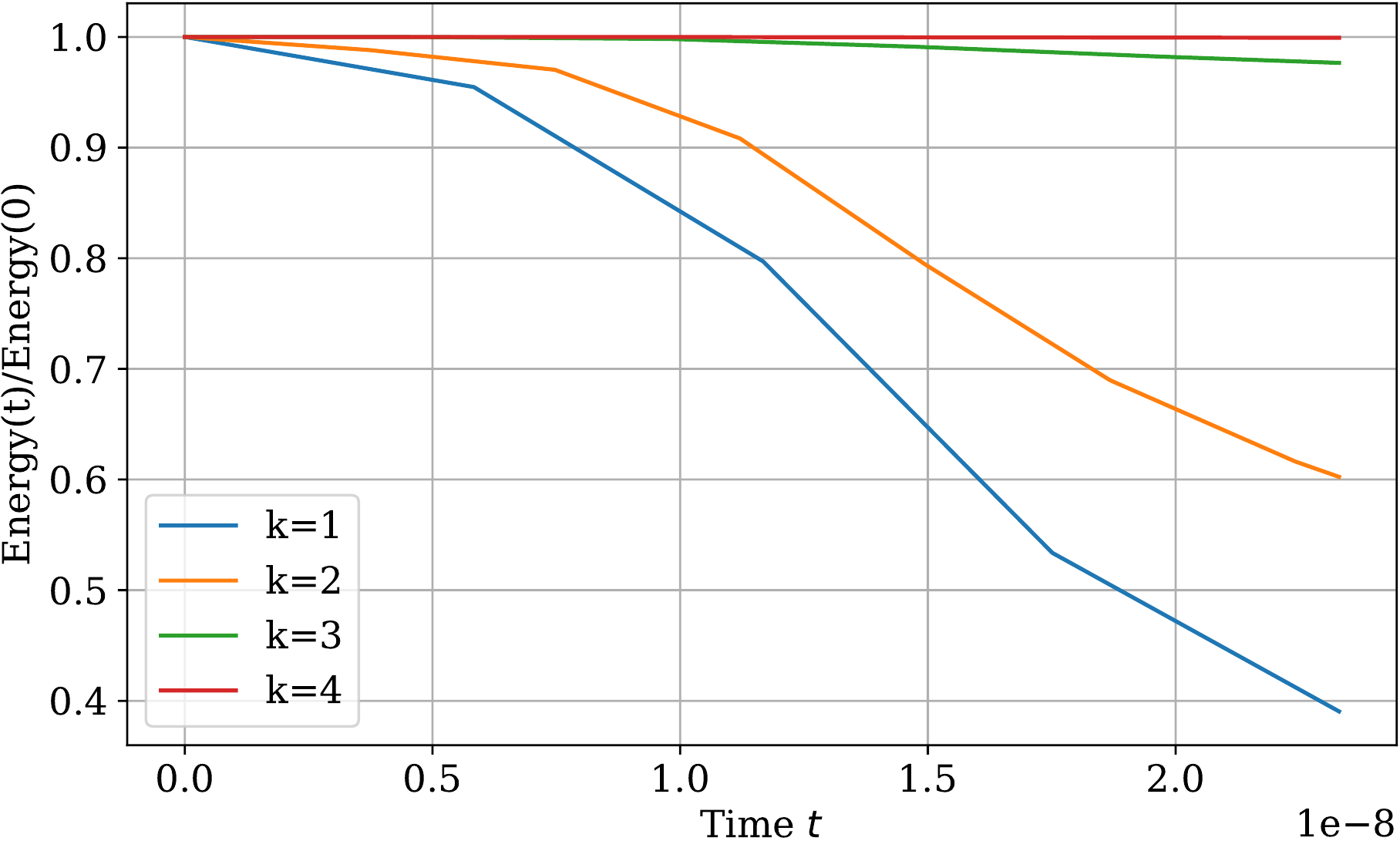} \\
(a) \hspace{6cm} (b) \\
\includegraphics[width=0.49\textwidth]{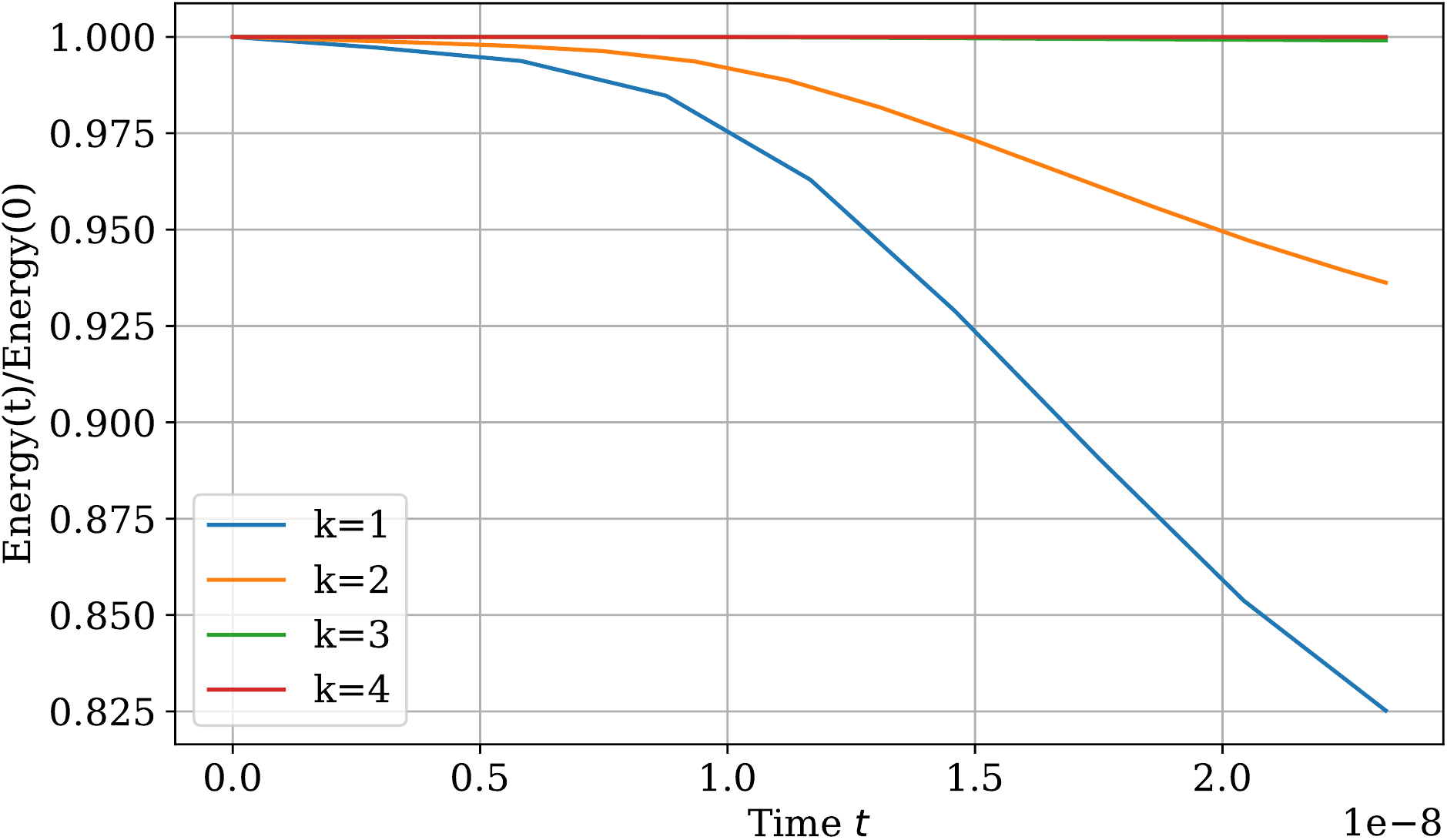} 
\includegraphics[width=0.49\textwidth]{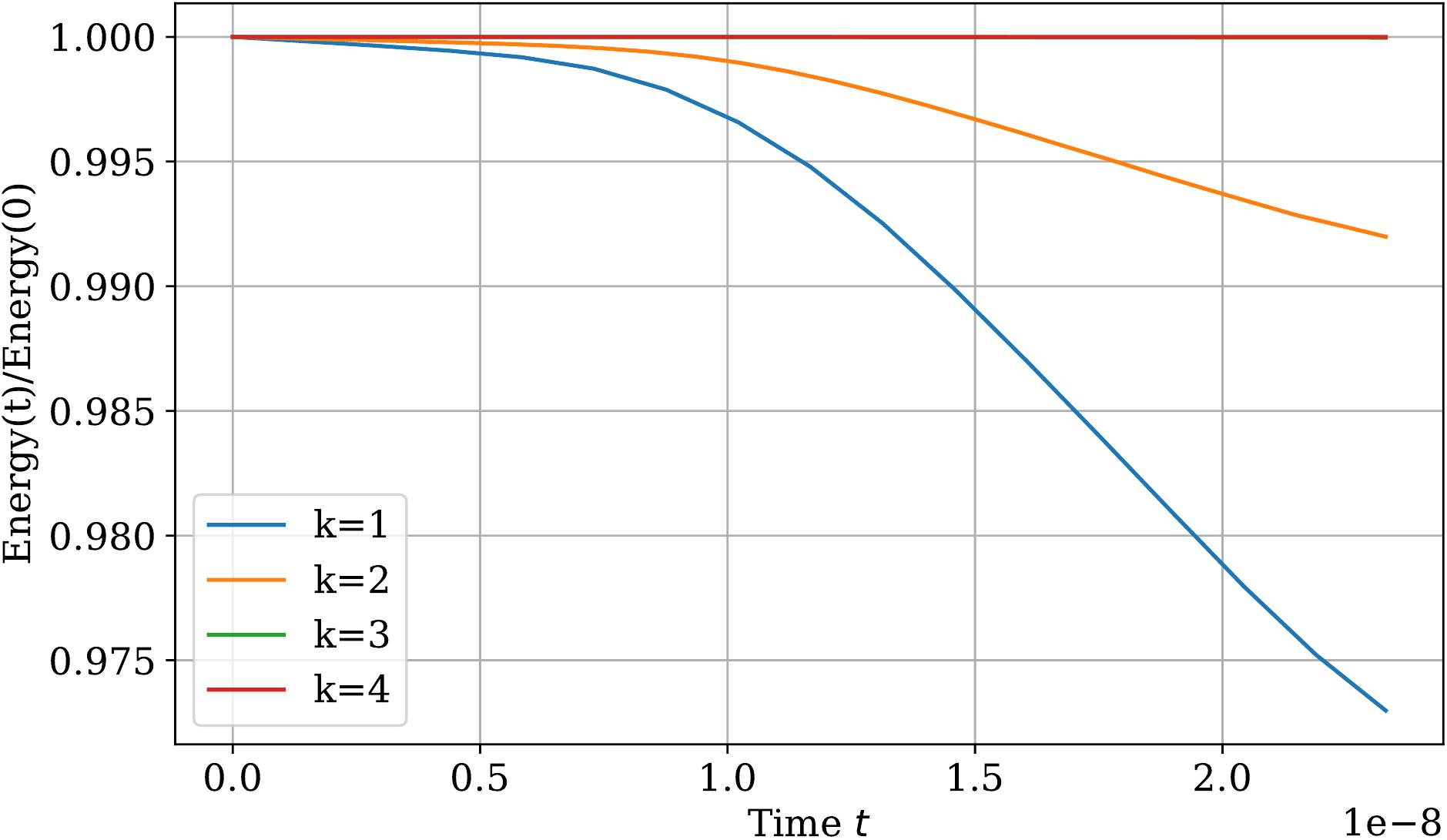} \\
(c) \hspace{6cm} (d) \\
\caption{Compact Gaussian electromagnetic pulse incident on a refractive disk. Evolution of total energy as a function of time. (a) $100 \times 100$ mesh, (b) $200 \times 200$ mesh, (c) $400 \times 400$, and (d) $800 \times 800$ mesh. The legends indicate the polynomial degree.}
\label{fig:gpulseE1}
\end{center}
\end{figure}


\begin{figure}
\begin{center}
\begin{tabular}{cc}
\includegraphics[width=0.49\textwidth]{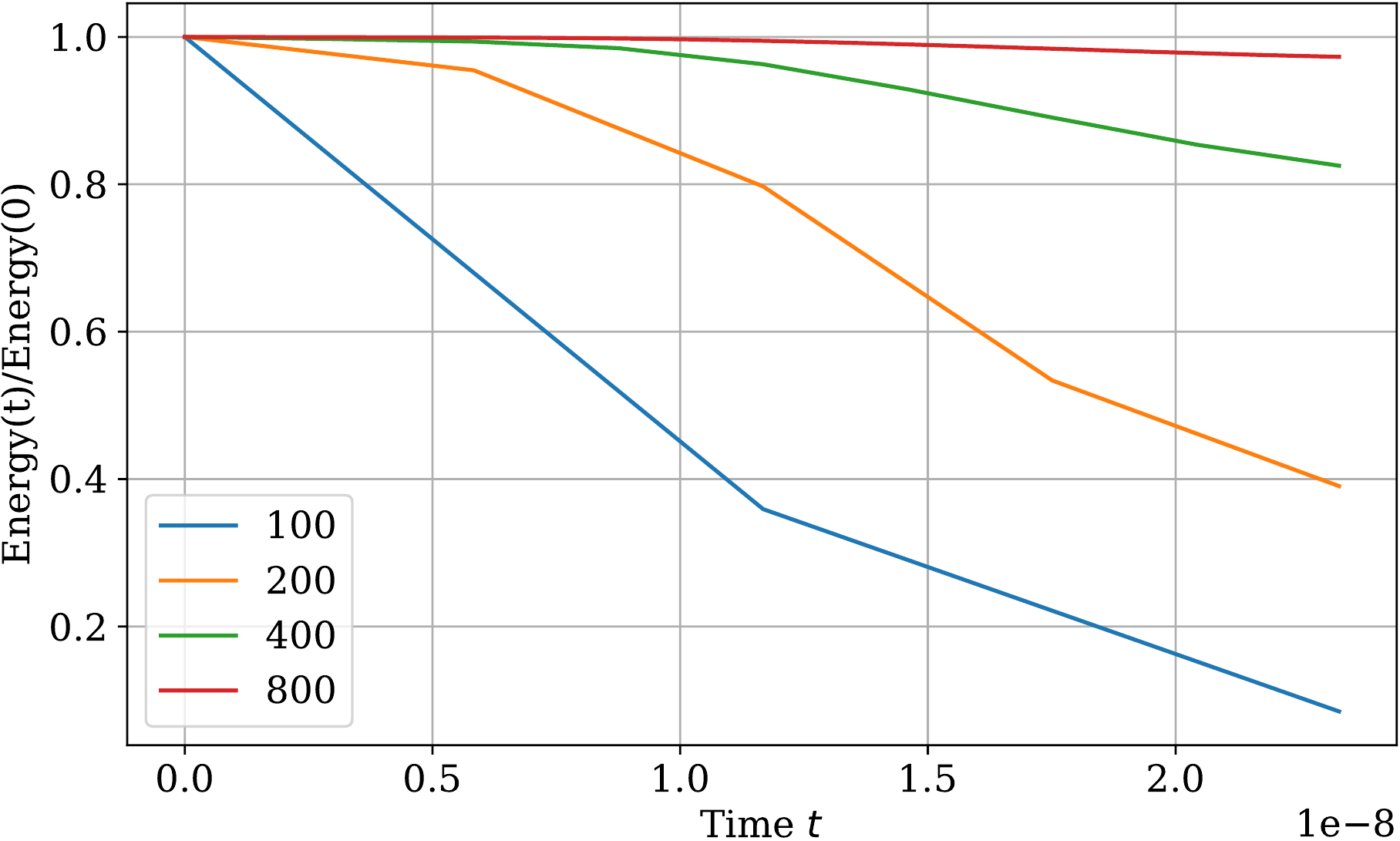} &
\includegraphics[width=0.49\textwidth]{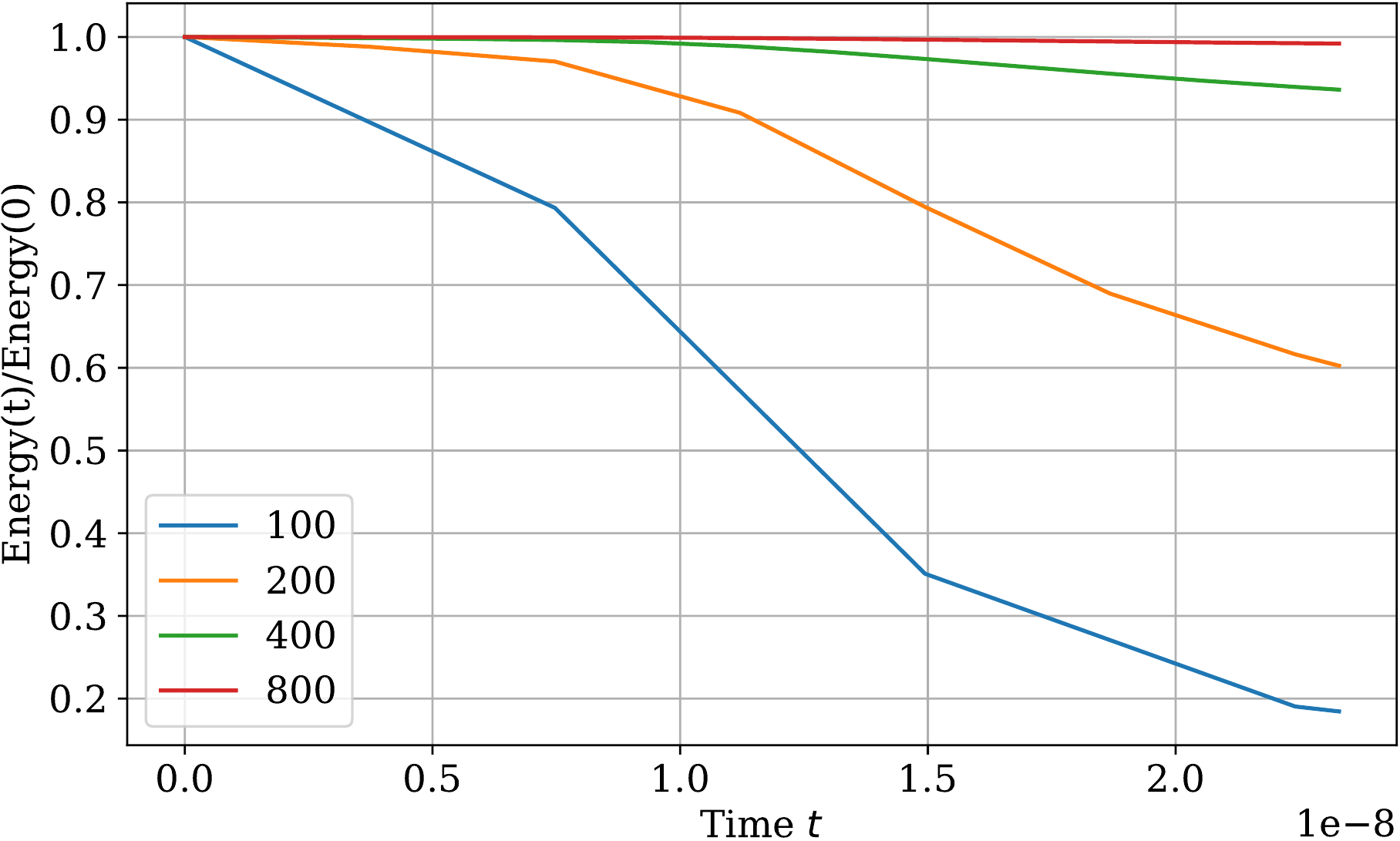} \\
(a) & (b) \\
\includegraphics[width=0.49\textwidth]{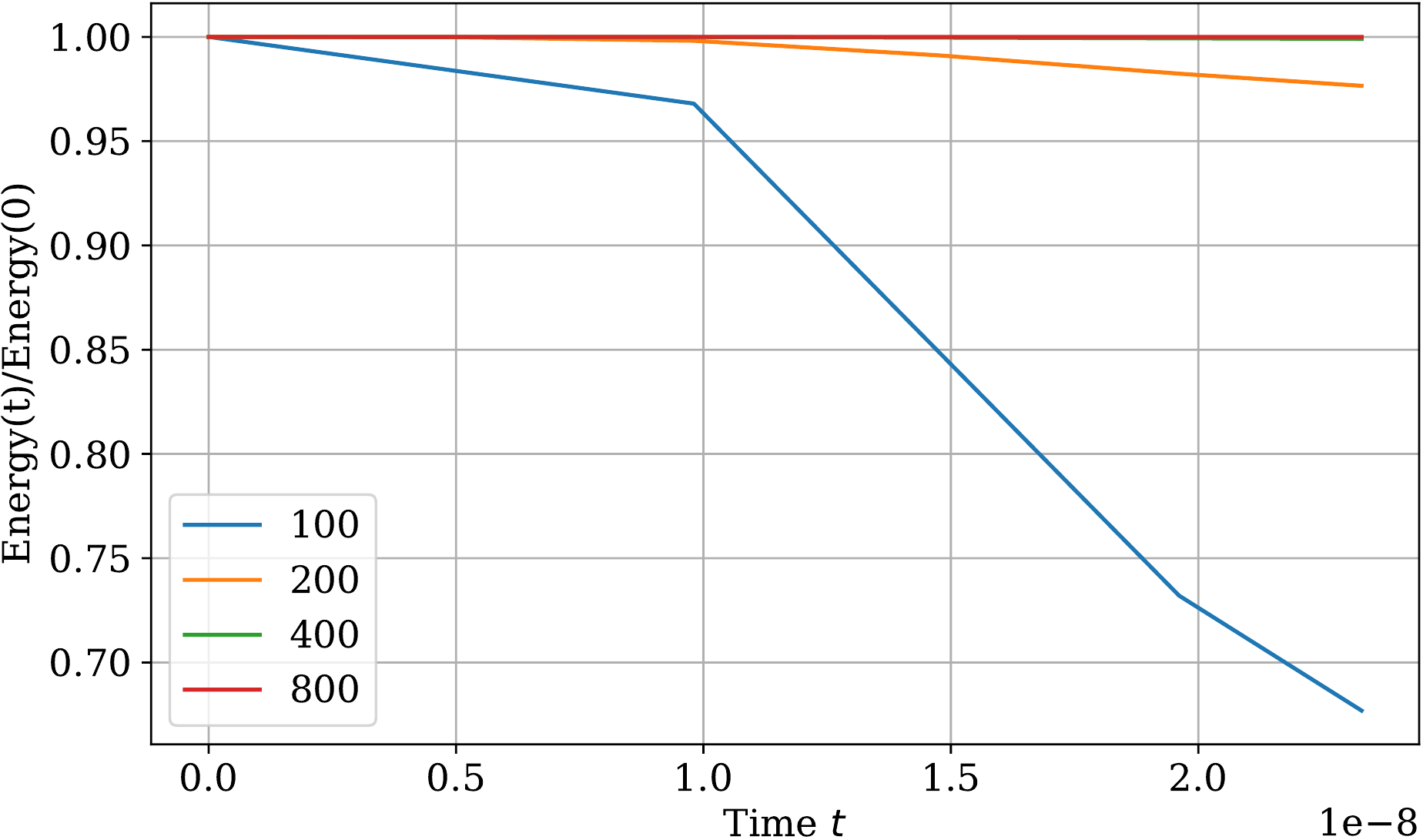} &
\includegraphics[width=0.49\textwidth]{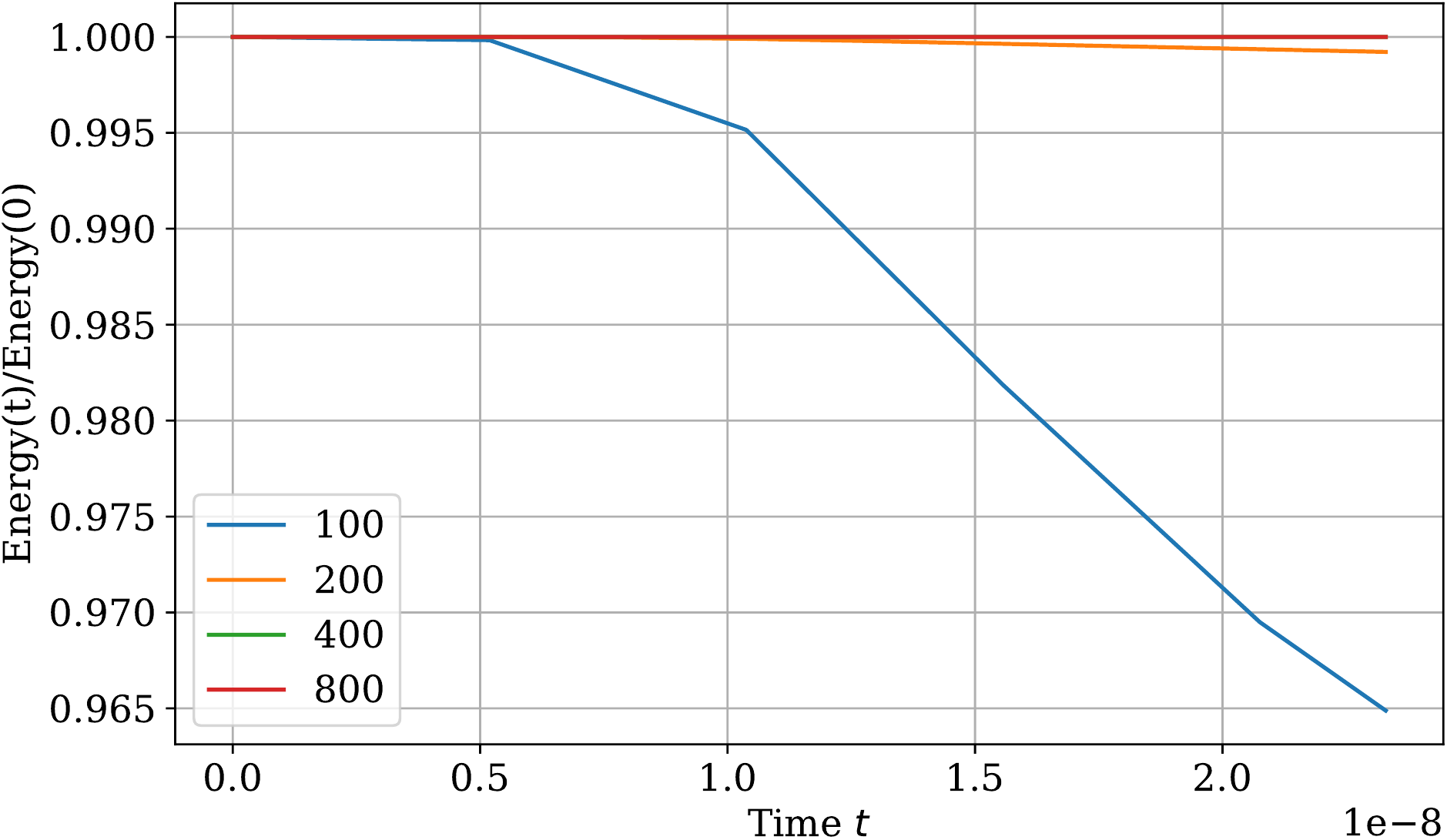} \\
(c) & (d)
\end{tabular}
\caption{Compact Gaussian electromagnetic pulse incident on a refractive disk. Evolution of total energy as a function of time for degree (a) $k=1$, (b) $k=2$, (a) $k=3$, (b) $k=4$. The legends indicate the mesh sizes.}
\label{fig:gpulseE2}
\end{center}
\end{figure}
In Section 6 we focused on the conservation of electromagnetic energy, showing that the DG schemes for CED are indeed energy stable. The DG methods presented here do not seem to need non-linear limiters and that is a very desirable trait for this class of schemes. It is known that when DG methods can operate without invoking limiters, the higher order DG schemes can come close to the ideal solution of the PDE. In practice, the use of Riemann solvers introduces some stabilization and some dissipation but in higher order DG schemes this dissipation is almost minimal.
	The conservation of electromagnetic energy, which depends quadratically on the primal variables that are evolved, is not guaranteed in numerical CED schemes. This is true for FDTD, FVTD and also for the DGTD schemes for CED designed here. In Balsara and Kappeli [14] we showed that the electromagnetic energy is, nevertheless, conserved extremely well by higher order DGTD schemes. However, Balsara and Kappeli showed this for electromagnetic radiation that is propagating in a vacuum. It is, therefore, interesting to try and quantify how well the electromagnetic energy is conserved on the mesh when electromagnetic radiation interacts with spatially varying dielectric properties in materials. To make this demonstration, we solve the problem of a compact Gaussian electromagnetic pulse incident on a dielectric disk on meshes with $100\times 100$, $200\times 200$, $400\times 400$ and $800\times 800$ zones. Some sample results are shown in Fig.~(\ref{fig:gpulsecont}) on a mesh of $200 \times 200$ cells using the fourth and fifth order schemes. We then plot the electromagnetic energy as a function of time. Because periodic boundary conditions were used, we hope that the more accurate DGTD schemes will conserve electromagnetic energy. 
	Fig.~(\ref{fig:gpulseE1}a) shows the energy evolution as a function of time on a $100\times 100$ zone mesh from second, third, fourth and fifth order DGTD schemes. From Fig.~(\ref{fig:gpulseE1}a) we see that only the fifth order scheme does a superlative job of energy conservation, with the fourth order scheme performing very well. Fig.~(\ref{fig:gpulseE1}b) shows the same information as Fig.~(\ref{fig:gpulseE1}a), but this time on a $200\times 200$ zone mesh. We see now that both the fourth and fifth order schemes show superlative energy conservation. Fig.~(\ref{fig:gpulseE1}c) shows the same information as Figs.~(\ref{fig:gpulseE1}a) and (\ref{fig:gpulseE1}b), but this time on a $400\times 400$ zone mesh. We now see that even the third order DGTD scheme has very good energy conservation properties. This trend continues in Fig.~(\ref{fig:gpulseE1}d) which shows the energy plots for the $800 \times 800$ mesh.
	Because Fig.~(\ref{fig:gpulseE1}) shows all the data, including the data from the second order DGTD scheme, it is not possible to fully appreciate how well the higher order DGTD schemes conserve energy. For that reason, Fig.~(\ref{fig:gpulseE2}) shows the energy evolution as a function of time for the second, third, fourth and fifth order schemes when a sequence of mesh resolutions are used for the same problem. The vertical scale in Figs.~(\ref{fig:gpulseE2}) show the extraordinarily good ability of the higher order DGTD schemes to conserve electromagnetic energy.

Tables~(\ref{tab:gpulsek1finalNONE})-(\ref{tab:gpulsek4finalNONE}) show the accuracy of the second, third, fourth and fifth order schemes for the Gaussian pulse problem. We see that the schemes reach their designed accuracies on relatively coarse meshes; which shows that our DG formulation is not just asymptotically very accurate but it also offers an accuracy advantage on poorly resolved meshes. This is especially significant because we allowed for an order of magnitude variation in the permittivity and did nothing special to treat that variation in permittivity. By comparing tables~(\ref{tab:gpulsek3finalNONE}) and (\ref{tab:gpulsek4finalNONE}) to tables (\ref{tab:gpulsek1finalNONE}) and (\ref{tab:gpulsek2finalNONE}) we see that the highest order schemes have reached their design accuracies on the coarsest meshes, which brings out another importance of very high order, globally constraint-preserving DG schemes for CED.

 \begin{table}
 \begin{center} 
 \begin{tabular}{|c|c|c|c|c|c|c|c|c|}
 \hline
 $N_x\times N_y$ & $\|\D^h-\D\|_{L^1}$ & Ord & $\|\D^h-\D\|_{L^2}$ &Ord & $\|B_z^h-B_z\|_{L^1}$ & Ord & $\|B_z^h-B_z\|_{L^2}$ & Ord\\ 
 \hline
$100\times 100$ & 8.7117e-05 & --- & 5.7716e-04&---&2.4173e-02&---&9.1668e-02&--- \\
$200\times 200$ & 4.1169e-05 & 1.08&3.6468e-04&0.66&1.0302e-02&1.23&5.5462e-02&0.72\\
$400\times 400$ & 8.9100e-06 & 2.21&8.6861e-05&2.07&2.1341e-03&2.27&1.2224e-02&2.18\\
$800\times 800$ & 1.2956e-06 & 2.78&1.1942e-05&2.86&3.2063e-04&2.73&1.6685e-03&2.87\\
 \hline 
 \end{tabular} 
 \end{center} 
 \caption{Compact Gaussian electromagnetic pulse incident on a refractive disk. Error convergence for degree $k=1$ under mesh refinement.}
 \label{tab:gpulsek1finalNONE} 
 \end{table}

 \begin{table}
 \begin{center} 
 \begin{tabular}{|c|c|c|c|c|c|c|c|c|}
 \hline
 $N_x\times N_y$ & $\|\D^h-\D\|_{L^1}$ & Ord & $\|\D^h-\D\|_{L^2}$ &Ord & $\|B_z^h-B_z\|_{L^1}$ & Ord & $\|B_z^h-B_z\|_{L^2}$ & Ord\\ 
 \hline
$100\times 100$ & 6.1017e-05 & --- & 4.9501e-04&---&1.5650e-02&---&7.8412e-02&--- \\
$200\times 200$ & 1.9119e-05 & 1.67&2.0730e-04&1.26&4.4152e-03&1.83&3.5235e-02&1.15\\
$400\times 400$ & 2.5029e-06 & 2.93&2.6761e-05&2.95&5.9512e-04&2.89&4.8243e-03&2.87\\
$800\times 800$ & 2.7078e-07 & 3.21&2.8929e-06&3.21&6.4388e-05&3.21&5.1386e-04&3.23\\
 \hline 
 \end{tabular} 
 \end{center} 
\caption{Compact Gaussian electromagnetic pulse incident on a refractive disk. Error convergence for degree $k=2$ under mesh refinement.}           
 \label{tab:gpulsek2finalNONE} 
 \end{table}

 \begin{table}
 \begin{center} 
 \begin{tabular}{|c|c|c|c|c|c|c|c|c|}
 \hline
 $N_x\times N_y$ & $\|\D^h-\D\|_{L^1}$ & Ord & $\|\D^h-\D\|_{L^2}$ &Ord & $\|B_z^h-B_z\|_{L^1}$ & Ord & $\|B_z^h-B_z\|_{L^2}$ & Ord\\ 
 \hline
$100\times 100$ & 2.1628e-05 & --- & 2.0508e-04&---&5.1752e-03&---&3.2412e-02&--- \\
$200\times 200$ & 1.0561e-06 & 4.36&1.2522e-05&4.03&2.2929e-04&4.50&2.0243e-03&4.00\\
$400\times 400$ & 4.3566e-08 & 4.60&4.7222e-07&4.73&1.0537e-05&4.44&8.8197e-05&4.52\\
$800\times 800$ & 1.8547e-09 & 4.55&1.6639e-08&4.83&5.3238e-07&4.31&4.4683e-06&4.30\\
 \hline 
 \end{tabular} 
 \end{center} 
\caption{Compact Gaussian electromagnetic pulse incident on a refractive disk. Error convergence for degree $k=3$ under mesh refinement.}
 \label{tab:gpulsek3finalNONE} 
 \end{table}

 \begin{table}
 \begin{center} 
 \begin{tabular}{|c|c|c|c|c|c|c|c|c|}
 \hline
 $N_x\times N_y$ & $\|\D^h-\D\|_{L^1}$ & Ord & $\|\D^h-\D\|_{L^2}$ &Ord & $\|B_z^h-B_z\|_{L^1}$ & Ord & $\|B_z^h-B_z\|_{L^2}$ & Ord\\ 
 \hline
$100\times 100$ & 4.0374e-06 & --- & 3.8018e-05&---&8.9893e-04&---&5.3665e-03&--- \\
$200\times 200$ & 6.5567e-08 & 5.94&6.2380e-07&5.93&2.1939e-05&5.36&2.7262e-04&4.30\\
$400\times 400$ & 2.4957e-09 & 4.72&1.7238e-08&5.18&9.1588e-07&4.58&7.8449e-06&5.12\\
 \hline 
 \end{tabular} 
 \end{center} 
\caption{Compact Gaussian electromagnetic pulse incident on a refractive disk. Error convergence for degree $k=4$ under mesh refinement.}       
 \label{tab:gpulsek4finalNONE} 
 \end{table}

\subsection{Refraction of a compact electromagnetic beam by a dielectric slab}
\begin{figure}
\begin{center}
\begin{tabular}{ccc}
\includegraphics[width=0.33\textwidth]{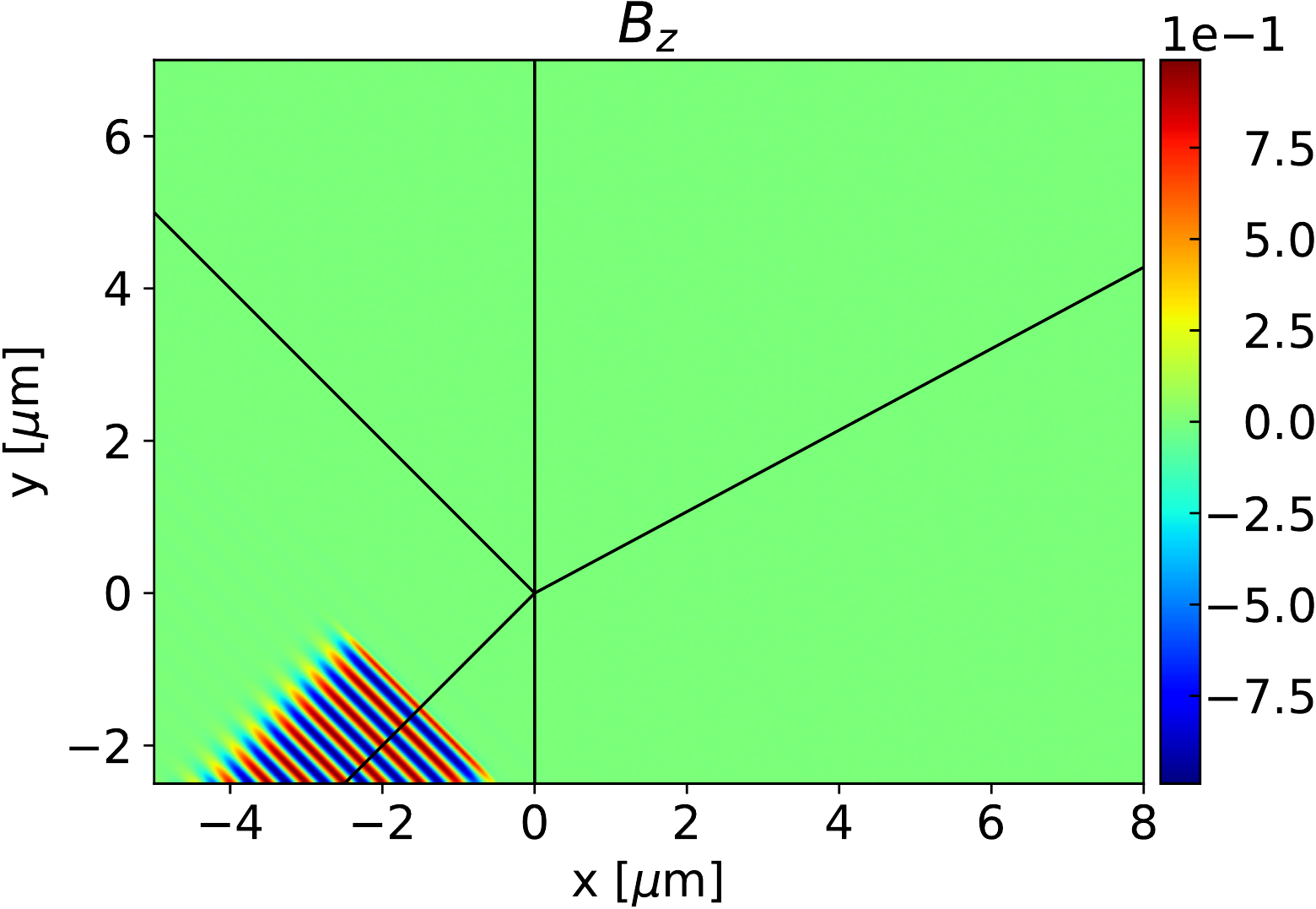} &
\includegraphics[width=0.33\textwidth]{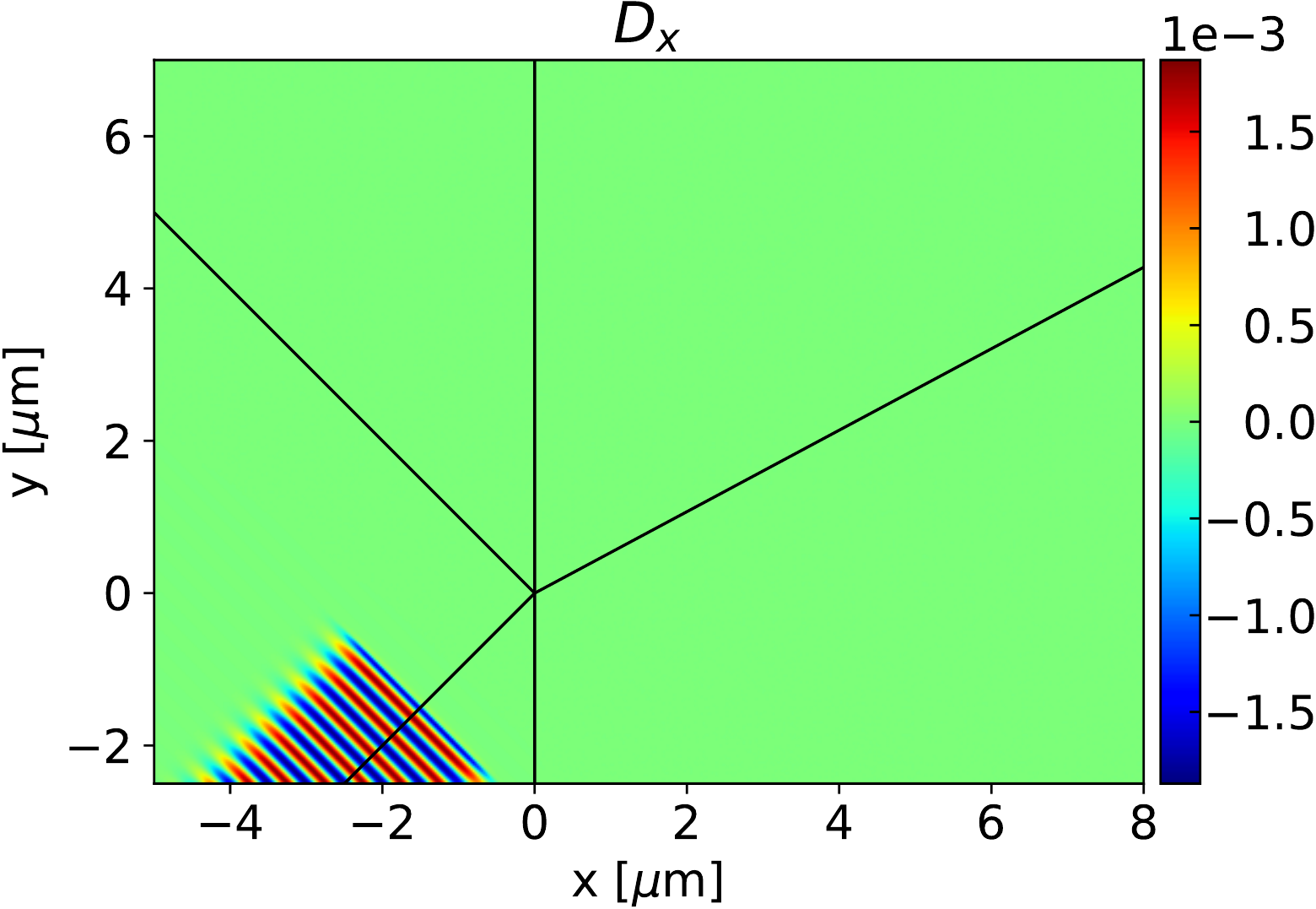} &
\includegraphics[width=0.33\textwidth]{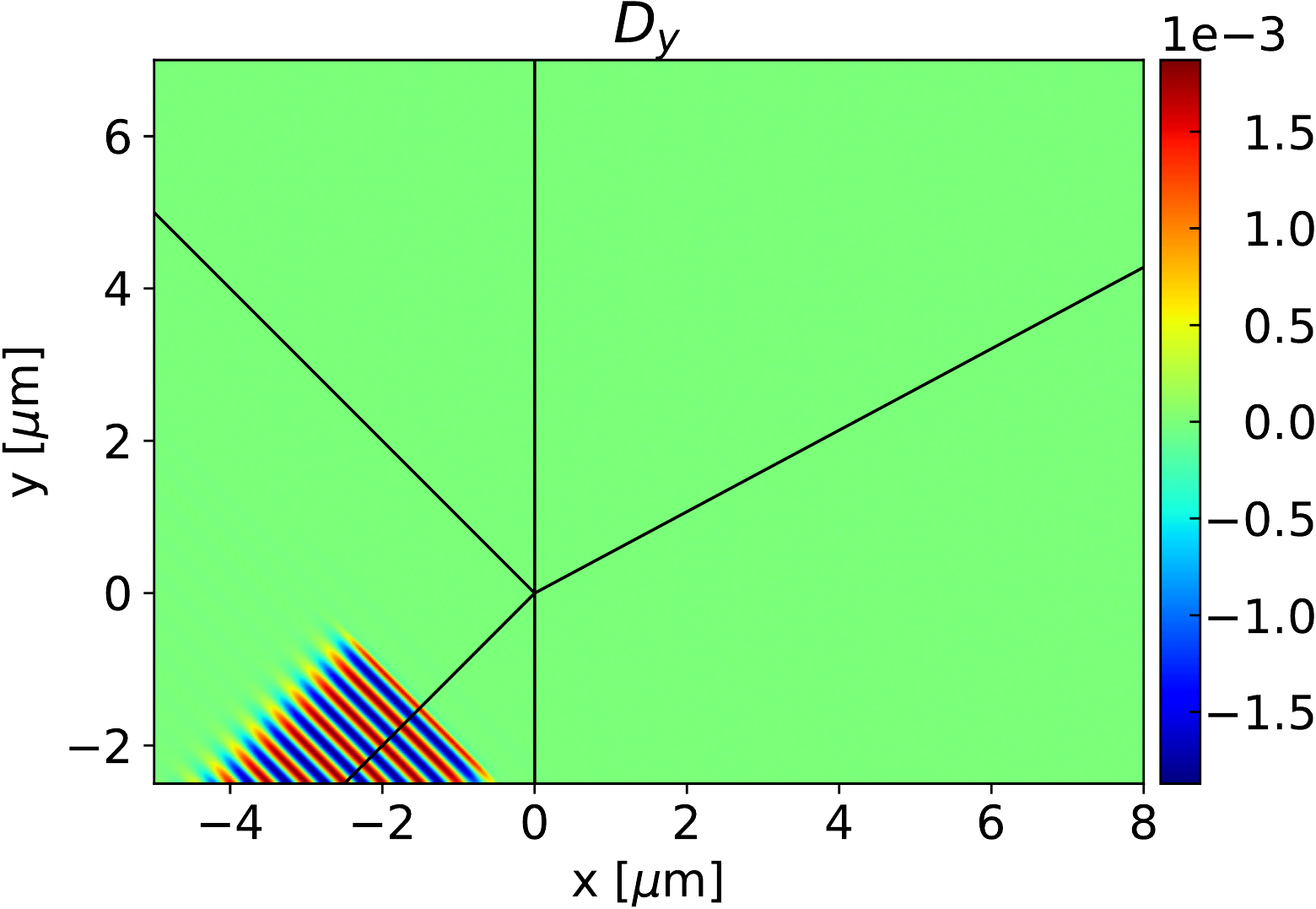} \\
\includegraphics[width=0.33\textwidth]{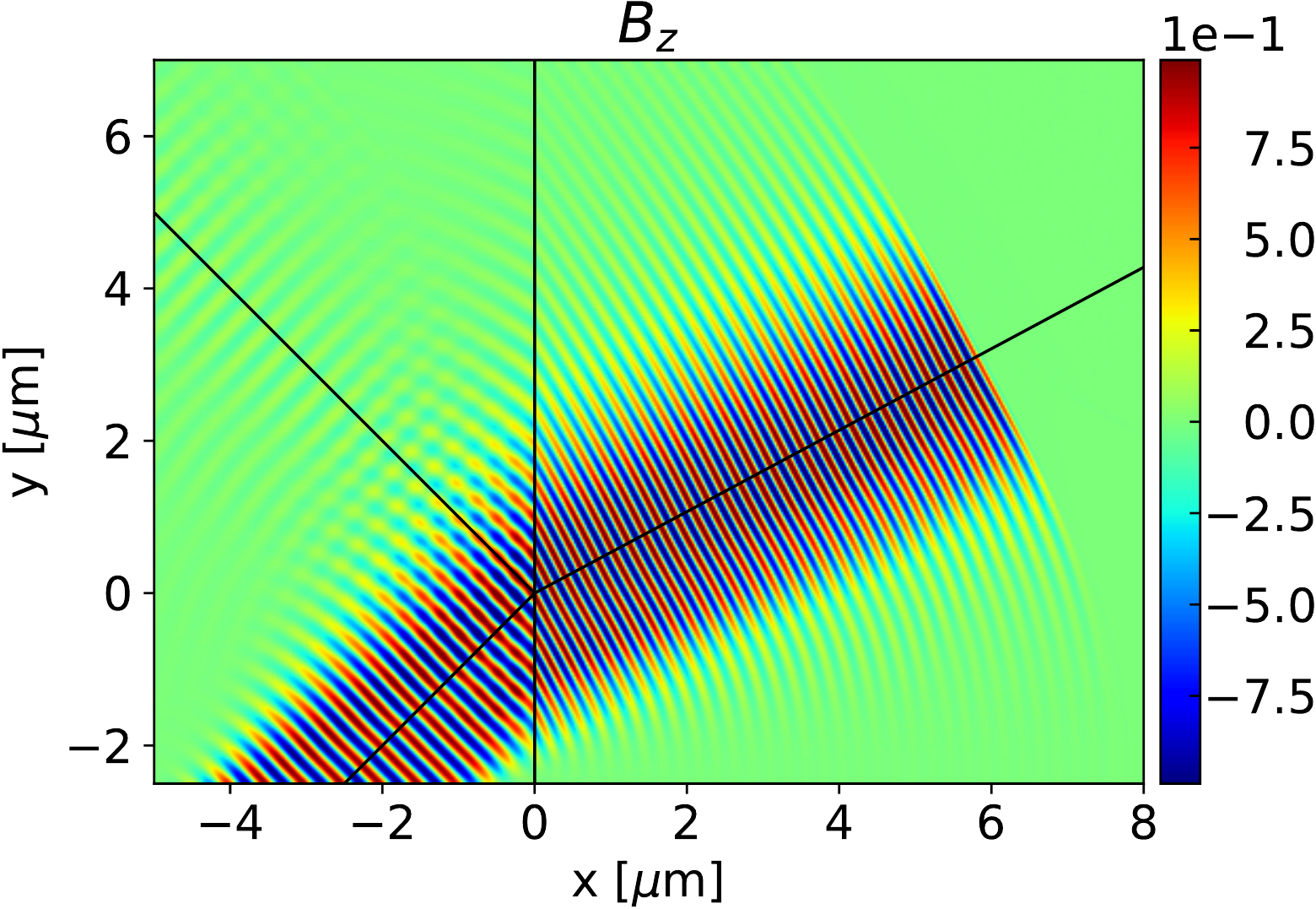} &
\includegraphics[width=0.33\textwidth]{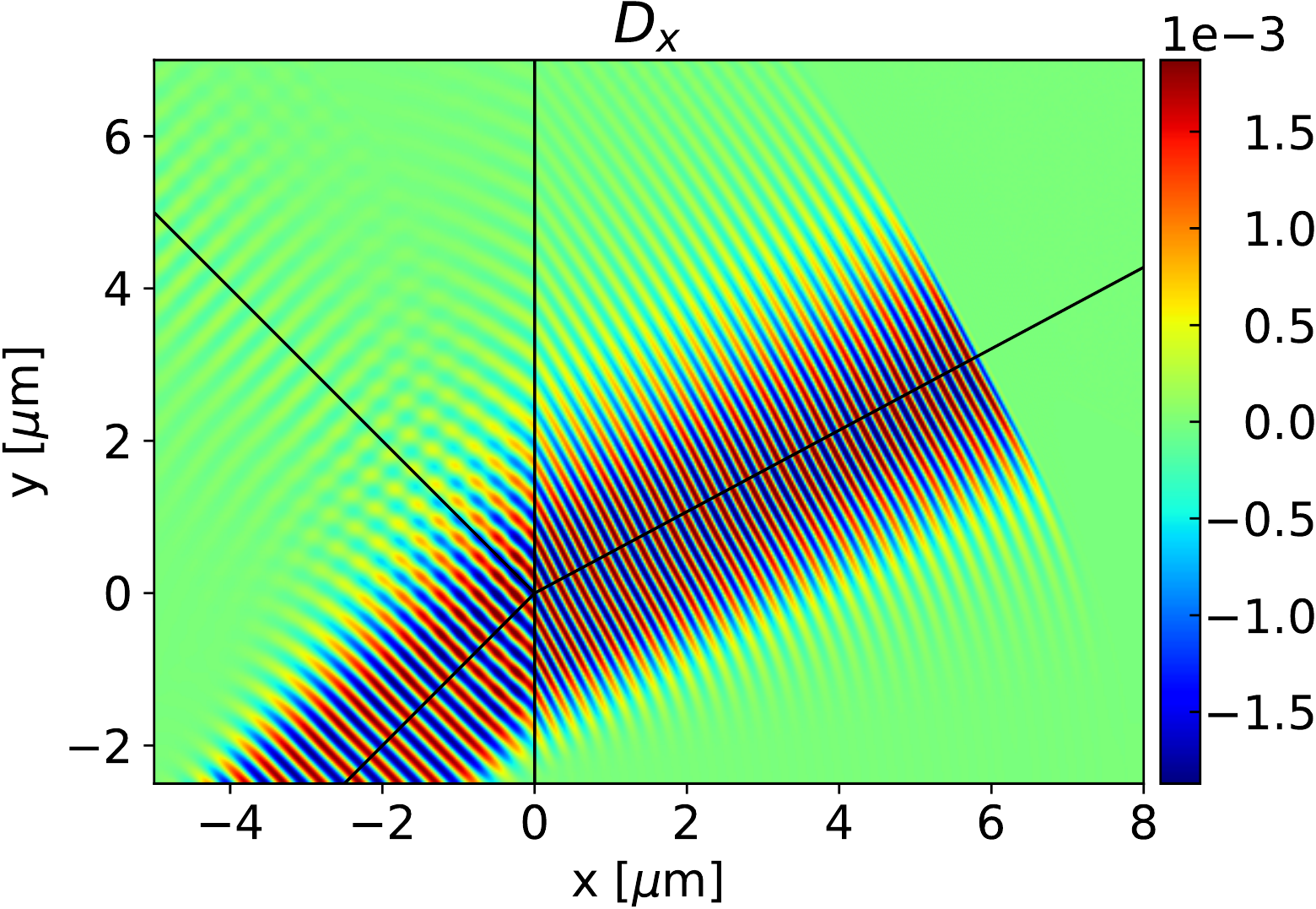} &
\includegraphics[width=0.33\textwidth]{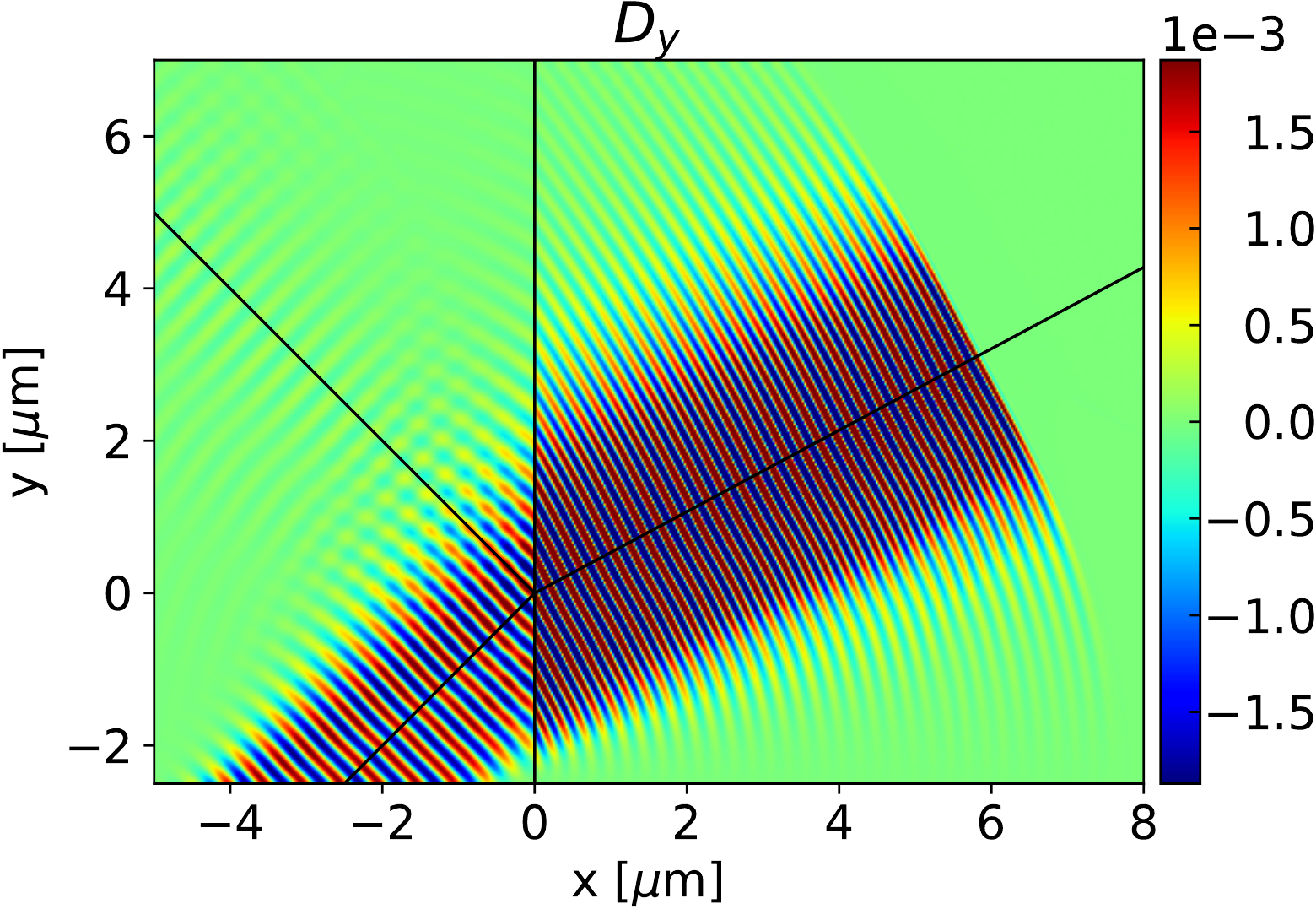} \\
\includegraphics[width=0.33\textwidth]{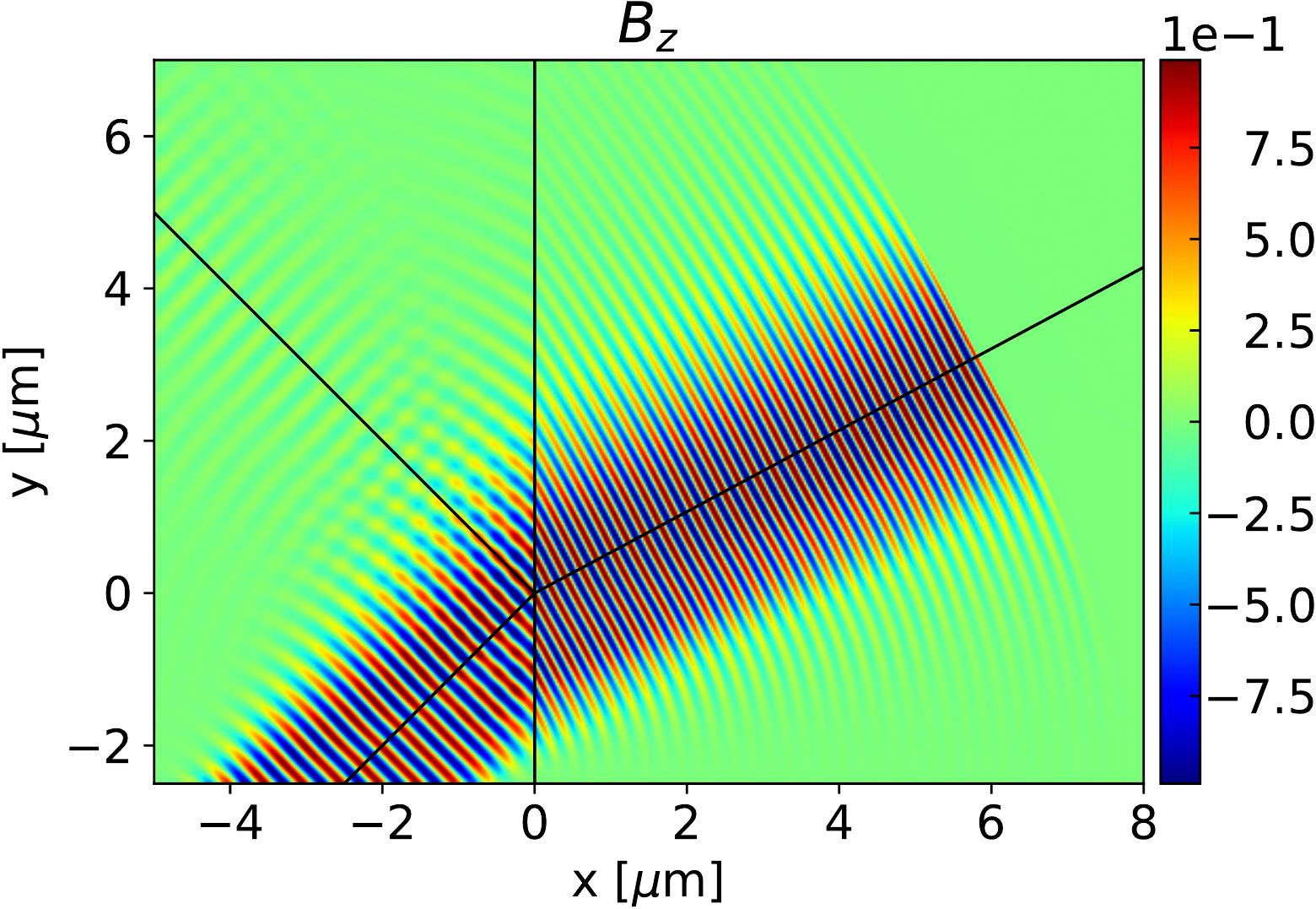} &
\includegraphics[width=0.33\textwidth]{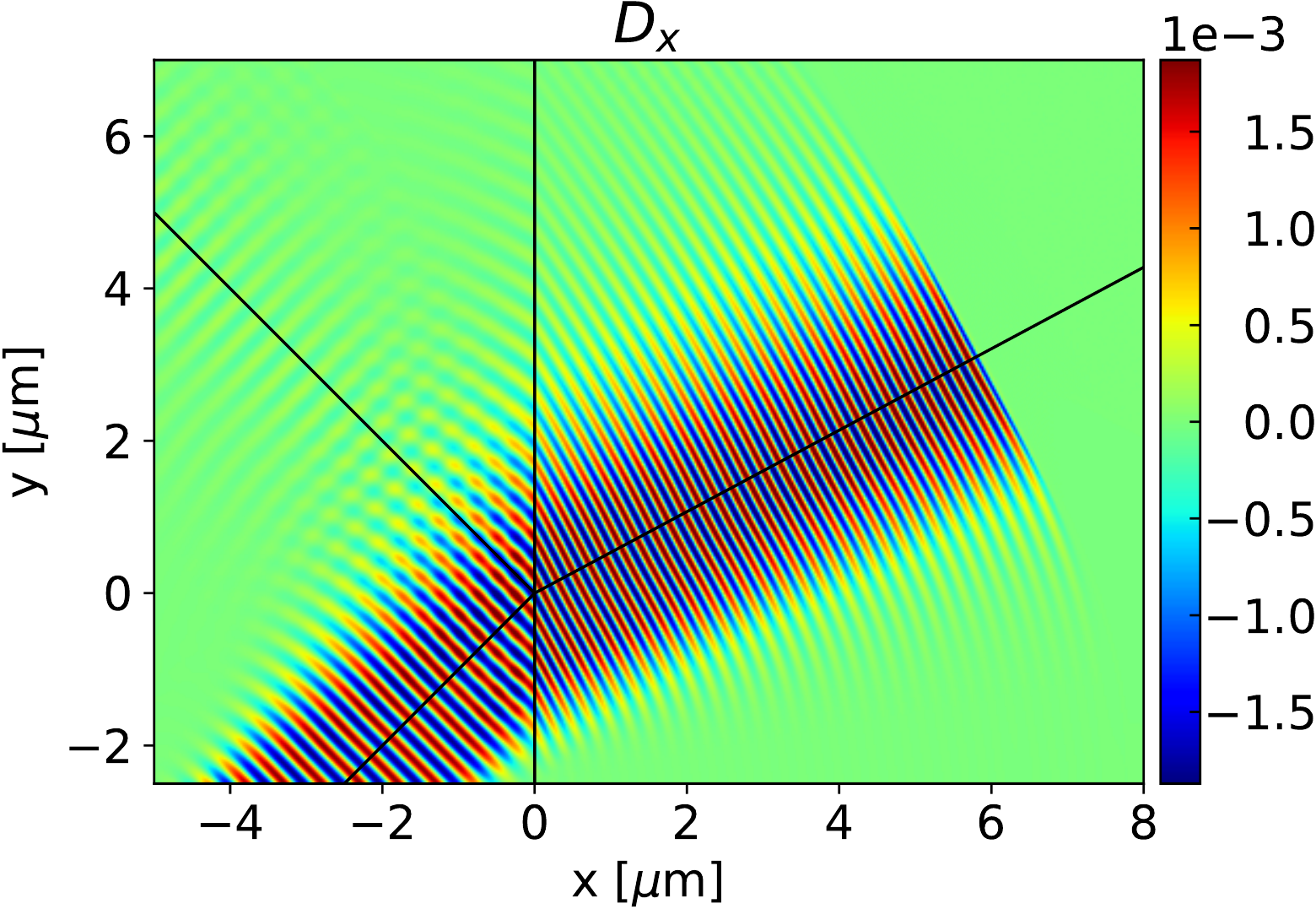} &
\includegraphics[width=0.33\textwidth]{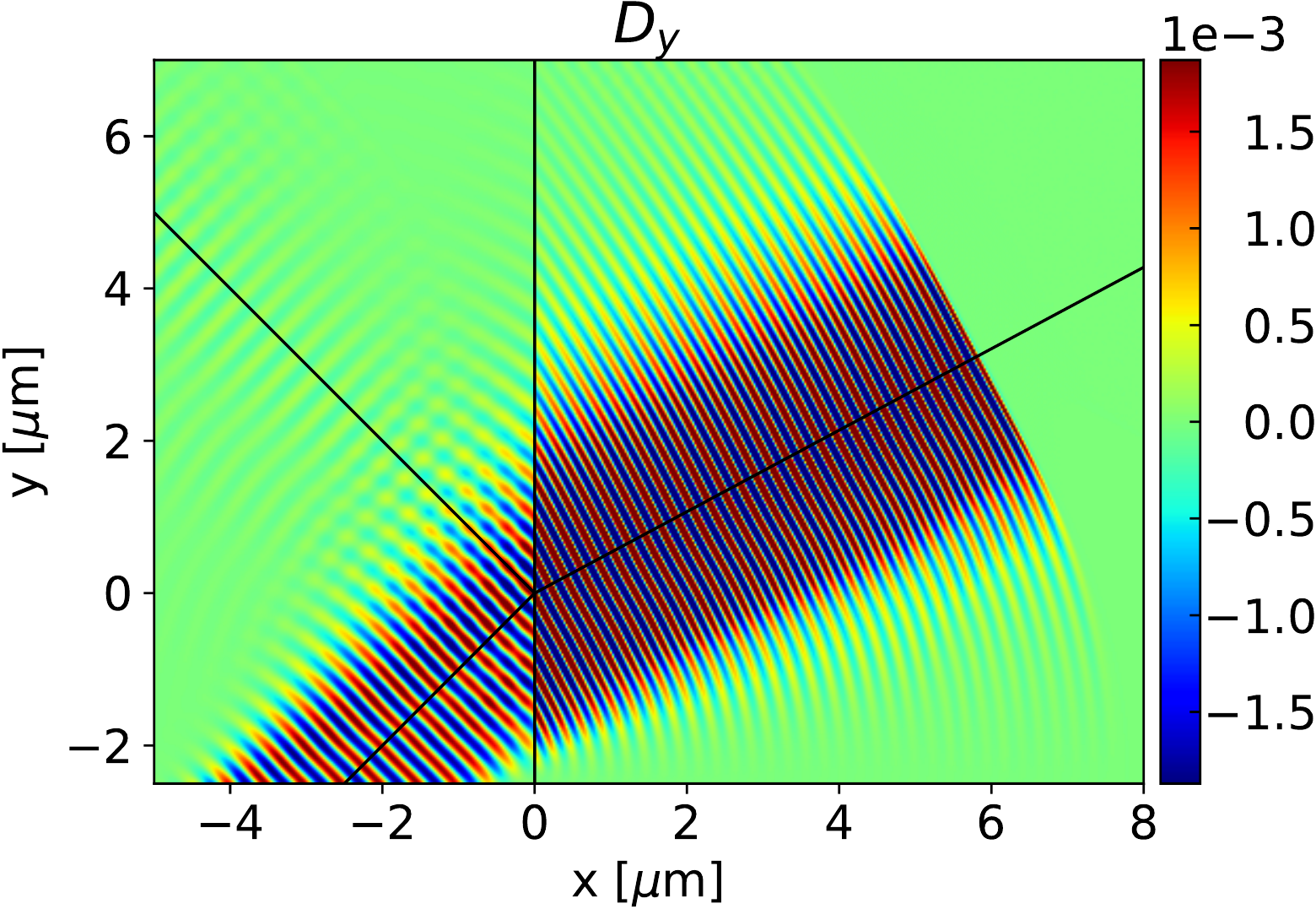}
\end{tabular}
\caption{Refraction of a compact electromagnetic beam by a dielectric slab on a mesh of $650 \times 475$ cells. Top row: initial condition, middle row: $k=3$, bottom row: $k=4$}
\label{fig:refrac}
\end{center}
\end{figure}

This test case describes the refraction of an electromagnetic beam by a two dimensional dielectric slab that spans over a domain $[-5.0, 8.0]\times[-2.5, 7.0]~\si{\square \micro \meter}$. The domain is divided into $\num{650x475}$ cells and has a constant permeability $\mu_0$ and the permittivity is given by $\epsilon(x,y) = 1.625 \epsilon_0 + 0.625 \epsilon_0 \tanh(\num{e08} x)$ so as to model a dielectric slab where  the permittivity changes from $\epsilon= 2.25\epsilon_0$ for $x~\geq~0$ to $\epsilon_0$ for $x<0$. The magnetic and electric vector potential of the incident electromagnetic beam are given by,
\begin{align}
\bm{A}(x,y,t)=&\frac{\lambda}{8\pi}\sin[2\pi(x+y-\sqrt{2}ct)]\bigg[1-\tanh\bigg(\frac{(x-a)+(y-b)-\sqrt{2}ct}{0.1\lambda}\bigg)\bigg]\notag\\
                &\bigg[1-\tanh\bigg(\frac{\abs{y-x}-\sqrt{2}{d}}{\sqrt{2}\delta}\bigg)\bigg]\hat{e}_y \label{eq:refracmag}\\
\bm{C}(x,y,t)=&-\frac{\lambda}{8\pi\sqrt{2}}\sin[2\pi(x+y-\sqrt{2}ct)]\bigg[1-\tanh\bigg(\frac{(x-a)+(y-b)-\sqrt{2}ct}{0.1\lambda}\bigg)\bigg]\notag\\
                &\bigg[1-\tanh\bigg(\frac{\abs{y-x}-\sqrt{2}{d}}{\sqrt{2}\delta}\bigg)\bigg]\hat{e}_z \label{eq:refracelec}
\end{align}
where wavelength of the incident beam $\lambda = \SI{0.5}{\micro\m}$ and other geometrical parameters are taken to be $d= 2.5\lambda, \delta = 0.5 \lambda$ and $(a,b)=(-3.0\lambda,-3.0\lambda)$. The impinging beam of radiation is incident on the surface of the dielectric slab at an angle of \ang{45}. The simulation was run to a time of \SI{4.0E-14}{\s}. Figure~(\ref{fig:refrac}) shows the initial condition and the solution at the final time obtained from fourth and fifth order schemes.
This problem, involving the refraction of a beam of radiation, was presented in Balsara et al.~\cite{Balsara2017},~\cite{Balsara2018}. Figure~(\ref{fig:refrac}) shows our results for DGTD schemes at fourth and fifth orders. Despite the use of high order DG schemes with rapidly varying dielectric properties, we never needed to use limiters in this problem. We see that our results are very consistent with the reference solution presented in those papers.

\subsection{Total internal reflection of a compact electromagnetic beam by a dielectric slab}
This test case is designed to simulate total internal reflection of an electromagnetic beam by a dielectric slab which has a constant permeability of $\mu_0$ and the permittivity is given by $\epsilon(x,y,z) = 2.5\epsilon_0 -1.5\epsilon_0\tanh(\num{4.0e08}x)$. Across the dielectric slab, the permittivity changes from $\epsilon= 4.0\epsilon_0$ for $x~\leq~0$ to $\epsilon_0$ for $x>0$ which implies that the refractive index is \num{2} for the dielectric slab. Therefore, following Snell's law, the critical angle for internal reflection in this dielectric slab is \ang{30}.

\begin{figure}
\begin{center}
\begin{tabular}{ccc}
\includegraphics[width=0.33\textwidth]{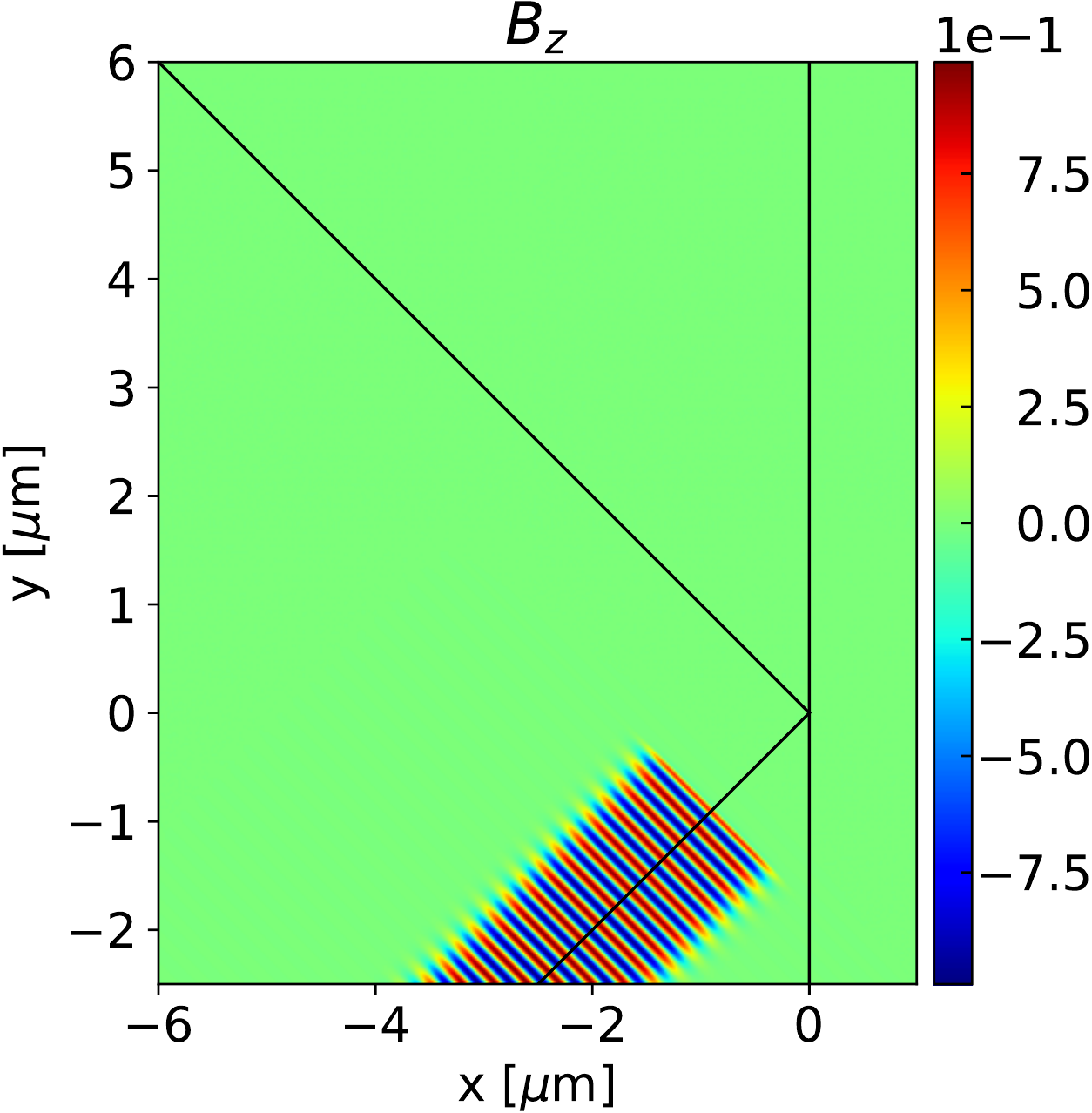} &
\includegraphics[width=0.33\textwidth]{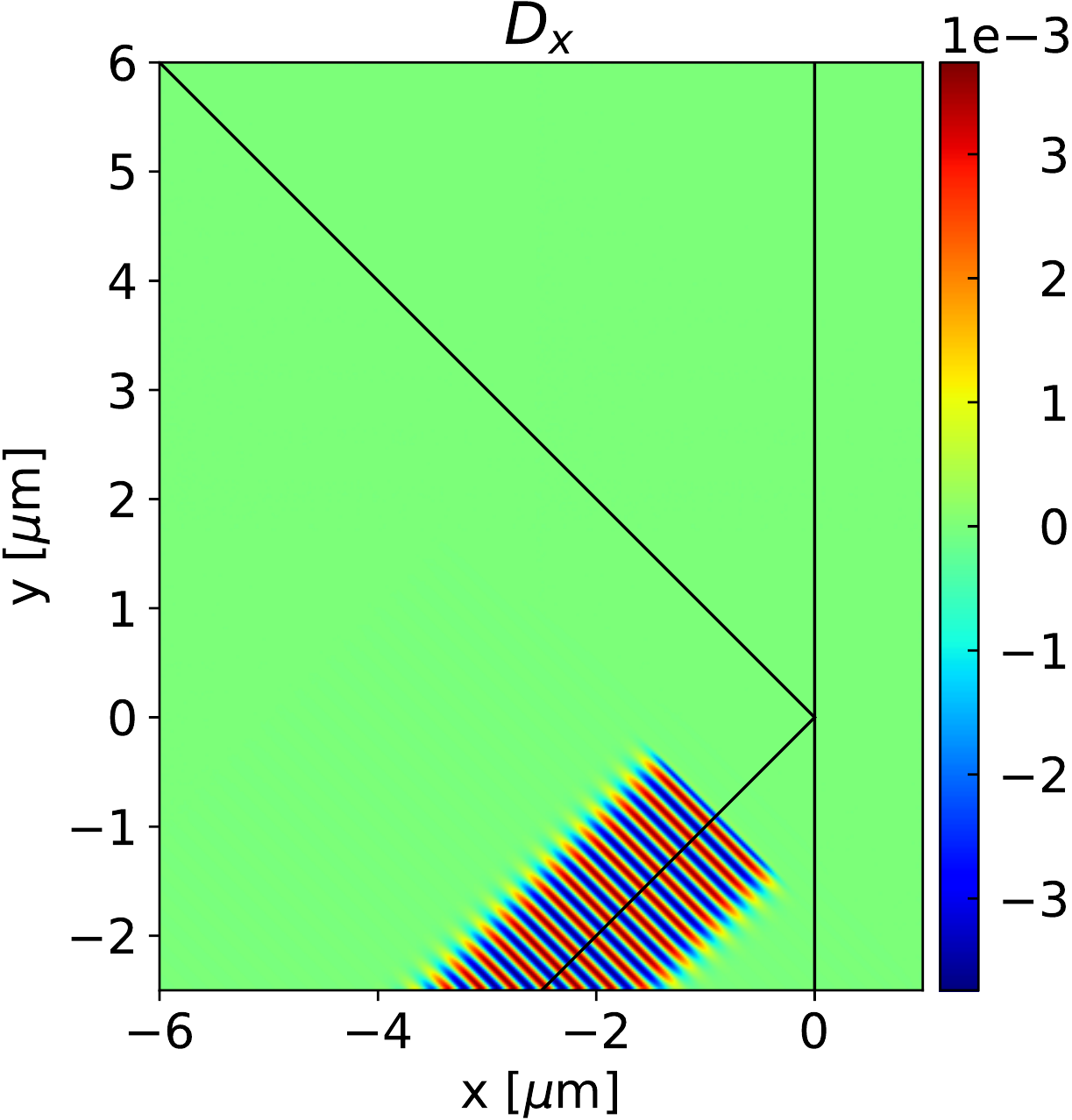} &
\includegraphics[width=0.33\textwidth]{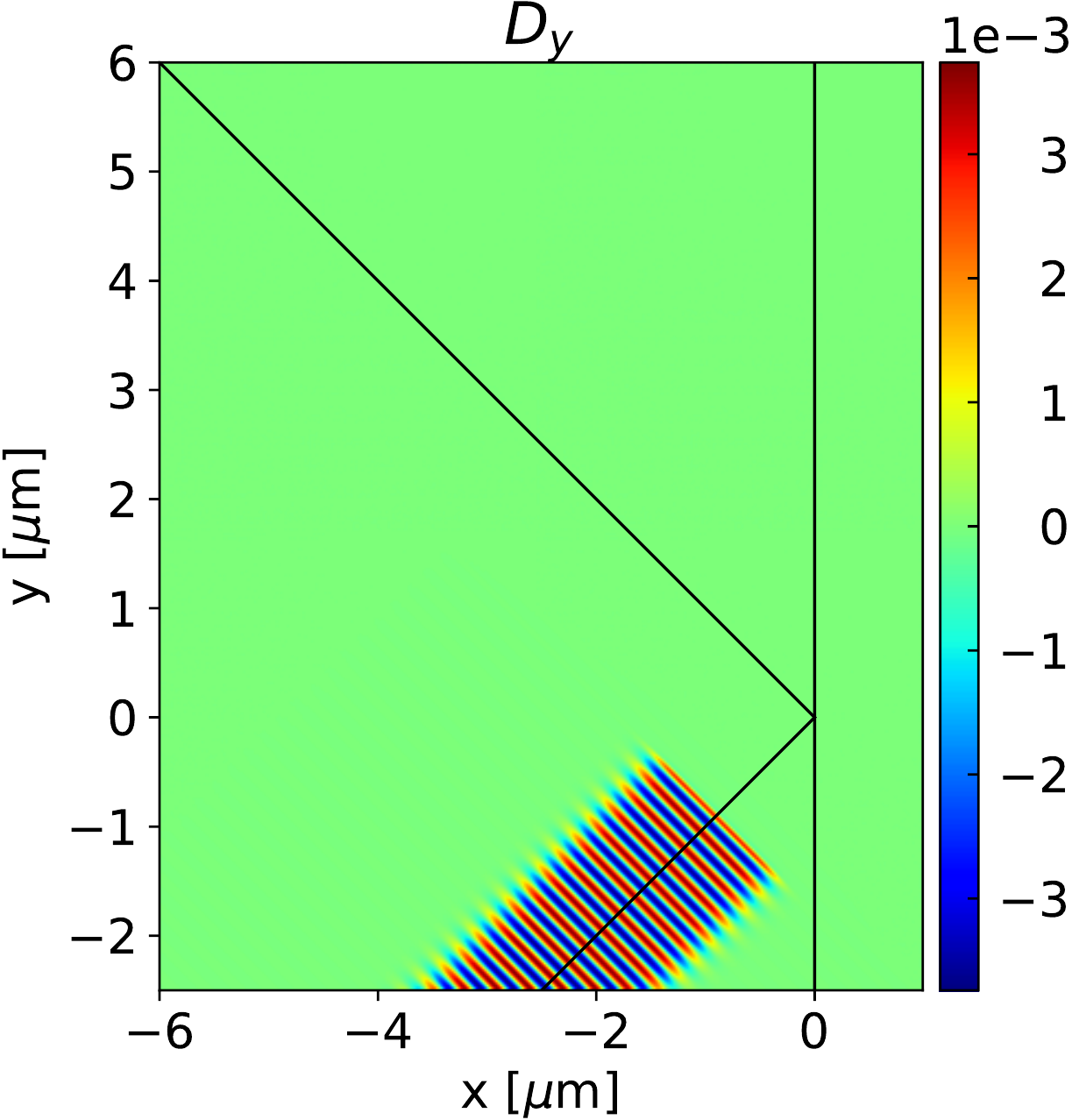} \\
\includegraphics[width=0.33\textwidth]{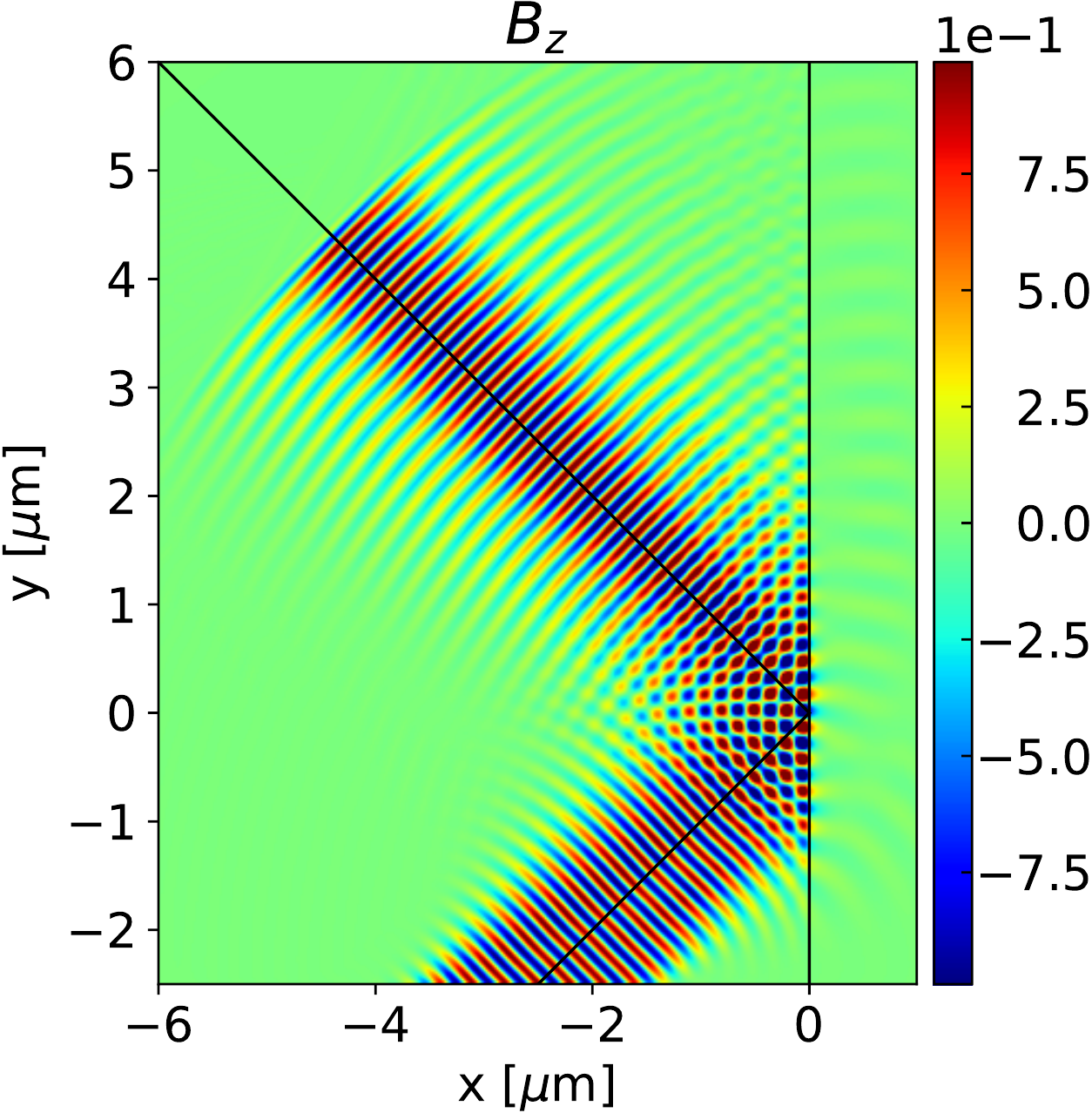} &
\includegraphics[width=0.33\textwidth]{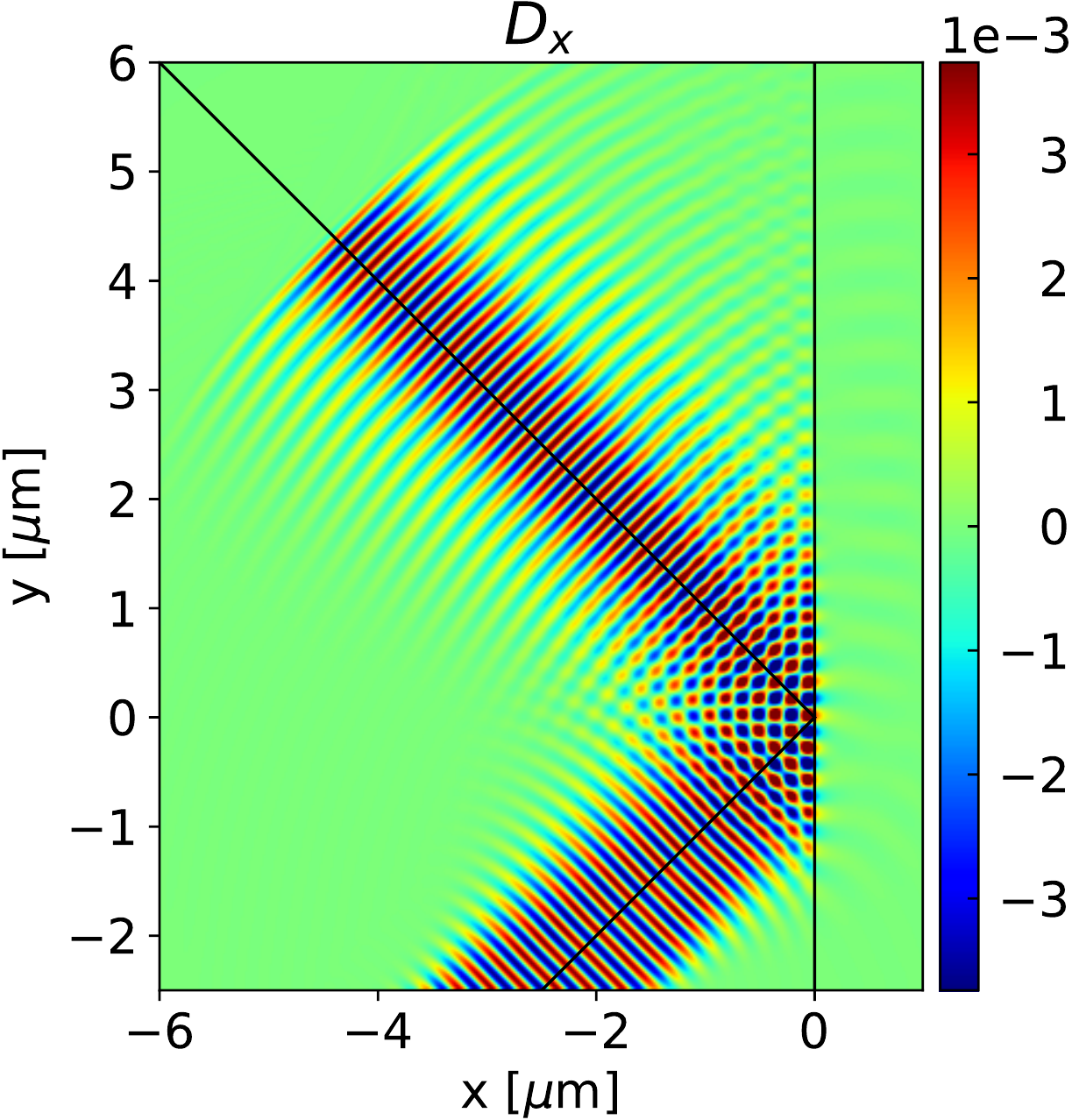} &
\includegraphics[width=0.33\textwidth]{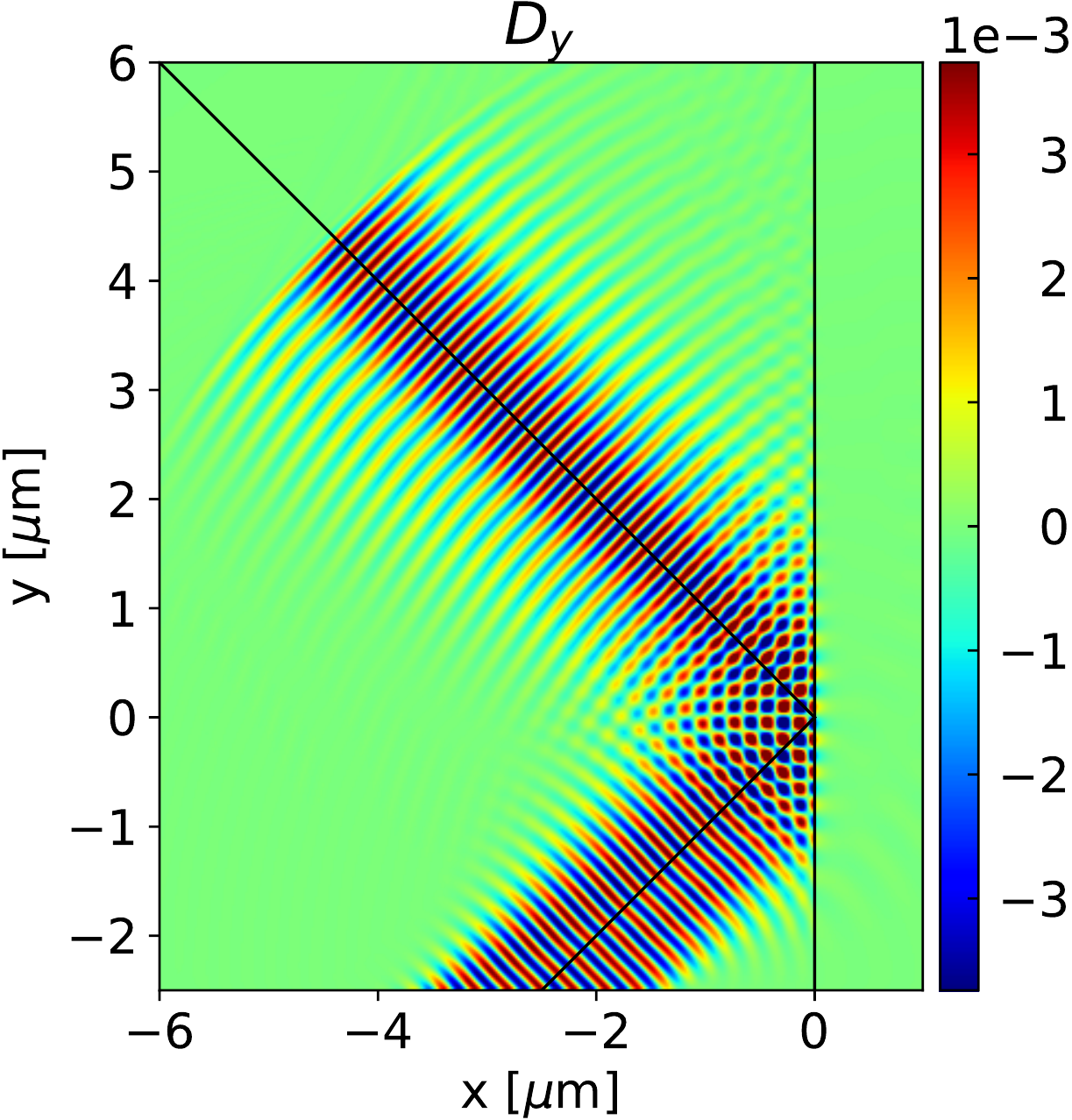} \\
\includegraphics[width=0.33\textwidth]{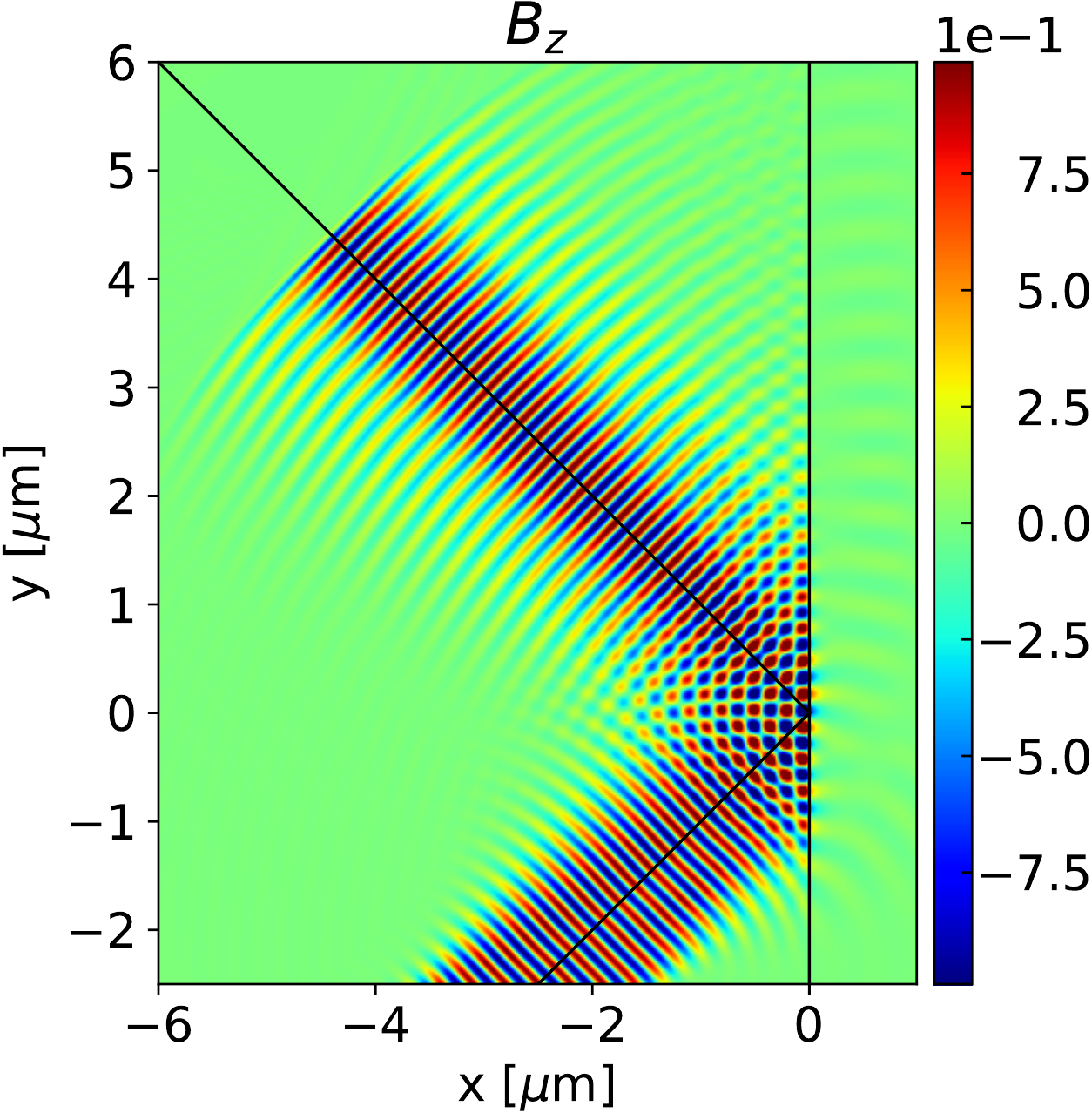} &
\includegraphics[width=0.33\textwidth]{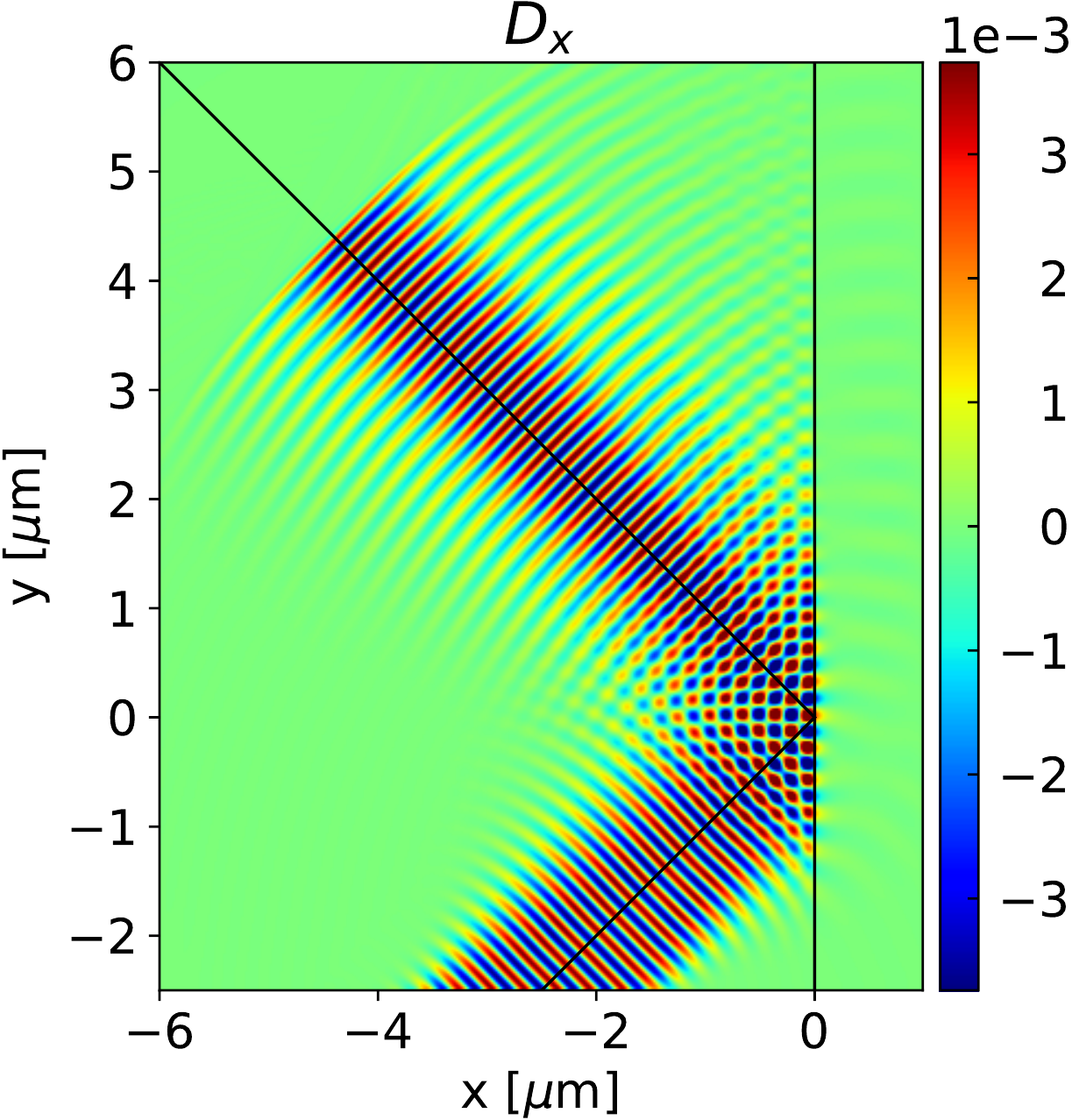} &
\includegraphics[width=0.33\textwidth]{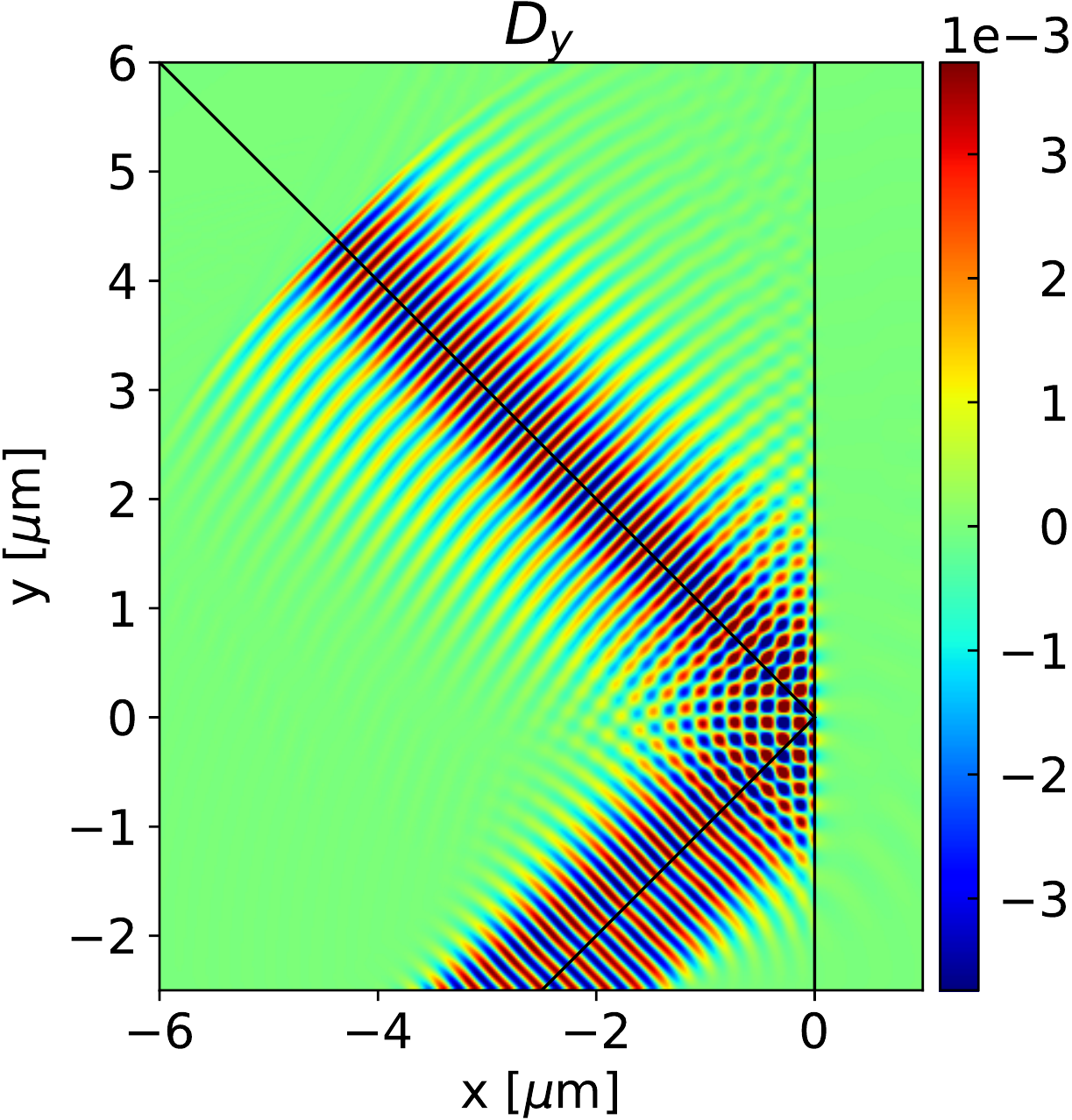}
\end{tabular}
\caption{Total internal reflection of a compact electromagnetic beam by a dielectric slab on a mesh of $350 \times 425$ cells. Top row: initial condition, middle row: $k=3$, bottom row: $k=4$}
\label{fig:reflec}
\end{center}
\end{figure}

The simulation is performed in a domain $[-6.0, 1.0]\times[-2.5, 6.0]~\si{\square \micro \meter}$ of such a dielectric slab 
divided in $\num{350x425}$ cells over a duration of \SI{5.0E-14}{\s}. The initial  conditions are similar to the previous case. However, for this problem, we chose $\lambda = \SI{0.3}{\micro\m}$, $d= 2.5\lambda$, $\delta = 0.5 \lambda$ and $(a,b)=(-3.0\lambda,-3.0\lambda)$. Figure~(\ref{fig:reflec}) shows the initial condition and the final solution obtained using the fourth and fifth order schemes. 
This problem, involving the total internal reflection of a beam of radiation, was presented in Balsara et al.~\cite{Balsara2017},~\cite{Balsara2018}. Figure~(\ref{fig:reflec}) shows our results for DGTD schemes at fourth and fifth orders. As in the previous test problem, we found that despite the use of high order DG schemes with rapidly varying dielectric properties, we never needed to use limiters in this problem. We see that our results are very consistent with the reference solution presented in those papers.

\section{Summary and conclusions}
In Balsara and K\"appeli~\cite{Balsara2018a} globally constraint preserving DG schemes (up to fourth order) for CED had been presented and their von Neumann stability analysis had been carried out. The von Neumann analysis gave many insights as regards to the CFL condition and indicated that superlative propagation of electromagnetic radiation could be achieved. However, that paper did not explore the role of varying permittivity and permeability. That aspect of numerical CED has been explored here. Besides, since we have designed fifth order DG schemes in this paper, we are able to get a clearer view of all the ingredients of a globally constraint preserving DG scheme for CED at all orders.

	Our DG schemes are novel because they can be viewed as retaining many of the best aspects of the FDTD schemes; while finding a pathway to higher order extensions. Our first conclusion is that at fourth and higher orders of accuracy, one has to evolve some zone-centered modes in addition to the face-centered modes.
	
	It is also well-known that the best way to operate a DG scheme is to use it without reliance on limiters; if the physical problem and the system of equations permit this. We have carried out several tests where the permittivity varied by almost an order of magnitude. For all these tests, we were able to run the DG scheme without using limiters. In fact, we never had to use limiters for any of the problems that are presented in this paper. This leads us to our second important conclusion that DG schemes of the sort designed here do not seem to require limiting in order to stabilize them. Just the logic of the Riemann solvers, along with the natural global constraint preservation that is in-built into the scheme, proved sufficient to keep our DG schemes stable at all orders.
	
	All our DG schemes evolve the facial electric displacement and magnetic induction and their higher moments as the primal variables. Maxwell's equations, without the presence of conductivity, ensure the conservation of electromagnetic energy. In our schemes we do not do anything special to conserve electromagnetic energy. However, the DG philosophy helps because it provides evolution for all the modes, ensuring that our DG schemes retain the spirit of the governing PDE, i.e. Maxwell's equations. Out third conclusion is that the DG schemes developed here have excellent ability to conserve electromagnetic energy on the computational mesh even when permittivity and permeability vary strongly in space; as long as the conductivity is zero. This is especially true for the fourth and higher order DG schemes presented here.
	
	The three conclusions establish DG schemes as strong performers for globally constraint-preserving CED with several very favorable properties.
\section*{Acknowledgements}
DSB acknowledges support via NSF grants NSF-ACI-1533850, NSF-DMS-1622457 , NSF\_ACI-1713765 and NSF-DMS-1821242. Support from a grant by Notre Dame International and the Airbus Chair on Mathematics of Complex Systems at TIFR-CAM, Bangalore, is also gratefully acknowledged.

\appendix
\section{Divergence-free reconstruction at second and third order}
The reconstruction at fourth order accuracy ($k=3$) has been detailed in section~(\ref{sec:rec3}). Using this we can obtain the divergence-free reconstruction at second and third orders as follows. At degree $k=1$, the field $\D$ is approximated on the faces and inside the cell as
\[
D_x^\pm(\eta) = a_0^\pm + a_1^\pm \phi_1(\eta), \qquad D_y^\pm(\xi) = b_0^\pm + b_1^\pm \phi_1(\xi)
\]
\[
D_x(\xi,\eta) = a_{00} + a_{10}\phi_1(\xi) + a_{01}\phi_1(\eta) + a_{20}\phi_2(\xi) + a_{11} \phi_1(\xi) \phi_1(\eta)
\]
\[
D_y(\xi,\eta) = b_{00} + b_{10}\phi_1(\xi) + b_{01}\phi_1(\eta) + b_{11}\phi_1(\xi) \phi_1(\eta) + b_{02}\phi_2(\eta)
\]
At degree $k=2$, the field $\D$ on the faces and inside the cell is approximated by
\[
D_x^\pm(\eta) = a_0^\pm + a_1^\pm \phi_1(\eta) + a_2^\pm \phi_2(\eta), \qquad D_y^\pm(\xi) = b_0^\pm + b_1^\pm \phi_1(\xi) + b_2^\pm \phi_2(\xi)
\]
\[
\begin{aligned}
D_x(\xi,\eta) = \ & a_{00} + a_{10}\phi_1(\xi) + a_{01}\phi_1(\eta) + a_{20}\phi_2(\xi) + a_{11} \phi_1(\xi) \phi_1(\eta) + a_{02}\phi_2(\eta) + \\
& a_{30} \phi_3(\xi)  + a_{12}\phi_1(\xi)\phi_2(\eta) \\
D_y(\xi,\eta) = \ & b_{00} + b_{10} \phi_1(\xi) + b_{01}\phi_1(\eta) + b_{20}\phi_2(\xi) + b_{11}\phi_1(\xi)\phi_1(\eta) + b_{02}\phi_2(\eta) + \\
& b_{21}\phi_2(\xi) \phi_1(\eta) + b_{03} \phi_3(\eta)
\end{aligned}
\]
The coefficients $a_{ij}$, $b_{ij}$ in the cell solution can be obtained from the formulae in  section~(\ref{sec:rec3}) by setting those coefficients which do not appear above to zero. Note that at these orders, the face solution completely determines the divergence-free reconstruction inside the cells, unlike at fourth and fifth orders, where additional information in terms of $\omega_i$ is required to complete the divergence-free reconstruction.
\section{Divergence-free reconstruction at fifth order}
\label{sec:rec4}
Let us start by assuming a form for the electric displacement using BDFM polynomial. Hence we take $D_x \in \poly_5 \setminus \{ \eta^5\}$ and $D_y \in \poly_5 \setminus \{ \xi^5\}$ which is the form of a BDFM polynomial~\cite{Brezzi1987}. The components of the vector field can be written in terms of the orthogonal basis functions as follows
\begin{align*}
D_x(\xi,\eta) = &\ a_{00} + a_{10}\phi_1(\xi) + a_{01}\phi_1(\eta) + a_{20}\phi_2(\xi) + a_{11} \phi_1(\xi)\phi_1(\eta) + a_{02}\phi_2(\eta) + \\
& a_{30} \phi_3(\xi) + a_{21}\phi_2(\xi)\phi_1(\eta) + a_{12}\phi_1(\xi)\phi_2(\eta) + a_{03}\phi_3(\eta) + a_{40}\phi_4(\xi) +\\
& a_{31}\phi_3(\xi)\phi_1(\eta) + a_{22}\phi_2(\xi)\phi_2(\eta) + a_{13}\phi_1(\xi)\phi_3(\eta)+ a_{04}\phi_4(\eta)  + a_{50}\phi_5(\xi)+\\
& a_{41}\phi_4(\xi)\phi_1(\eta) + a_{32}\phi_3(\xi)\phi_2(\eta)+ a_{23}\phi_2(\xi)\phi_3(\eta) + a_{14}\phi_1(\xi)\phi_4(\eta) \\
D_y(\xi,\eta) = &\ b_{00} + b_{10}\phi_1(\xi) + b_{01}\phi_1(\eta) + b_{20}\phi_2(\xi) + b_{11}\phi_1(\xi)\phi_1(\eta) + b_{02}\phi_2(\eta) + \\
& b_{30}\phi_3(\xi) + b_{21}\phi_2(\xi)\phi_1(\eta) + b_{12}\phi_1(\xi) \phi_2(\eta) + b_{03} \phi_3(\eta) +  b_{31}\phi_3(\xi)\phi_1(\eta) + \\
& b_{22}\phi_2(\xi)\phi_2(\eta) + b_{13}\phi_1(\xi)\phi_3(\eta) + b_{04}\phi_4(\eta)+ b_{40} \phi_4(\xi) + b_{05} \phi_5(\eta)+ \\
& b_{41} \phi_4(\xi)\phi_1(\eta) + b_{32}\phi_3(\xi)\phi_2(\eta)+ b_{23}\phi_2(\xi)\phi_3(\eta) + b_{14}\phi_1(\xi)\phi_4(\eta)
\end{align*}
which has a total of 40 coefficients.  By matching the above polynomial to the solution on the faces, we obtain the following set of 20 equations.
\begin{align*}
a_{00} \pm \shalf a_{10} + \tfrac{1}{6}a_{20} \pm \tfrac{1}{20}a_{30} + \tfrac{1}{70}a_{40} \pm \tfrac{1}{252}a_{50} & = a_0^\pm \\
a_{01} \pm \shalf a_{11} + \tfrac{1}{6}a_{21} \pm \tfrac{1}{20}a_{31} + \tfrac{1}{70}a_{41}  & = a_1^\pm \\
a_{02} \pm \shalf a_{12} + \tfrac{1}{6}a_{22} \pm \tfrac{1}{20}a_{32} & = a_2^\pm \\
a_{03} \pm \shalf a_{13} + \tfrac{1}{6}a_{23} & = a_3^\pm \\
a_{04} \pm \shalf a_{14} & = a_4^\pm \\
b_{00} \pm \shalf b_{01} + \tfrac{1}{6}b_{02} \pm \tfrac{1}{20}b_{03} + \tfrac{1}{70}b_{04} \pm \tfrac{1}{252}b_{05}  & = b_0^\pm \\
b_{10} \pm \shalf b_{11} + \tfrac{1}{6}b_{12} \pm \tfrac{1}{20}b_{13} + \tfrac{1}{70}b_{14}  & = b_1^\pm \\
b_{20} \pm \shalf b_{21} + \tfrac{1}{6}b_{22} \pm \tfrac{1}{20}b_{23}  & = b_2^\pm \\
b_{30} \pm \shalf b_{31} + \tfrac{1}{6}b_{32} & = b_3^\pm \\
b_{40} \pm \shalf b_{41} & = b_4^\pm
\end{align*}                                                                    
Setting the divergence to zero yields the following set of 15 equations
\begin{align*}
(a_{10}+\tfrac{1}{10}a_{30} + \tfrac{1}{126} a_{50})\dy+(b_{01} +\tfrac{1}{10}b_{03} + \tfrac{1}{126} b_{05})\dx &= 0\\
(2 a_{20}+\tfrac{6}{35}a_{40})\dy + (b_{11} +b_{13}/10)\dx) &= 0\\
(a_{11}+a_{31}/10) \dy + (2 b_{02} +\tfrac{6}{35} b_{04}) \dx &= 0\\
(3 a_{30} + \tfrac{5}{21} a_{50}) \dy + (b_{21} + \tfrac{1}{10} b_{23}) \dx &= 0\\
(2 a_{21} + \tfrac{6}{35} a_{41}) \dy +(2 b_{12} + \tfrac{6}{35} b_{14}) \dx &= 0\\
(a_{12} + \tfrac{1}{10} a_{32}) \dy + (3 b_{03} + \tfrac{5}{21} b_{05}) \dx &= 0 \\
 4 a_{40}\dy +   b_{31} \dx &= 0 \\
 3 a_{31}\dy + 2 b_{22} \dx &= 0 \\
 2 a_{22}\dy + 3 b_{13} \dx &= 0 \\
   a_{13}\dy + 4 b_{04} \dx &= 0 \\
  5a_{50}\dy +   b_{41} \dx &=0  \\
  4a_{41}\dy + 2b_{32}  \dx &=0  \\
  3a_{32}\dy + 3b_{23}  \dx &=0  \\
  2a_{23}\dy + 4b_{14}  \dx &=0  \\
   a_{14}\dy + 5b_{05}  \dx &=0 
\end{align*}
The first equation in the above set of equations is redundant as it is included in the other equations due to the divergence-free constraint~(\ref{eq:compat}). Ignoring this equation, we can solve for some of the coefficients in terms of the face solution as follows.

\boxed{
\begin{minipage}[h]{0.5\textwidth}
\begin{align*}
a_{00}=& \shalf (a_0^- + a_0^+)  + \tfrac{1}{12} (b_1^+ - b_1^-) \tfrac{\dx}{\dy} \\
a_{10}=& a_0^+ - a_0^- + \tfrac{1}{30}(b_2^+ - b_2^-) \tfrac{\dx}{\dy} \\ 
a_{20}=& -\shalf (b_1^+ - b_1^-) \tfrac{\dx}{\dy} + \tfrac{3}{140}(b_3^+ - b_3^-) \tfrac{\dx}{\dy}\\
a_{30}=& -\tfrac{1}{3}(b_2^+ - b_2^-) \tfrac{\dx}{\dy} + \tfrac{1}{63}(b_4^+ - b_4^-) \tfrac{\dx}{\dy} \\
a_{40}=& -\tfrac{1}{4}(b_3^+ - b_3^-) \tfrac{\dx}{\dy} \\
a_{50}=& -\tfrac{1}{5}(b_4^+ - b_4^-) \tfrac{\dx}{\dy} \\
a_{04}=& \shalf (a_4^- + a_4^+) \\
a_{13}=& a_3^+ - a_3^- \\
a_{14}=& a_4^+ - a_4^- \\
\end{align*}
\end{minipage}
\noindent
\begin{minipage}[h]{0.5\textwidth}
\begin{align*}
b_{00}=& \shalf (b_0^- + b_0^+) + \tfrac{1}{12} (a_1^+ - a_1^-) \tfrac{\dy}{\dx} \\
b_{01}=& b_0^+ - b_0^- + \tfrac{1}{30}(a_2^+ - a_2^-) \tfrac{\dy}{\dx} \\
b_{02}=& -\shalf (a_1^+ - a_1^-) \tfrac{\dy}{\dx} + \tfrac{3}{140}(a_3^+ - a_3^-) \tfrac{\dy}{\dx} \\
b_{03}=& -\tfrac{1}{3} (a_2^+ - a_2^-) \tfrac{\dy}{\dx}+\tfrac{1}{63}(a_4^+ - a_4^-) \tfrac{\dy}{\dx} \\
b_{04}=& -\tfrac{1}{4} (a_3^+ - a_3^-) \tfrac{\dy}{\dx} \\
b_{05}=& -\tfrac{1}{5} (a_4^+ - a_4^-) \tfrac{\dy}{\dx} \\
b_{40}=& \shalf (b_4^- + b_4^+) \\
b_{31}=& b_3^+ - b_3^- \\
b_{41}=& b_4^+ - b_4^- \\
\end{align*}

\end{minipage}
}

\noindent
The remaining coefficients satisfy the following equations

\begin{minipage}[h]{0.5\textwidth}
\begin{align*}
a_{01} + \tfrac{1}{6} a_{21} +\tfrac{1}{70} a_{41} &=\shalf (a_1^+ + a_1^-) \\
a_{03} + \tfrac{1}{6} a_{23} &= \shalf (a_3^+ + a_3^-) \\
a_{02} + \tfrac{1}{6} a_{22} &= \shalf (a_2^+ + a_2^-) \\
a_{12} + \tfrac{1}{10} a_{32} &= a_2^+ - a_2^- \\
a_{11} + \tfrac{1}{10}a_{31} &= a_1^+ - a_1^- \\
b_{10} + \tfrac{1}{6} b_{12} +\tfrac{1}{70} b_{14} &= \shalf (b_1^+ + b_1^-) \\
b_{20} + \tfrac{1}{6} b_{22} &= \shalf (b_2^+ + b_2^-)\\
b_{30} + \tfrac{1}{6} b_{32} &= \shalf (b_3^+ + b_3^-) \\
b_{21} + \tfrac{1}{10} b_{23}&= b_2^+ - b_2^- \\
b_{11} + \tfrac{1}{10}b_{13} &= b_1^+ - b_1^- 
\end{align*}
\end{minipage}
\begin{minipage}[h]{0.5\textwidth}
\begin{align*}
(2 a_{21} + \tfrac{6}{35} a_{41}) \dy +(2 b_{12} + \tfrac{6}{35} b_{14}) \dx &= 0\\
 3 a_{31}\dy + 2 b_{22} \dx &=0\\
 2 a_{22}\dy + 3 b_{13} \dx &=0\\
  4a_{41}\dy + 2b_{32}  \dx &=0  \\
  3a_{32}\dy + 3b_{23}  \dx &=0  \\
  2a_{23}\dy + 4b_{14}  \dx &=0  \\
\end{align*}
\end{minipage}

\noindent
We have more unknowns than equations, so we have to make some further assumptions on the remaining coefficient. Let set all the coefficients at degree five  to zero since they are not required to get fifth order accuracy, i.e.,
\[
a_{41} = a_{32} = b_{23} = b_{14} = a_{23} = b_{32} = 0
\]
and we can immediately obtain the solution for the following coefficients
\[
\boxed{
a_{03} = \shalf (a_3^+ + a_3^-), \quad
a_{12} = a_2^+ - a_2^-, \quad
b_{30} = \shalf (b_3^+ + b_3^-), \quad
b_{21} = b_2^+ - b_2^-
}
\]
The remaining equations and unknowns can be broken into two sets of equations. The first set is of the form
\begin{eqnarray}
\label{eq:bdfm31}
a_{01} + \frac{1}{6}a_{21} &=& \half (a_1^- + a_1^+) =: r_1 \\
\label{eq:bdfm32}
b_{10} + \frac{1}{6} b_{12} &=& \half (b_1^- + b_1^+) =: r_2 \\
\label{eq:bdfm33}
b_{12} \dx + a_{21} \dy &=& 0
\end{eqnarray}
Here we have four unknowns but only three equations. We can solve these equations by introducing the additional variable $b_{10} - a_{01}  = \omega_1$ and the solution is same as in the fourth order case given in section~(\ref{sec:rec3}).

We are now left with the following set of four equations
\begin{align}
a_{11} + \frac{1}{10} a_{31} &= (a_1^+ - a_1^-)=:r_3        \quad & a_{02} + \frac{1}{6}  a_{22} &= \shalf (a_2^+ + a_2^-)=:r_5  \\
b_{20} +  \frac{1}{6} b_{22} &= \shalf (b_2^+ + b_2^-)=:r_4 \quad & b_{11} + \frac{1}{10} b_{13} &= (b_1^+ - b_1^-)=:r_6 
\end{align}
but there are six unknown coefficients. All of these coefficients are at or below degree four and have to be retained for fifth order accuracy. We need additional information to solve for all the coefficients and we introduce the following two equations
\begin{align*}
b_{20} - a_{11}  = \omega_2, \qquad b_{11} - a_{02}  = \omega_3
\label{eq:omg2}
\end{align*}
Then we can solve for all the remaining unknown coefficients to obtain the following solution
\begin{align*}
a_{11} &=\frac{1}{2 + 5\frac{\dy}{\dx}}\bigg(5 r_3 \frac{\dy}{\dx} + 2 r_4 - 2 \omega_2\bigg) && a_{02} =\frac{1}{5 + 2\frac{\dy}{\dx}}\bigg(2 r_5 \frac{\dy}{\dx} + 5 r_6 - 5 \omega_3\bigg)\\
a_{31} &= 10(r_3 - a_{11})                                                                    && a_{22} =6(r_5 - a_{02})\\
b_{20} &= \omega_2 + a_{11}                                                                   && b_{11} =\omega_3 + a_{02} \\
b_{22} &= 6(r_4 - b_{20})                                                                     && b_{13} =10(r_6 - b_{11})
\end{align*}
This completely specifies the reconstructed field inside the cell. The evolution equations for $\omega_2, \omega_3$ are explained in section~(\ref{sec:five}).
\section{Numerical fluxes}
The fluxes in the conservative form of the 2-D Maxwell's equations~(\ref{eq:claw}) have the form $\fx = A_x \con$, $\fy = A_y \con$ where the matrices $A_x$, $A_y$ may depend on the spatial coordinate due to varying material properties and are given by
\[
A_x = \begin{bmatrix}
0 & 0 & 0 \\
0 & 0 & 1/\mu \\
0 & 1/\varepsilon & 0 \end{bmatrix}, \qquad A_y = \begin{bmatrix}
0 & 0 & -1/\mu \\
0 & 0 & 0 \\
-1/\varepsilon & 0 & 0 \end{bmatrix}
\]
For any unit vector $\un = (n_x, n_y)$, the matrix 
\[
A_n = A_x n_x + A_y n_y = \begin{bmatrix}
0 & 0 & -n_y/\mu \\
0 & 0 & n_x/\mu \\
-n_y/\varepsilon & n_x/\varepsilon & 0 \end{bmatrix}
\]
has real eigenvalues given by $\{-c, 0, +c\}$ where $c = 1/\sqrt{\varepsilon\mu}$, and a complete set of eigenvectors given by
\[
\begin{bmatrix}
+\sqrt{\varepsilon} n_y \\
-\sqrt{\varepsilon} n_x \\
\sqrt{\mu} \end{bmatrix}, \qquad 
\begin{bmatrix}
n_x \\
n_y \\
0 \end{bmatrix}, \qquad
\begin{bmatrix}
-\sqrt{\varepsilon} n_y \\
+\sqrt{\varepsilon} n_x \\
\sqrt{\mu} \end{bmatrix}
\]
\subsection{Solution of the 1-D Riemann problem}
\label{sec:riem1d}
Consider the two states $\con^L$, $\con^R$ separated across an interface with normal vector $\un$ which points from $L$ to $R$, and such that $\D^L \cdot \un = \D^R \cdot \un$, i.e., the normal component of $\D$ is continuous. The flux in the direction $\un$ is 
\[
\fx_n = \fx n_x + \fy n_y = \begin{bmatrix}
-\frac{n_y}{\mu} B_z \\
+\frac{n_x}{\mu} B_z \\
\frac{1}{\varepsilon}(D_y n_x - D_x n_y) \end{bmatrix}
\]
Let the intermediate states be denoted by $\con^*$, $\con^{**}$. Then the jump conditions across the three waves are \\
1) Across the $-c$ wave
\begin{eqnarray}
\label{eq:mc1}
-\frac{n_y}{\mu}(B_z^* - B_z^L) &=& -c (D_x^* - D_x^L) \\
\label{eq:mc2}
+\frac{n_x}{\mu}(B_z^* - B_z^L) &=& -c (D_y^* - D_y^L) \\
\label{eq:mc3}
\frac{1}{\varepsilon}[ (D_y^* n_x - D_x^* n_y) - (D_y^L n_x - D_x^Ln_y)] &=& -c (B_z^* - B_z^L)
\end{eqnarray}
2) Across the $0$ wave
\begin{eqnarray}
\label{eq:z1}
-\frac{n_y}{\mu}(B_z^{**} - B_z^*) &=& 0 \\
\label{eq:z2}
+\frac{n_x}{\mu}(B_z^{**} - B_z^*) &=& 0 \\
\label{eq:z3}
\frac{1}{\varepsilon}[ (D_y^{**} n_x - D_x^{**} n_y) - (D_y^* n_x - D_x^* n_y)] &=& 0
\end{eqnarray}
3) Across the $+c$ wave
\begin{eqnarray}
\label{eq:pc1}
-\frac{n_y}{\mu}(B_z^{**} - B_z^L) &=& +c (D_x^{**} - D_x^R) \\
\label{eq:pc2}
+\frac{n_x}{\mu}(B_z^{**} - B_z^L) &=& +c (D_y^{**} - D_y^R) \\
\label{eq:pc3}
\frac{1}{\varepsilon}[ (D_y^{**} n_x - D_x^{**} n_y) - (D_y^R n_x - D_x^R n_y)] &=& +c (B_z^{**} - B_z^R)
\end{eqnarray}
From (\ref{eq:mc1}), (\ref{eq:mc2}) we get $\D^* \cdot \un = \D^L \cdot \un$ while (\ref{eq:pc1}), (\ref{eq:pc2}) yields $\D^{**}\cdot\un = \D^R \cdot\un$, and hence the normal component of $\D$ is continuous throughout the Riemann fan. From (\ref{eq:z1}), (\ref{eq:z2}) we get $B_z^* = B_z^{**}$ while (\ref{eq:z3}) shows that $(D_y^{**} n_x - D_x^{**} n_y) = (D_y^* n_x - D_x^* n_y)$, which implies that $\con^* = \con^{**}$. Adding (\ref{eq:mc3}) and (\ref{eq:pc3}) we get
\[
D_y^* n_x - D_x^* n_y = \frac{(D_y^L n_x - D_x^Ln_y) + (D_y^R n_x - D_x^R n_y)}{2} - \frac{\varepsilon c}{2}(B_z^R - B_z^L)
\]
By combining the equations in the form $(\ref{eq:pc2}) n_x - (\ref{eq:pc1}) n_y$ we get
\[
B_z^* = \frac{1}{2}(B_z^L + B_z^R) - \frac{\mu c}{2}[ (D_y^R n_x - D_x^R n_y) - (D_y^L n_x - D_x^Ln_y)]
\]
The numerical flux is given by
\[
\hat{\fx}_n  = \begin{bmatrix}
-\frac{n_y}{\mu} B_z^* \\
+\frac{n_x}{\mu} B_z^* \\
\frac{1}{\varepsilon}(D_y^* n_x - D_x^* n_y) \end{bmatrix}
\]
From the above solution, we can extract the information required in the FR and DG schemes. For a vertical face $(n_x,n_y) = (1,0)$ and the numerical fluxes are
\[
\hHz = \frac{1}{\mu} B_z^*, \qquad \hEy = \frac{1}{\varepsilon} D_y^*
\]
and on a horizontal face $(n_x, n_y) = (0,1)$
\[
\hHz = \frac{1}{\mu} B_z^*, \qquad \hEx = \frac{1}{\varepsilon} D_x^*
\]

\subsection{Solution of the 2-D Riemann problem}
\label{sec:riem2d}
The update of the normal component of $\D$ stored on the faces requires the knowledge of the magnetic field $H_z$ at the vertices. We need a unique value of this quantity in order to obtain a constraint perserving scheme. On a 2-D Cartesian mesh, at each vertex of the mesh, we have four states that come together to define a 2-D Riemann problem. The solution of this Riemann problem leads to the generation of four 1-D wave structures and a strongly interacting state around the vertex. More details on the solution of the 2-D Riemann problem can be found in~\cite{Balsara2017}. For the sake of completeness, we summarize the formulae for the specific case of scalar material properties. The magnetic field at the vertex is given by
\[
\tHz = \frac{1}{\mu} B_z^*
\]
where
\[
\begin{aligned}
B_z^* = \frac{1}{4}(B_z^{DL} + B_z^{DR} + B_z^{UL} + B_z^{UR}) + & \frac{\mu c}{2}\left[ \half(D_x^{UR} + D_x^{UL}) - \half (D_x^{DR} + D_x^{DL}) \right] \\
- & \frac{\mu c}{2} \left[ \half(D_y^{UR} + D_y^{DR}) - \half (D_y^{UL} + D_y^{DL}) \right]
\end{aligned}
\]
and the superscripts DL, DR, UL, UR refer to the four states meeting at the vertex. Note that due to the divergence conforming nature of our approximating polynomials, we actually have $D_x^{DL} = D_x^{DR}$, $D_x^{UL} = D_x^{UR}$, $D_y^{DL} = D_y^{UL}$ and $D_y^{DR} = D_y^{UR}$. The flux $\tHz$ has the usual structure of a central flux which is the average of the four fluxes at the vertex and some jump terms. As shown in section~(\ref{sec:stab}), the jump terms add some dissipation to the numerical scheme and hence are important to obtain a stable scheme.
\bibliographystyle{siam}
\bibliography{bbibtex}
\end{document}